\documentclass{amsart}
\usepackage[utf8]{inputenc}

\usepackage{textcmds}
\usepackage{mathtools}
\usepackage{amssymb}
\usepackage{amsthm}
\usepackage{amsmath}
\usepackage{bm}
\usepackage{bbm}
\usepackage{enumerate}
\usepackage{cite}
\usepackage{xfrac}
\usepackage{xcolor}
\usepackage{tikz-cd}
\usepackage{lipsum}
\usepackage{adjustbox}
\usepackage{url}
\usepackage{hyperref}
\usepackage{multicol}

\usepackage[left=3cm,right=3cm,top=3cm,bottom=4.5cm]{geometry}

\theoremstyle{plain}
\newtheorem{theorem}{Theorem}[section]
\newtheorem{corollary}[theorem]{Corollary}
\newtheorem{lemma}[theorem]{Lemma}

\theoremstyle{definition}
\newtheorem{definition}[theorem]{Definition}
\newtheorem*{definition*}{Definition}

\theoremstyle{plain}

\theoremstyle{plain}

\theoremstyle{plain}

\theoremstyle{plain}
\newtheorem{fact}[theorem]{Fact}

\theoremstyle{remark}
\newtheorem{remark}[theorem]{Remark}

\theoremstyle{remark}
\newtheorem{notation}[theorem]{Notation}

\theoremstyle{remark}
\newtheorem{example}[theorem]{Example}

\theoremstyle{remark}
\newtheorem{convention}[theorem]{Convention}

\theoremstyle{plain}
\newtheorem{proposition}[theorem]{Proposition}

\theoremstyle{plain}

\newtheorem{theoremintro}{Theorem}

\newtheorem{corollaryintro}{Corollary}

\theoremstyle{definition}

\DeclareMathOperator{\projection}{proj}

\def\dotminussym#1{%
  \setbox0=\hbox{$-$}%
  \kern.5\wd0%
  \hbox to 0pt{\hss\hbox{$-$}\hss}%
  \raise.6\ht0\hbox to 0pt{\hss$.$\hss}%
  \kern.5\wd0%
}
\newcommand{\dotminus}{\mathbin{\mathpalette\dotminussym{}}}

\def\Ind#1#2{#1\setbox0=\hbox{$#1x$}\kern\wd0\hbox to 0pt{\hss$#1\mid$\hss}
\lower.9\ht0\hbox to 0pt{\hss$#1\smile$\hss}\kern\wd0}

\def\ind{\mathop{\mathpalette\Ind{}}}

\def\notind#1#2{#1\setbox0=\hbox{$#1x$}\kern\wd0
\hbox to 0pt{\mathchardef\nn=12854\hss$#1\nn$\kern1.4\wd0\hss}
\hbox to 0pt{\hss$#1\mid$\hss}\lower.9\ht0 \hbox to 0pt{\hss$#1\smile$\hss}\kern\wd0}

\title{Model theory of difference fields with an additive character on the fixed field}

\author{Stefan Marian Ludwig}
\address{Albert-Ludwigs-Universität Freiburg,
Mathematisches Institut,
Abteilung für Mathematische Logik,
Ernst-Zermelo-Straße 1,
79104 Freiburg i. B., Germany}
\email{stefan.ludwig@mathematik.uni-freiburg.de}	
\thanks{SML has received funding from the European Union's Horizon 2020 research and innovation programme under the Marie Sk\l{}odowska-Curie grant agreement N\textsuperscript{\underline{o}} 945322.  \includegraphics[scale=0.025]{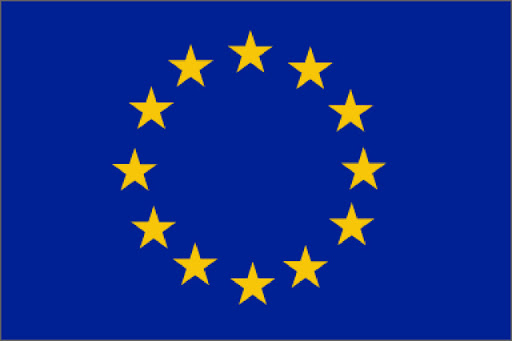}	Moreover, he was partially supported by GeoMod AAPG2019 (ANR-DFG), `Geometric and Combinatorial Configurations in Model Theory'.}



\begin{document}

\begin{abstract} 
Following a research line proposed by Hrushovski in his work on pseudofinite fields with an additive character, we investigate the theory $\mathrm{ACFA}^{+}$ which is the model companion of the theory of difference fields with an additive character on the fixed field added as a continuous logic predicate. $\mathrm{ACFA}^{+}$ is the common theory (in characteristic $0$) of the algebraic closure of finite fields with the Frobenius automorphism and the standard character on the fixed field and turns out to be a simple theory. We fully characterise 3-amalgamation and deduce that the connected component of the Kim-Pillay group (for any completion of $\mathrm{ACFA}^{+}$) is abelian as conjectured by Hrushovski. Finally, we describe a natural expansion of $\mathrm{ACFA}^{+}$ in which geometric elimination of continuous logic imaginaries holds.
\end{abstract}

\maketitle




\section{Introduction}

The common theory of all finite fields in the language of rings was determined and studied by Ax in his seminal work \cite{Hrushovski2021AxsTW}. The infinite models of this theory are called pseudofinite fields and are characterised by being perfect, PAC and having Galois group $\hat{\mathbb{Z}}$. The essential number-theoretic ingredient in Ax's work are the Čebotarev density theorem and the Lang-Weil bounds on the number of rational points of absolutely irreducible varieties over finite fields. They were famously generalised to difference varieties (with the automorphism given by the Frobenius) by Hrushovski in \cite{elementaryfrobenius} with a later purely algebro-geometric proof due to Shuddhodan and Varshavsky \cite{HrushLangWeilbyShuddhodanVayarslan}. As a consequence Hrushovski generalised Ax's result. He showed that a difference field $(K,\sigma)$ (with infinite fixed field $F:=\mathrm{Fix}(\sigma)$) is a model of the common theory of the algebraic closure of finite fields with the Frobenius if and only if $(K,\sigma)\models\mathrm{ACFA}$. Here, $\mathrm{ACFA}$ is the model companion of the theory of difference fields as introduced in \cite{Macintyre1997-MACGAO} and \cite{acfa}. The in-depth study of $\mathrm{ACFA}$ by Chatzidakis and Hrushovski in \cite{acfa} was continued in \cite{differencefields2} and then applied to questions in algebraic dynamics \cite{Chatzidakis_Hrushovski_2008-descent1}, \cite{chatzhrush_descent2},\cite{scanlonmedvedev}.
See also Hrushovski's work on the Manin-Mumford conjecture \cite{HRUSHOVSKIManinMumford}.\\
More recently, motivated by work of Kowalski on exponential sums over definable sets in finite fields \cite{Kowalski2005ExponentialSO}, Hrushovski generalised Ax's result in yet another direction determining the common theory $\mathrm{PF}^{+}$ of finite fields with a non-trivial additive character added as a continuous logic predicate \cite{Hrushovski2021AxsTW}. In \cite{ludwig2025pseudofinitefieldsadditivemultiplicative} the author of this article extended this work to finite fields with both suffgeneric additive and multiplicative character. This generalisation is partly motivated by the fact that the main number-theoretic ingredient used in the context of $\mathrm{PF}^{+}$, the so-called Weil bounds on exponential sums over algebraic curves, also hold for character sums involving a multiplicative character such as the ubiquitous Gauss sums.
Hrushovski does not only give a natural axiomatisation of $\mathrm{PF}^{+}$ but moreover shows quantifier elimination in a natural definitional expansion of the language. Moreover, $\mathrm{PF}^{+}$ turns out to be a simple theory, as $n$-amalgamation (and thus in particular 3-amalgamation) holds over all algebraically closed sets. In addition a generalisation of the definability of the Chatzidakis-van den Dries-Macintyre counting measure to this context is obtained.
In this article we consider the combination of both of Hrushovski's generalisations to $\mathrm{ACFA}$ and $\mathrm{PF}^{+}$ of Ax's work to thus treat the following questions. \textit{What is the common theory of the algebraic closure of finite fields with the Frobenius and an additive character on the fixed field? What are its model-theoretic properties? In particular, which model-theoretic phenomena occur due to the interaction of the character with the Frobenius?}\\
It was already suggested by Hrushovski in \cite[Section 6.3]{Hrushovski2021AxsTW} to study the theory $\mathrm{ACFA}^{+}$, the amalgam of $\mathrm{ACFA}$ and $\mathrm{PF}^{+}$ and here we will follow this proposed line of research. Hrushovski already observes that 3-amalgamation over algebraically closed sets $A$ in $\mathrm{ACFA}^{+}$ may fail. He conjectures that from a complete understanding of it one should be able to deduce that the identity component $H_{\mathrm{KP}}$ of the Kim-Pillay group is abelian (but possibly nontrivial). We will give a complete description of $H_{\mathrm{KP}}$ and thus in particular confirm the conjectured abelianity. The main step will indeed be to give a full characterisation of 3-amalgamation which we will then also use to classify the CL-imaginaries (hyperimaginaries) present in $\mathrm{ACFA}^{+}$ as well as to deduce that $\mathrm{ACFA}^{+}$ is simple.
The occurrence of a  non-trivial $H_{\mathrm{KP}}$ in a natural theory such as $\mathrm{ACFA}^{+}$ constitutes an interesting phenomenon from a model-theoretic point of view.
In simple theories the full Lascar group is already given by the Kim-Pillay group \cite{noteonlstp},\cite{simplicityCATs}. It is a longstanding open question in classical discrete logic whether a simple theory exists whose Kim-Pillay-group has a non-trivial connected component. While in more restrictive cases in classical logic (e.g. \cite{Buechler1999-BUELST}) it can be shown that such a group can not exist, Hrushovski already observes in \cite{Hrushovski1998SimplicityAT} that in the context of Robinson theories and thus outside of the classical first-order context, such a group may be found. $\mathrm{ACFA}^{+}$ will be a natural example in continuous logic.\\
Another motivation to study expansions of $\mathrm{ACFA}$ involving characters stems from the following. In \cite{katz} Katz raised the question of a possible theory of exponential sums over $\mathbb{Z}$. Kowalski suggested (Remark 17 in \cite{Kowalski2005ExponentialSO}) that model theory might be helpful and Hrushovski points out (6.2 in \cite{Hrushovski2021AxsTW}) that $\mathrm{ACFA}$ and extensions thereof might yield a natural context. He also remarks that \emph{because of the probable abelian nature of the connected components of $\mathrm{ACFA}^{+}$ ($\dots$) it seems clear that $\mathrm{ACFA}^{+}$ is not the full answer.}\footnote{6.3 (iii) in \cite{Hrushovski2021AxsTW}} So, in some sense our result on $H_{\mathrm{KP}}$ also showcases the limitations of $\mathrm{ACFA}^{+}$. However, while not making any direct progress in this direction, we hope that our work on $\mathrm{ACFA}^{+}$ can serve as a preparation of grounds for some future developments in this context.\\
Finally, while (the failure of) 3-amalgamation that we treat in this article already gives an interesting model-theoretic account of a non-trivial interaction of the field automorphism and the character on the fixed field, we will show in a forthcoming article (based on Chapter 4 of the author's PhD thesis) that higher amalgamation in
 $\mathrm{ACFA}^{+}$ behaves rather surprisingly (at least to the author). On one hand, one can construct a set satisfying the criterion for 3-amalgamation but for which  4-amalgamation does not hold. On the other hand, using a technical stability-theoretic argument one can show that $n$-amalgamation still holds for all $n\in\mathbb{N}$ over substructures that are models of $\mathrm{ACFA}$.
\subsection*{Presentation of main results}
We will work in the language $\mathcal{L}_{\sigma}^{+}$ consisting of $\mathcal{L}_{\mathrm{ring}}$ together with a unary function symbol $\sigma$ and a continuous logic predicate $\Psi$ with range $S^{1}\cup\{0\}\subseteq\mathbb{C}$. We write $F$ for the fixed field $\mathrm{Fix}(\sigma)$ of a model $K\models \mathrm{ACFA}$. From now on we restrict to work in characteristic $0$, i.e., $\mathrm{ACFA}$ and $\mathrm{PF}$ will be assumed to contain a set of axioms specifying that the characteristic is $0$. 

\begin{definition*}[See Definition \ref{definitionacfaplus}]
    An $\mathcal{L}_{\sigma}^{+}$-structure $(K,\sigma,\Psi)$ is a model of $\mathrm{ACFA}^{+}$, if
    \begin{itemize}
        \item $(K,\sigma)\models \mathrm{ACFA}$,
        \item $\Psi(K\backslash F)=0$ and $\Psi\restriction_{F}:(F,+)\rightarrow(S^{1},\cdot)$ is a group homomorphism,
        \item For every absolutely irreducible curve $C\subset \mathbb{A}^{n}$ that is defined over $F$ and not contained in a rational hyperplane over $F$, the set $\Psi^{(n)}(C(F))\subseteq\mathbb{T}$ is a dense subset.
    \end{itemize}
\end{definition*}
The last axiom is expressible by a set of first-order sentences (see Fact \ref{factaxiomisexpressible}) and is already contained in the theory $\mathrm{PF}^{+}$. In fact, the models of $\mathrm{ACFA}^{+}$ are exactly the models of $\mathrm{ACFA}$
with the enriched structure on the fixed field being a model of $\mathrm{PF}^{+}$.
Regarding this axiomatisation, one might wonder whether all model-theoretic results on $\mathrm{ACFA}^{+}$ immediately follow from what is known from $\mathrm{ACFA}$ and $\mathrm{PF}^{+}$. We will later see that this is not the case, in particular, when type-amalgamation comes into play. However, a quantifier elimination result for $\mathrm{ACFA}^{+}$ can be obtained from the corresponding results for $\mathrm{ACFA}$ and $\mathrm{PF}^{+}$ by using the fact that $F$ is stably embedded in any model of $\mathrm{ACFA}^{+}$.
In \hyperref[sectionstructureexpansionsonstablyembeded]{Appendix II} we verify that known statements from classical logic on structure expansions on a definable set work similarly when we allow the expansion to contain continuous logic predicates. Then, apart from taking care of the below described subtlety, we obtain the following rather straightforwardly. 

\begin{theoremintro}[See Corollary \ref{characterisationeltequivacfaplus},
Theorem \ref{modelcompanion} and
Lemma \ref{inducedstructureplusfixedfield}]
   The $\mathcal{L}_{\sigma}^{+}$-theory $\mathrm{ACFA}^{+}$ satisfies the following properties. 
\begin{itemize}
     \item Every difference field $(K,\sigma)$ with additive character on the fixed field embeds in a model of $\mathrm{ACFA}^{+}$. Moreover $\mathrm{ACFA}^{+}$ is model-complete.
     \item Let $(K_{1},\sigma_{1},\Psi_{1}),(K_{2},\sigma_{2},\Psi_{2})$ be two models of $\mathrm{ACFA}^{+}$
 and $E$ a common algebraically closed difference subfield on which $\sigma_{1}$ and $\sigma_{2}$ agree as automorphisms and such that $\Psi_{1}, \Psi_{2}$ agree on $\mathrm{Fix}(E)$, then
 \[(K_{1},\sigma_{1},\Psi_{1})\equiv_{E} (K_{2},\sigma_{2},\Psi_{2}).
 \]
    \item $(F,\Psi)$ is stably embedded. 
 \end{itemize}
\end{theoremintro}
The small but important subtlety that we have to take care of is that the induced structure on $F$ (in $\mathrm{ACFA}$) is slightly larger than the pure field structure. More precisely, the induced structure distinguishes a generator of the Galois group (the restriction of $\sigma$). In other words, a sequence of algebraic imaginaries is named. To be able to apply the above reasoning, we therefore have to slightly adjust the corresponding results in $\mathrm{PF}^{+}$ (see Lemma \ref{lemmainducedstructuredefinableclosure} and Corollary \ref{corollaryqepfplusind}).\\
We denote by $\Psi_{q}$ the \textit{standard additive character} on $\mathbb{F}_{q}$, that is, $\Psi_{q}=\mathrm{exp}(2\pi i \mathrm{Tr}(x))$ where $\mathrm{Tr}:\mathbb{F}_{q}\rightarrow\mathbb{F}_{p}$ is the trace to the corresponding prime field. 
Combining the corresponding results of Hrushovski on $\mathrm{ACFA}$ \cite{elementaryfrobenius} and $\mathrm{PF}^{+}$ \cite{Hrushovski2021AxsTW} we instantly obtain that any characteristic 0 ultraproduct of the structures $(\overline{\mathbb{F}}_{q},\mathrm{Frob}_{q},\Psi_{q})$ is a model of $\mathrm{ACFA}^{+}$.
We also prove the other direction, that is, for any completion of $\mathrm{ACFA}^{+}$ there is a model consisting of an ultraproduct of the above structures.
Note, however, that if we restrict to work with the class of prime fields and standard character this hinges on an open number-theoretic question (as already in $\mathrm{PF}^{+}$, Section 5 in \cite{Hrushovski2021AxsTW}). We thus obtain the following.
\begin{theoremintro}[See Theorem \ref{maintheoremlimittheory}]
     $\mathrm{ACFA}^{+}$ is the characteristic $0$-asymptotic theory of the algebraic closure of finite fields with Frobenius and standard character and similarly for the algebraic closure of prime fields with Frobenius and arbitrary character.  
\end{theoremintro}
We use the same proof as Hrushovski for $\mathrm{PF}^{+}$ but again the above-mentioned subtlety on the choice of the generator of the Galois group (in this case of $F\cap\overline{\mathbb{Q}}$) has to be taken into account. It can be remedied using a basic number-theoretic calculation (Lemma \ref{primecomputationlemma}).\\
We denote by $\mathrm{acl}_{\sigma}(A)$ the model-theoretic algebraic closure of some $A\subseteq\mathcal{M}\models \mathrm{ACFA}$ (which is the field-theoretic algebraic closure of the difference field generated by $A$). Independence in $\mathrm{ACFA}^{+}$ will be defined as in $\mathrm{ACFA}$, that is, two sets $A,B$ are independent over a common subset $C$ if $\mathrm{acl}_{\sigma}(A)$ is algebraically independent from $\mathrm{acl}_{\sigma}(B)$ over $\mathrm{acl}_{\sigma}(C)$. The independence relation then inherits its good properties from $\mathrm{ACFA}$, except, possibly, the Independence Theorem/3-amalgamation over algebraically closed sets. We then determine the exact condition on $A$ for 3-amalgamation to hold over $A$.
\begin{theoremintro}[See Theorem \ref{theorem3amalgamationcharacterisation}]\label{introductiontheorem3amalgamation}
    In $\mathrm{ACFA}^{+}$ 3-amalgamation holds over $A=\mathrm{acl}_{\sigma}(A)$ if and only if $A$ is $\sigma$-AS-closed, that is, for all $a\in A$, there is $b\in A$ such that $\sigma(b)-b=a$.
\end{theoremintro}
To see that being $\sigma$-AS-closed is indeed required for 3-amalgamation, one first observes that $\mathfrak{T}_{a}:=\{x\,|\,\sigma(x)-x=a\}$ defines an $(F,+)$-Torsor. Then, one can find $\alpha_{1},\alpha_{2},\alpha_{3}$ in $\mathfrak{T}_{a}$, pairwise independent over $A$. Whenever $\mathfrak{T}_{a}\cap A=\emptyset$, for any $r\in S^{1}$, it is consistent with $tp(\alpha_{1}/A)\cup tp(\alpha_{2}/A)$ that $\Psi(\alpha_{1}-\alpha_{2})=r$. Using this, one can then construct a 3-amalgamation problem over $A$ which does not have a solution. For the other direction, one shows that any obstruction to 3-amalgamation can only come from additive constraints on elements of $F$ which, by a stability-theoretic argument, are shown to originate from elements of $\mathfrak{T}_{a}$ with $\mathfrak{T}_{a}\cap A=\emptyset$.
An important consequence we can then directly deduce from Theorem \ref{introductiontheorem3amalgamation} is the following.
\setcounter{corollaryintro}{3}
\addtocounter{theoremintro}{1}
\begin{corollaryintro}[See Corollary \ref{corollaryacfaplusissimple}]
    The theory $\mathrm{ACFA}^{+}$ is simple.
\end{corollaryintro}
Recall that the Kim-Pillay group is the quotient $\mathrm{Aut}(\mathcal{M})/\mathrm{Aut}_{\mathrm{KP}}(\mathcal{M})$ where $\mathcal{M}$ is a monster model of the ambient theory and $\mathrm{Aut}_{\mathrm{KP}}(\mathcal{M})$ denotes those automorphisms that fix all classes of all bounded type-definable equivalence relations. The Kim-Pillay group is an invariant of the theory and can be turned into a topological group via the so-called logic topology. Its identity component is isomorphic to the quotient $\mathrm{Gal}_{\mathrm{KP}}/\mathrm{Gal}_{\mathrm{Sh}}$ where the latter is the subgroup of automorphisms that stabilise all classes of definable equivalence relations with finitely many classes. In \hyperref[sectionmodeltheoreticgaloisgroups]{Appendix I} we give an account of how those facts work out in the realm of continuous logic. We write $\mathrm{Gal}_{\mathrm{KP}}^{A}$ if we work over a set of parameters $A$ and denote by $H_{\mathrm{KP}}^{A}$ its identity component.
In (any completion of) the theory $\mathrm{ACFA}^{+}$ the topological group $H_{\mathrm{KP}}$ can be determined by its action on the set of hyperimaginaries (CL-imaginaries) $\mathfrak{T}_{a}/E_{a}\cong S^{1}$ where $E_{a}$ is given by the equivalence relation defined by $E_{a}(x,y)\iff \Psi(x-y)=1$. The group $H_{\mathrm{KP}}$ acts on those objects via rotations. The action turns out to be transitive whenever $\mathfrak{T}_{a}\cap A=\emptyset$ and this can be generalised to products. Altogether we obtain the following.
\begin{theoremintro}[See Theorem \ref{theoremidentycompnenKPgrouptopequiv}]
    The identity component $H_{\mathrm{KP}}^{A}$ is isomorphic as a topological group to $\prod_{B}\varprojlim \mathbb{R}/n\mathbb{Z}$  where $B$ is a basis of a $\mathbb{Q}$-vector space complement of $\wp_{\sigma}(A):=\{\sigma(a)-a\,|\,a\in A\}$. In particular, it is abelian.
\end{theoremintro}
Hyperimaginaries are interdefinable with (sequences of) continuous logic imaginaries (referred to as CL-imaginaries) which allows us to work inside the continuous analogue of $T^{\mathrm{eq}}$. Now, another consequence that one can obtain from Theorem \ref{introductiontheorem3amalgamation} is a (weak) classification of imaginaries. Concretely, we obtain the following.
\begin{theoremintro}[See Theorem \ref{theoremweakEIACFAplusEI}]
    The theory $\mathrm{ACFA}^{+}_{\mathrm{EI}}$ which one obtains from $\mathrm{ACFA}^{+}$ by adding the CL-imaginaries $\mathfrak{T}_{a}/E_{a}$ (see Definition \ref{definitionEIextension}) has weak elimination of imaginaries.
\end{theoremintro}
The proof is similar to the one for $\mathrm{ACFA}$ in \cite{acfa} with Theorem \ref{introductiontheorem3amalgamation} as the main ingredient
but in our context we use a metric version of Neumann's Lemma \cite{Conant2022SeparationFI}. Unlike in $\mathrm{ACFA}$, due to the imaginaries we added, we cannot obtain strong elimination of imaginaries. We can, however, slightly strengthen the above result and eliminate imaginaries up to finite sets of finite tuples instead of compact ones (which weak elimination of imaginaries in continuous logic would only yield). This is Corollary \ref{corollaryEIuptofinitesetsoffinitetuples}.

\subsection*{Structure of the article}After reviewing some properties of the fixed field in a model of $\mathrm{ACFA}$ in Section \ref{sectionpreliminaries} we introduce the theory $\mathrm{ACFA}^{+}$ in Section \ref{sectionintroacfapluschaptertwo} and prove a quantifier elimination result as well as its characterisation as the limit theory of the algebraic closure of finite fields equipped with the Frobenius and an additive character (on the fixed field). In Section \ref{section3amalgam} we characterise 3-amalgamation and deduce that $\mathrm{ACFA}^{+}$ is simple. Then, we use these results in Section \ref{sectionautomgrouptorsors} to give a concrete description of the connected component of the Kim-Pillay group. Finally, we characterise the CL-imaginaries in Section \ref{sectioneliminationofimaginaries} and conclude the main part of the article with some remarks on forthcoming work in Section \ref{sectionfurtherremarks}. In \hyperref[sectionmodeltheoreticgaloisgroups]{Appendix I} we give an account of the Kim-Pillay group in continuous logic presenting the continuous version of several well-known results. In \hyperref[sectionstructureexpansionsonstablyembeded]{Appendix II} we present results on continuous logic expansions on stably embedded sets in a classical logic structures that are used in Section \ref{sectionintroacfapluschaptertwo}.

\subsection*{Acknowledgments}
This article is based on Chapters 2 and 3 (and the appendix on Chapter 1) of my PhD thesis \cite{ludwig:tel-05236078} written at the École Normale Supérieure Paris. I would like to thank my supervisors Zoé Chatzidakis and Martin Hils for their invaluable support during my thesis as well as Ehud Hrushovski for discussing the topic with me during a secondment I did at the University of Oxford during my PhD thesis. Further, I would like to thank Silvain Rideau-Kikuchi for pointing out the question on elimination of imaginaries (Section \ref{sectioneliminationofimaginaries}) and Mira Tartarotti for helpful comments on the topic of Appendix I.

\section{Preliminaries}\label{sectionpreliminaries}

\subsection{Continuous logic}
In this work we will not have to make use of the full generality of continuous logic, but we will restrict ourselves to working in the setting where the metric is discrete (except of Section \ref{sectioneliminationofimaginaries}). In this case we can use the equality symbol as we are used to from classical logic\footnote{which we consider as a $\{0,1\}$-valued metric if we want to employ the classical continuous logic formalism.} but predicates are allowed to take complex values. This can be easily coded in the real-valued continuous logic formalism as in \cite{mtfms} by using predicates for the real as well as for the imaginary part. See the introduction of \cite{Hrushovski2021AxsTW} or 2.3 in \cite{ludwig2025pseudofinitefieldsadditivemultiplicative} for an account thereof.\\
As one of the main results of this article is a confirmation of a conjecture of Hrushovski on the Kim-Pillay group of (completions of) $\mathrm{ACFA}^{+}$ we will briefly recall its definition. The Kim-Pillay group comes with a natural topology. In \hyperref[sectionmodeltheoreticgaloisgroups]{Appendix I} we give an account thereof in the context of continuous logic and transfer several results from the classical context.

\begin{definition}
An equivalence relation on $\mathcal{M}$ is called \textit{bounded} if it has at most $2^{|T|}$ classes. 
\end{definition}

\begin{definition}
We say that two (possibly infinite) tuples $\bar{a},\bar{b}$ have the same \textit{Kim-Pillay strong type}, or simply \textit{KP-strong type} over some small set $A$, written $\bar{a}\equiv_{A}^{\mathrm{KP}}\bar{b}$, if they are equivalent in any bounded type-definable over $A$ equivalence relation on the sort of $\bar{a},\bar{b}$. (We omit $A$ in the case of $A=\emptyset$).\\
We define the group $\mathrm{Aut}_{\mathrm{KP}}(M)$ as the subgroup of the automorphism group consisting of those $f\in \mathrm{Aut}(M)$ fixing the class of any bounded $\emptyset$-type-definable equivalence relation. In other words, $f\in \mathrm{Aut}(M)$ lies in $\mathrm{Aut}_{\mathrm{KP}}(M)$ if and only if $f(\bar{a})\equiv^{\mathrm{KP}}\bar{a}$ for any tuple $\bar{a}$ in $M$. Similarly, we define the group $\mathrm{Aut}_{\mathrm{KP}}^{A}(M)$ for bounded $A$-type-definable equivalence relations.
The \textit{Kim-Pillay Galois group} is then the group $\mathrm{Aut}(M)/\mathrm{Aut}_{\mathrm{KP}}(M)$ and will be denoted by $\mathrm{Gal}_{\mathrm{KP}}(M)$. Correspondingly, we define the group $\mathrm{Gal}_{\mathrm{KP}}^{A}(M)$ as $\mathrm{Aut}(M/A)/\mathrm{Aut}_{\mathrm{KP}}^{A}(M)$.
\end{definition}

\subsection{Some facts around $\mathrm{ACFA}$ and $\mathrm{PF}$}\label{sectionfactsfromclassicalcontextchaptertwo}
In this section we recall some of the results around the model theory of pseudofinite fields and difference fields. Pseudofinite fields, the infinite models of the theory of all finite fields, were first introduced and studied by Ax in \cite{ax-elttheoryoffinitefields}. In \cite{ChatzidakisMacvdd} the theory was further developed and the authors obtained a generalisation of the Lang-Weil to definable sets as well as the definability of the (nonstandard) counting measure. To this day, the theory of pseudofinite fields has been actively studied in model theory and beyond, yielding many applications (see, for example, \cite{HrushovskiPillay1995} or \cite{tao-expandingpolynomials}).

\begin{notation}
    Given a field $K$ we will use interchangeably (for the rest of the article) the notations $\overline{K}$ and $K^{\mathrm{alg}}$ for the (field theoretic) algebraic closure. The model theoretic algebraic closure will be denoted by $\mathrm{acl}(K)$ or $\mathrm{acl}_{T}(K)$, $\mathrm{acl}_{\mathcal{L}}(K)$ when the ambient theory or language is not clear from the context.
\end{notation}

\begin{definition}
    A field $F$ is called \textit{pseudofinite}, if it satisfies the following axioms.
    \begin{itemize}
        \item $F$ is perfect.
        \item $F$ is PAC, that is, every absolutely irreducible variety over $F$ has an $F$-rational point.
        \item $\mathrm{Gal}(F)=\hat{\mathbb{Z}}$.
    \end{itemize}
    The $\mathcal{L}_{\mathrm{ring}}$-theory that expresses these axioms will be denoted by $\mathrm{PF}$.
\end{definition}

\begin{fact}\label{theoremaxtheorempsf}(Theorem 9 in \cite{ax-elttheoryoffinitefields})
    A field $F$ is pseudofinite if and only if it is infinite and satisfies all $\mathcal{L}_{\mathrm{ring}}$-sentences true in all finite fields.
\end{fact}

\begin{fact}\label{kiefeQE}(See Theorem 2 in \cite{kiefedefbsetszetafunction}.) Every $\mathcal{L}_{\mathrm{ring}}$-formula $\phi(\bar{x})$ is equivalent modulo $\mathrm{PF}$ to a boolean combination of formulas of the form $\exists t\,g(\bar{x},t)=0$ where $g\in\mathbb{Z}[\bar{X},T]$.
    
\end{fact}
Now we give some facts around the model theory of difference fields. For us a difference field will be a pair $(K,\sigma)$ consisting of a field together with an automorphism $\sigma:K\rightarrow K$. We recall some properties of the theory $\mathrm{ACFA}$ of existentially closed difference fields. It was first studied in \cite{Macintyre1997-MACGAO} and was then further developed by Chatzidakis and Hrushovski in \cite{acfa} and subsequent papers (e.g.  \cite{chatzhrush_descent2},\cite{Chatzidakis_Hrushovski_2008-descent1}) which also contain several applications to algebraic dynamics. We work in the language $\mathcal{L}_{\sigma}$ that consists of the ring language together with a unary function symbol $\sigma$.

\begin{definition}(1.1 in \cite{acfa})
    Let $\mathrm{ACFA}$ be the theory axiomatised by the scheme of axioms expressing
the following properties of the $\mathcal{L}_{\sigma}$-structure $(K,\sigma)$:
\begin{itemize}
    \item $\sigma$ is an automorphism of $K$.
    \item $K$ is an algebraically closed field.
    \item For every absolutely irreducible varieties $V,U$ with $V\subseteq U\times \sigma(U)$ projecting generically onto $U$ and $\sigma(U)$ and every proper algebraic subset $W\subset V$ there is $a\in U(K)$ such that $(a,\sigma(a))\in V\backslash W$.
\end{itemize}

\end{definition}

\begin{fact}\label{basicfactsonacfa}(1.1-1.3 in \cite{acfa})
\begin{itemize}
    \item Every difference field embeds in a model of $\mathrm{ACFA}$ and $\mathrm{ACFA}$ is model-complete.
    \item Let $(K_{1},\sigma_{1}),(K_{2},\sigma_{2})$ be models of $\mathrm{ACFA}$ and $E$ a common difference subfield. Then \[(K_{1},\sigma_{1})\equiv_{E}(K_{2},\sigma_{2})\;\iff\; \left(E^{\mathrm{alg}},\sigma_{1}\restriction_{E^{\mathrm{alg}}}\right)\cong\left(E^{\mathrm{alg}},\sigma_{2}\restriction_{E^{\mathrm{alg}}}\right).\]
    \item The fixed field $F$ of a model $(K,\sigma)\models\mathrm{ACFA}$ is a pseudofinite field. 
\end{itemize}
    
\end{fact}

\begin{corollary}(1.6 in \cite{acfa})\label{corollaryexplicitqeforacfa}
Every formula $\phi(\bar{x})$ is equivalent modulo $\mathrm{ACFA}$ to a boolean combination of formula of the form $\exists y\,\psi(\bar{x},y)$ for $y$ a single variable, $\psi(\bar{x},y)$ quantifier-free, and $\psi(\bar{a},b)$ implying for every $\bar{a},b$ that b is algebraic (in the pure field sense) over (the field generated by) $\{\bar{a},\sigma(\bar{a}),\dots,\sigma^{m}(\bar{a})\}$ for some $m\in\mathbb{N}$.
    
\end{corollary}

\subsection{The fixed field in a model of $\mathrm{ACFA}$}\label{sectionfixedfieldsclassical}

We recall some more facts and notions concerning the fixed field $F$ of a model $K\models\mathrm{ACFA}$. The references are (1.11) and (1.13) in \cite{acfa}.

\begin{fact}\label{factstableembeddednessacfaclassicallogic}(See (1.11) in \cite{acfa}.)
The fixed field $F$ is stably embedded. We even have that the restriction of every $\mathcal{L}_{\sigma}(K)$-definable subset to $F$ is definable with parameters from $F$ only using the language $\mathcal{L}_{\mathrm{ring}}$.
    
\end{fact}

Next, we want to describe the induced structure on the fixed field of a model of $\mathrm{ACFA}$. Whereas the fixed field $F$ is stably embedded, the induced structure on $F$ is not exactly the ring structure, but something slightly bigger. The reason is that if $f_{\bar{a}}(X):=X^{n}+a_{1}X^{n-1}+\dots+a_{n}$ is irreducible over $F$ for $\bar{a}\in F^{n}$ and $\alpha$ a root of $f_{\bar{a}}$, then $\sigma(\alpha)=\alpha^{n}+b_{1}\alpha^{n-1}+\dots+b_{n}$ for some $\bar{b}\in F^{n}$ and the datum of the tuple $\bar{b}$ does not depend on the choice of the root $\alpha$ (using that the Galois group of $F$ is abelian). However, this will be the only thing we have to add.

\begin{definition}\label{inducedstructurelanguage}
    Let $\mathcal{L}_{\sigma-\mathrm{ind}}$ be the language $\mathcal{L}_{\mathrm{ring}}$ enriched by a $2n$-ary predicate $C_{n}$ for each $n\in\mathbb{N}$. As introduced in (1.13) of \cite{acfa} we consider the $\mathcal{L}_{\sigma-\mathrm{ind}}$-extension $\mathrm{PF}_{\sigma-\mathrm{ind}}$ of $\mathrm{PF}$ where we add the following axioms describing the new predicates: 
    \begin{itemize}
        \item If $(a_{1},\dots,a_{n},b_{1},\dots,b_{n})\in C_{n}$, then the polynomial $f_{\bar{a}}(X):=X^{n}+a_{1}X^{n-1}+\dots+a_{n}$ is irreducible over $F\models\mathrm{PF}$.
        \item If $(a_{1},\dots,a_{n},b_{1},\dots,b_{n})\in C_{n}$, then given some root $\alpha$ of $f_{\bar{a}}$ we have that $\tau_{\bar{b}}({\alpha}):=b_{1}\alpha^{n-1}+b_{2}\alpha^{n-2}+\dots+b_{n}$ is a root of $f$ different from $\alpha$.
        \item If $(a_{1},\dots,a_{n},b_{1},\dots,b_{n})\in C_{n}$ and $(a_{1},\dots,a_{n},d_{1},\dots,d_{n})\in C_{n}$, then $\bar{b}=\bar{d}$.
        \item For any root $\alpha$ of $f_{\bar{a}}$ we have $\tau^{n}_{\bar{b}}({\alpha})=\alpha$ and $\tau^{k}_{\bar{b}}({\alpha})\neq\alpha$ for $1\leq k\leq n$.
        \item The $\tau_{\bar{b}}$ are all compatible:\\ Let $(a_{1},\dots,a_{n},b_{1},\dots,b_{n})\in C_{n}$, $(e_{1},\dots,e_{k},g_{1},\dots,g_{k})\in C_{k}$, $k\vert n$, $h$ a polynomial over $\mathbb{Z}$, $\bar{u}$ a tuple in $F$ and $\alpha$ a root of $f_{\bar{a}}$ such that $\beta=h(\alpha,\bar{u})$ is a root of $f_{\bar{e}}$. Then $h(\tau_{\bar{b}}(\alpha),\bar{u})=\tau_{\bar{g}}(\beta)$.
       
    \end{itemize}
\end{definition}
\begin{remark}
    Note that in models of $\mathrm{PF}_{\sigma-\mathrm{ind}}$ the predicates $C_{n}$ are $\mathcal{L}_{\mathrm{ring}}$-definable with parameters. More concretely, $C_{n}$ is $\mathcal{L}_{\mathrm{ring}}$-definable from any $2n$-ary tuple which itself lies in $C_{n}$.
\end{remark}

\begin{fact}(1.13 Proposition A in \cite{acfa})\label{inducedstructureonfixedfield} Let $K\models\mathrm{ACFA}$ and let $F$ be its fixed field. Then we obtain exactly the induced structure on $F$ by considering $F$ as an $\mathcal{L}_{\sigma-\mathrm{ind}}$-structure where the $C_{n}$ are chosen such that $(\bar{a},\bar{b})\in C_{n}$ if $\sigma(\alpha)=\tau_{\bar{b}}({\alpha})$ (in the notation of Definition \ref{inducedstructurelanguage}). Moreover, $F\models\mathrm{PF}_{\sigma-\mathrm{ind}}$ as an $\mathcal{L}_{\sigma-\mathrm{ind}}$-theory.
\end{fact}
\begin{proof}
    Every predicate $C_{n}$ is clearly $\mathcal{L}_{\sigma}$-definable over $\emptyset$ and thus the induced structure necessarily contains the $C_{n}$. Now, we show that the $C_{n}$ indeed suffice to describe the induced structure.\\
    To do so, we fix $(K_{1},\sigma_{1})\models\mathrm{ACFA}$ and $(K_{2},\sigma_{2})\models\mathrm{ACFA}$. Further, we assume that $F_{1}\cong_{\mathcal{L}_{\sigma-\mathrm{ind}}} F_{2}$, i.e., we assume that the $F_{i}$ are isomorphic as $\mathcal{L}_{\sigma-\mathrm{ind}}$-structures with the interpretations of the $C_{n}$ as given in the statement of the fact. First, we want to show that then $(F_{1}^{\mathrm{alg}},\sigma_{1})\cong(F_{2}^{\mathrm{alg}},\sigma_{2})$ holds. But this follows from compactness: It suffices to consider the Galois extension $L_{n}$ of $F$ of degree $n$. Here we use that the predicate $C_{n}$ interprets the structure $(L_{n},\sigma\restriction_{L_{n}})$ as $\sigma\restriction_{L_{n}}$ is determined by its action on the roots of $f_{\bar{a}}$ where $f_{\bar{a}}$ is an irreducible polynomial of degree $n$ whose splitting field is $L_{n}$. Hence, we obtain $(L_{n}\cap F_{1}^{\mathrm{alg}},\sigma_{1}\restriction_{L_{n}})\cong (L_{n}\cap F_{2}^{\mathrm{alg}},\sigma_{2}\restriction_{L_{n}})$ and, consequently, $(F_{1}^{\mathrm{alg}},\sigma_{1})\cong(F_{2}^{\mathrm{alg}},\sigma_{2})$. Lemma 1.3 in \cite{acfa} then already implies that $(K_{1},\sigma_{1})\equiv_{F}(K_{2},\sigma_{2})$ after naming $F_{1}$, respectively $F_{2}$ by $F$. Thus, $F_{i}$ as an $\mathcal{L}_{\sigma-\mathrm{ind}}$ structure carries the induced structure from $(K_{i},\sigma)$ for $i=1,2$.
\end{proof}

\begin{lemma}\label{lemmamodelofofindisfixedfieldofacfa}
    For every model $F\models\mathrm{PF}_{\sigma-\mathrm{ind}}$ there is some $K\models\mathrm{ACFA}$ such that $F$ is the fixed field of $K$ together with the induced structure.
\end{lemma}

\begin{proof}
    First, we note that for all $n$ and $F_{n}$ the (unique) field extensions of degree $n$ of $F$, the predicate $C_{n}$ interprets a generator of the cyclic Galois group $\mathrm{Gal}(F_{n}/F)$ by the first four axioms. The last axiom then by compactness ensures that the $C_{n}$ together interpret a generator $\sigma$ of $\mathrm{Gal}(F^{\mathrm{alg}}/F)\cong\hat{\mathbb{Z}}$. Finally, it suffices to use that $(F^{\mathrm{alg}},\sigma)$ embeds in a model of $\mathrm{ACFA}$ with fixed field $F$ (see, e.g., Theorem 5.2 of \cite{afshordel}) together with Fact \ref{inducedstructureonfixedfield}.
\end{proof}


\begin{lemma}\label{elimofimaginariesinducedstructure}
    The $\mathcal{L}_{\sigma-\mathrm{ind}}$-theory $\mathrm{PF}_{\sigma-\mathrm{ind}}$ has elimination of imaginaries.
\end{lemma}
\begin{proof}
We can consider $F\models\mathrm{PF}_{\sigma-\mathrm{ind}}$ as the fixed field of a model of $\mathrm{ACFA}$ by Lemma \ref{lemmamodelofofindisfixedfieldofacfa}. We have elimination of imaginaries in $\mathrm{ACFA}$ by (1.10) in \cite{acfa}. We consider some imaginary $b$ defined over $F$ in $\mathcal{L}_{\sigma}$ and find some real tuple $c\in K^{n}$ interdefinable with $b$ (in $\mathcal{L}_{\sigma}$). Then $c\in F$, since $b$ is definable over $F$ and thus $\sigma(b)=b$. Now, as by Fact \ref{inducedstructureonfixedfield} the $\mathcal{L}_{\sigma-\mathrm{ind}}$-structure describes exactly the induced structure on $F$ and, moreover, $F$ is stably embedded, this interdefinability is witnessed by an $\mathcal{L}_{\sigma-\mathrm{ind}}$-formula which completes the proof.
\end{proof}
The fact that the induced structure strictly extends the (pure) field structure on the fixed field manifests itself as follows. The definable closure of a subset $A$ of the fixed field can be larger when computed in $\mathcal{L}_{\sigma}$ than when it is computed in $\mathcal{L}_{\mathrm{ring}}$. 

\begin{fact}(1.13 and 3.7 in \cite{acfa})\label{factdclwithandwithoutsigma}
    Let $A$ be a subfield of $F$. We then have $\mathrm{acl}_{\mathrm{PF}}(A)=\mathrm{acl}_{\mathrm{ACFA}}(A)\cap F=A^{\mathrm{alg}}\cap F=:\Tilde{A}$. Further, let $\mathcal{G}$ be the absolute Galois group of $A$, $N_{\mathcal{G}}(\overline{\langle\sigma\rangle})$, the normaliser in $\mathcal{G}$ of the closed subgroup generated by $\sigma$ and $C_{\mathcal{G}}(\sigma)$, the centraliser of $\sigma$.
    Finally, let $\mathrm{Aut}_{\mathcal{L}_{\mathrm{ring}}}(\Tilde{A}/A)$ and $\mathrm{Aut}_{\mathcal{L}_{\sigma-\mathrm{ind}}}(\Tilde{A}/A)$ be the group of the corresponding structure automorphisms over $A$. Then, the following hold:
    \begin{itemize}
        \item For all $\tau\in N_{\mathcal{G}}(\overline{\langle\sigma\rangle})$, the equality $\tau(\Tilde{A})=\Tilde{A}$ holds. 
        \item  $\mathrm{Aut}_{\mathcal{L}_{\mathrm{ring}}}(\Tilde{A}/A)\cong\Tilde{N}_{\mathcal{G}}(\overline{\langle\sigma\rangle})$, where the latter denotes the restriction to $\Tilde{A}$ of the elements from $N_{\mathcal{G}}(\overline{\langle\sigma\rangle})$. We can identify $\Tilde{N}_{\mathcal{G}}(\overline{\langle\sigma\rangle})\cong N_{\mathcal{G}}(\overline{\langle\sigma\rangle})/\overline{\langle\sigma\rangle}$.
        \item $\mathrm{dcl}_{\mathrm{PF}}(A)=\mathrm{Fix}(N_{\mathcal{G}}(\overline{\langle\sigma\rangle}))=\mathrm{Fix}(\Tilde{N}_{\mathcal{G}}(\overline{\langle\sigma\rangle}))$.
        \item  $\mathrm{Aut}_{\mathcal{L}_{\sigma-\mathrm{ind}}}(\Tilde{A}/A)\cong\Tilde{C}_{\mathcal{G}}(\sigma)$ where the latter denotes the restriction to $\Tilde{A}$ of the elements from $C_{\mathcal{G}}(\sigma)$. We can identify $\Tilde{C}_{\mathcal{G}}(\sigma)\cong C_{\mathcal{G}}(\sigma)/\overline{\langle\sigma\rangle}$.
        \item $\mathrm{dcl}_{\mathrm{ACFA}}(A)=\mathrm{dcl}_{\mathrm{PF}_{\sigma-\mathrm{ind}}}(A)=\mathrm{Fix}(C_{\mathcal{G}}(\sigma))=\mathrm{Fix}(\Tilde{C}_{\mathcal{G}}(\sigma))$.
    \end{itemize}
\end{fact}

We want to describe definable functions that we add to $\mathcal{L}_{\mathrm{ring}}$ (resp. $\mathcal{L}_{\sigma-\mathrm{ind}}$) to obtain quantifier elimination for $\mathrm{PF}$ (resp. $\mathrm{PF}_{\sigma-\mathrm{ind}}$) as well as that substructures are definably closed in the enriched language. For the $\mathcal{L}_{\mathrm{ring}}$-part we use the description and notation as in 3.2 of \cite{Hrushovski2021AxsTW}.

\begin{definition}\label{definitiondefclosedfunctionsforPF}
We work in $F\models\mathrm{PF}$. Let $P(\bar{X},Y)$ and $Q(\bar{X},Y)$ be integral polynomials. We define $\kappa_{P,Q}(\bar{a})=b$, if there exists at least one $d\in F$ such that $P(\bar{a},d)=0$ and for any such $d$ we have $Q(\bar{a},d)=b$. If no such $b$ exists or if there is no root of $P(\bar{a},X)$ in $F$, we set $\kappa_{P,Q}(\bar{a})=0$.
\end{definition}

\begin{definition}\label{definitiondefclosedfunctionsforPFind}
    Let $F\models\mathrm{PF}_{\sigma-\mathrm{ind}}$. Let $n$ be a natural number. Let $P(\bar{X},Y),Q(\bar{X},Y)$ and $h_{i}(\bar{X},Y),r_{i}(\bar{X},Y),$ be integral polynomials for $1\leq i\leq n$.
    We define the set $\mathcal{D}_{P,\bar{h},\bar{r}}(\bar{a})$ to consist of all $d\in F$ such that $P(\bar{a},d)=0$ and such that the tuple $(h_{1}(\bar{a},d),\dots, h_{n}(\bar{a},d),r_{1}(\bar{a},d),\dots, r_{n}(\bar{a},d))$ lies in $C_{n}$. Then, we define $\kappa_{P,Q,\bar{h},\bar{r}}:F^{n}\rightarrow F$ as the map that sends an $n$-tuple $\bar{a}$ to $b$, if there exists at least one $d\in \mathcal{D}_{P,\bar{h},\bar{r}}(\bar{a})$ and for any such $d\in \mathcal{D}_{P,\bar{h},\bar{r}}(\bar{a})$ we have $Q(\bar{a},d)=b$. If there is no such $d$ or $b$, we set $\kappa_{P,Q,\bar{h},\bar{r}}(\bar{a})=0$.
    We allow the tuple $(\bar{h},\bar{r})$ to be empty, i.e., the set of all the $\kappa_{P,Q,\bar{h},\bar{r}}$ includes the functions $\kappa_{P,Q}$ from Definition \ref{definitiondefclosedfunctionsforPF}.

\end{definition}

The following result is stated in 3.2 of \cite{Hrushovski2021AxsTW}. We did not find a proof in the literature.

\begin{fact}\label{factqepseudofinitefieldsfunctional}
    Let $\mathcal{L}_{\mathrm{ring},\kappa}$ be the language $\mathcal{L}_{\mathrm{ring}}$ together with the functions $\kappa_{P,Q}$ from Definition \ref{definitiondefclosedfunctionsforPF}. The theory $\mathrm{PF}$ (together with the interpretations of the $\kappa_{P,Q}$ as above) has elimination of quantifiers, and substructures are definably closed. 
\end{fact}
\begin{proof}
    Quantifier elimination follows directly from Fact \ref{kiefeQE} as for $Q=1$ the equation $\kappa_{P,Q}(\bar{a})=1$ holds if and only if $F\models\exists t\,P(\bar{a},t)=0$ (and otherwise $\kappa_{P,Q}(\bar{b})=0$). Let $A$ be a substructure with relative algebraic closure $\Tilde{A}$. By Fact \ref{factdclwithandwithoutsigma} $b\in \Tilde{A}$ is in $\mathrm{dcl}(A)$ if and only if $b$ is fixed by any element in the group $N_{\mathcal{G}}(\overline{\langle\sigma\rangle})$. For every finite Galois extension $D/A$ we denote by $\Tilde{N}_{\mathrm{Gal}(D/A)}(\langle\sigma\rangle)$ the restriction of $N_{\mathrm{Gal}(D/A)}(\langle\sigma\rangle)$ to $D\cap \Tilde{A}$, i.e., $\Tilde{N}_{\mathrm{Gal}(D/A)}(\langle\sigma\rangle)\cong N_{\mathrm{Gal}(D/A)}(\langle\sigma\rangle)/\langle\sigma\rangle$. As $N_{\mathcal{G}}(\overline{\langle\sigma\rangle})$ is isomorphic to the profinite group $\varprojlim N_{\mathrm{Gal}(D/A)}(\langle\sigma\rangle)$ (with the inverse system ranging over Galois extensions $D/A$)\footnote{This follows from the fact that if $\tau\not\in N_{\mathcal{G}}(\overline{\langle\sigma\rangle})$, then $\tau\Tilde{\sigma}\tau^{-1}$ is not a power of $\Tilde{\sigma}$ in $\mathrm{Gal}(D/A)$ for some $D/A$ Galois where $\Tilde{\sigma}$ is the restriction of $\sigma$ to $D$.} we find some finite Galois extension $C/A$ with $c\in B$ such that $N_{\mathrm{Gal}(C/A)}(\langle\sigma\rangle)$ fixes $b$. Let $\Tilde{C}$ denote $C\cap\Tilde{A}$.\\
    Now $\Tilde{N}_{\mathrm{Gal}(C/A)}(\langle\sigma\rangle)\cong \mathrm{Aut}_{\mathcal{L}_{\mathrm{ring}}}(\Tilde{C}/A)$. Hence, it follows that $\mathrm{Aut}_{\mathcal{L}_{\mathrm{ring}}}(\Tilde{C}/A)$ fixes $b$. Let $e$ be such that $\Tilde{C}=A(e)$. Let $P(X)$ be the minimal polynomial of $e$ over $A$, let $\bar{a}$ be its tuple of coefficients and $Q$ be such that $Q(e,\bar{a})=b$ (after possibly expanding $\bar{a}$). There is a bijection between $\mathrm{Aut}_{\mathcal{L}_{\mathrm{ring}}}(\Tilde{C}/A)$ and the roots of $P(X)$ that lie in $\Tilde{C}$ via the map that sends $\tau\in \mathrm{Aut}_{\mathcal{L}_{\mathrm{ring}}}(\Tilde{C}/A)$ to $\tau(e)$. Then, it follows that $\kappa_{P,Q}(\bar{a})=b$ and thus $b\in A$.  
\end{proof}

\begin{lemma}\label{lemmainducedstructuredefinableclosure}
       Let $\mathcal{L}_{\sigma-\mathrm{ind},\kappa}$ be the language $\mathcal{L}_{\sigma-\mathrm{ind}}$ together with the functions $\kappa_{P,Q,\bar{h},\bar{r}}$ from Definition \ref{definitiondefclosedfunctionsforPFind}. The theory $\mathrm{PF}_{\sigma-\mathrm{ind}}$ (together with the interpretations of the $\kappa_{P,Q,\bar{h},\bar{r}}$ as above) has elimination of quantifiers and substructures are definably closed.
\end{lemma}
\begin{proof}
    Let $A$ be a substructure of a model $F\models\mathrm{PF}_{\sigma-\mathrm{ind}}$, $\Tilde{A}$ its relative algebraic closure, and $\bar{a}$ denote a tuple in $A$.
    For quantifier elimination, we first note that for $Q=1$ the equation $\kappa_{P,Q,\bar{h},\bar{r}}(\bar{a})=1$ holds if and only if $F\models \exists t\in \mathcal{D}_{P,\bar{h},\bar{r}}(\bar{a})$. Now, the formulas of the form $\exists t\in \mathcal{D}_{P,\bar{h},\bar{r}}(\bar{x})$ describe the isomorphism type of $(\mathrm{acl}_{\mathrm{ACFA}}(A),\sigma)$ as a difference field (where we interpret $\sigma$ using the $C_{n}$ as before). Thus, by Fact \ref{basicfactsonacfa} (quantifier reduction in $\mathrm{ACFA}$) and the fact that $F$ carries exactly the induced structure (Fact \ref{inducedstructureonfixedfield}) quantifier elimination of $\mathrm{PF}_{\sigma-\mathrm{ind}}$ in the language $\mathcal{L}_{\sigma-\mathrm{ind},\kappa}$ already follows.\\
    Let $A$ be a substructure and $b\in \mathrm{dcl}(A)$. With the same argument as in the proof of Fact \ref{factqepseudofinitefieldsfunctional} we find some Galois extension $B/A$ with $\Tilde{B}:=B\cap F$ such that $b$ is fixed by $C_{\mathrm{Gal}(B/A)}(\sigma)$ and $\mathrm{Aut}_{\mathcal{L}_{\sigma-\mathrm{ind}}}(\Tilde{B}/A)\cong C_{\mathrm{Gal}(B/A)}(\sigma)/\langle\sigma\rangle$.
    By compactness (and finiteness of $\mathrm{Aut}_{\mathcal{L}_{\sigma-\mathrm{ind}}}(B/A)$) we can find finitely many $n_{i}\in\mathbb{N}$ and a finite number of tuples from $\Tilde{B}$ such that $\mathrm{Aut}_{\mathcal{L}_{\sigma-\mathrm{ind}}}(\Tilde{B}/A)$ is given by those elements from $\mathrm{Aut}_{\mathcal{L}_{\mathrm{ring}}}(\Tilde{B}/A)$ that preserve the relations $C_{n_{i}}$ for the tuples from $B$ that we obtained. By the compatibility axiom in Definition \ref{inducedstructurelanguage} we can reduce to working with $n$ the least common multiple of the $n_{i}$ and one $2n$-tuple $\bar{\gamma}\in\Tilde{B}\cap C_{n}$. Let $e\in F$ be a generator of $B$ and let $P(X)$ be the minimal polynomial of $e$ over $A$, let $\bar{a}$ be its tuple of coefficients and $Q$ be such that $Q(e,\bar{a})=b$ (after possibly enlarging $\bar{a}$). Further, let $\bar{h},\bar{r}$ 
    be such that $\gamma_{i}=h_{i}(\bar{a},e)$ and $\gamma_{2i}=r_{i}(\bar{a},e)$ for all $1\leq i \leq n$ (again after possibly expanding $\bar{a}$). Then, we have $\kappa_{P,Q,\bar{h},\bar{r}}(\bar{a})=b$ and thus $b\in A$.
\end{proof}
\begin{corollary}\label{corollaryGaloisoverdefinableclosure}
    Let $F\models\mathrm{PF}_{\mathrm{ind}}$ and $A\subseteq F$ be a substructure of $F$ in $\mathcal{L}_{\sigma-\mathrm{ind},\kappa}$ with relative algebraic closure $\Tilde{A}$. Then $\Tilde{A}$ is a (possibly infinite) Galois extension of $A$ and $\mathrm{Gal}(\Tilde{A}/A)=\mathrm{Aut}_{\mathcal{L}_{\sigma-\mathrm{ind}}}(\Tilde{A}/A)$.
\end{corollary}
\begin{proof}
     By Lemma \ref{lemmainducedstructuredefinableclosure} $A$ is definably closed as a structure in the language $\mathcal{L}_{\sigma-\mathrm{ind}}$.\footnote{Recall that the $\kappa_{P,Q,\bar{h},\bar{r}}$ are $\mathcal{L}_{\sigma-\mathrm{ind}}$-definable functions.} In other words, the fixed field of the group of $\mathcal{L}_{\sigma-\mathrm{ind}}$-structure automorphisms $\mathrm{Aut}_{\sigma-\mathrm{ind}}(\Tilde{A}/A)$ is given by $A$ and thus $\Tilde{A}$ is Galois over $A$ and $\mathrm{Gal}(\Tilde{A}/A)=\mathrm{Aut}_{\sigma-\mathrm{ind}}(\Tilde{A}/A)$.
\end{proof}


\section{The theory $\mathrm{ACFA}^{+}$}\label{sectionintroacfapluschaptertwo}

In this section we introduce the theory $\mathrm{ACFA}^{+}$ as well as the theory $\mathrm{PF}^{+}$ that Hrushovski studies in \cite{Hrushovski2021AxsTW}. From his results we deduce a quantifier elimination result for the theory $\mathrm{PF}^{+}_{\sigma-\mathrm{ind}}$ (see Definition \ref{definitionPFplusinduces}).

\begin{convention}
    From now on and for the rest of this article we restrict ourselves to characteristic $0$ in the following sense: $\mathrm{ACFA}$ and $\mathrm{PF}$ will by convention always contain a set of axioms stating that the characteristic is $0$.
\end{convention}

\begin{notation}
    We write $S^{1}$ for the unit circle and $\mathbb{T}^{n}\cong S^{1}\times\cdots\times S^{1}$ for the $n$-dimensional complex torus.
\end{notation}

We define the language $\mathcal{L}_{\sigma,+}$ to consist of a sort $K$ (with equality treated in the usual way) for the difference field with the ring operations and a symbol $\sigma$ for the automorphism. In addition, there is one unary relation $\Psi:K\rightarrow S^{1}\cup\{0\}$.
We denote by $\Psi^{(n)}:K^{n}\rightarrow\{S^{1}\cup\{0\}\}^{n}$ the map given by $(x_{1},\dots,x_{n})\rightarrow(\Psi(x_{1}),\dots,\Psi(x_{n}))$.

\begin{definition}
    We say that a rational hyperplane over $F$ in $\mathbb{A}^{n}$ is a variety defined by an equation of the form $\sum_{1\leq i\leq n}z_{i}X_{i}=b$ where $b\in F$ and $z_{i}\in\mathbb{Z}$ for $1\leq i\leq n$ and $z_{i}\neq 0$ for some $1\leq i\leq n$. We say that it has height $\leq m$, if $|z_{i}|\leq m$ for all $1\leq i\leq n$.
\end{definition}

\begin{definition}\label{definitionacfaplus}
We define the theory $\mathrm{ACFA}^{-}$ to consist of the $\mathcal{L}_{\sigma}$-theory $\mathrm{ACFA}$ together with axioms stating that $\Psi(K\backslash F)=0$ where $F:=\mathrm{Fix}(\sigma)$ and moreover that $\Psi\restriction_{F}:(F,+)\rightarrow(S^{1},\cdot)$ is a group homomorphism. The theory $\mathrm{ACFA}^{+}$ then extends $\mathrm{ACFA}^{-}$ by the following set of axioms:\\

\noindent
($\star$) Let $n,m\in\mathbb{N}$. Let $h\in\mathbb{Q}[z_{1},z_{1}^{-1},\dots,z_{n},z_{n}^{-1}]$ be a finite Fourier series (Laurent polynomial) with degrees $\leq m$, real-valued on $\mathbb{T}^{n}$, with no constant term. For any absolutely irreducible curve $C$ over $F$ with $C\subset\mathbb{A}^{n}$, not contained in any rational hyperplane over $F$ of height at most $m$,
\[\sup\{h(\Psi^{(n)}(\bar{x}))\;:\;\bar{x}\in C(F)\}\geq 0.\]

\end{definition}

\begin{fact}\label{factaxiomisexpressible}(Lemma 3.5 in \cite{Hrushovski2021AxsTW}.)\label{factequivalenceaxiomtodensevaluesoncruve}
   Let $F$ be a field of characteristic $0$ with an additive character $\Psi$. Then, condition $(\star)$ from Definition \ref{definitionacfaplus} is true for $F$ if and only if the following holds:\\
   For any absolutely irreducible curve $C\subset \mathbb{A}^{n}$ defined over $F$ that is not contained in a rational hyperplane over $F$, the set $\Psi^{(n)}(C(F))$ is dense in $\mathbb{T}^{n}$ (in the euclidean topology).
   
\end{fact}

\begin{definition}\label{definitionPFplusinduces}
    Let $\mathcal{L}_{+}$ be the ring language extended by a unary relation symbol $\Psi$ (with values in $\mathbb{T}$). Then Definition \ref{definitionacfaplus} contains the theory $\mathrm{PF}^{+}$ from \cite{Hrushovski2021AxsTW} in the following sense: $\mathrm{PF}^{+}$ states that $F$ is pseudofinite\footnote{It suffices to give axioms for $\mathrm{Gal}(F)=\hat{\mathbb{Z}}$, pseudo-algebraic closedness follows from $(\star)$.}, $\Psi$ is a homomorphism for the additive group of the field and ($\star$) holds as above.\\
    Similarly we define $\mathcal{L}_{\sigma-\mathrm{ind},+}$ and $\mathrm{PF}_{\sigma-\mathrm{ind}}^{+}$ as the extension of $\mathrm{PF}_{\sigma-\mathrm{ind}}$ by axioms stating that $\Psi$ is a homomorphism for the additive group of the field and ($\star$) holds as above.
    
\end{definition}

\begin{definition}\label{definitionPsinsym}(See 3.2 in \cite{Hrushovski2021AxsTW}.)
    Let $\bar{a}=(a_{1},\dots,a_{n})$ be an n-tuple in $F\models\mathrm{PF}^{+}$. We denote by $Z(\bar{a})$ the unordered tuple containing the elements from $\{x\in F\,|\,x^{n}+a_{1}x^{n-1}+\cdots+a_{n}=0\}$ such that every element occurs with its multiplicity.
    We define the terms $\Psi_{\mathrm{sym}}^{n}(y_{1},\dots,y_{n})$ as follows:
    \[\Psi_{\mathrm{sym}}^{n}(a_{1},\dots,a_{n}):=\sum_{x\in Z(\bar{a})}\Psi(x).\]
    We denote by $\mathcal{L}_{\kappa,+}^{\mathrm{sym}}$ the language consisting of $\mathcal{L}_{+}$ together with n-ary predicate symbols $\Psi_{\mathrm{sym}}^{n}$ for all $n\in \mathbb{N}$ as well as symbols for the $\kappa_{P,Q}$ as in Definition \ref{definitiondefclosedfunctionsforPF}. The $\mathcal{L}_{\kappa,+}^{\mathrm{sym}}$-theory $\mathrm{PF}^{+}$ extends the $\mathcal{L}_{+}$-theory $\mathrm{PF}^{+}$ by the interpretations of the new symbols as described above.
\end{definition}

\begin{definition}
We denote by $\mathcal{L}_{\sigma-\mathrm{ind},\kappa,+}^{\mathrm{sym}}$ the language $\mathcal{L}_{\sigma-\mathrm{ind},+}$ together with n-ary predicate symbols symbols $\Psi_{\mathrm{sym}}^{n}$ for all $n\in \mathbb{N}$ as well as symbols for the $\kappa_{P,Q,\bar{h},\bar{r}}$ as in Definition \ref{definitiondefclosedfunctionsforPFind} and again we consider $\mathrm{PF}^{+}_{\mathrm{ind}}$ as an $\mathcal{L}_{\sigma-\mathrm{ind},\kappa,+}^{\mathrm{sym}}$-theory by adding the interpretations of the new symbols.
\end{definition}

The following is the quantifier elimination given by Hrushovski in \cite{Hrushovski2021AxsTW}. 

\begin{fact}(Proposition 3.9 in \cite{Hrushovski2021AxsTW}.)\label{factqepfplus}
    The theory $\mathrm{PF}^{+}$ has quantifier elimination in the language $\mathcal{L}_{\kappa,+}^{\mathrm{sym}}$.
\end{fact}

Our goal will be to obtain quantifier elimination in $\mathrm{PF}^{+}_{\mathrm{ind}}$ working in $\mathcal{L}_{\sigma-\mathrm{ind},\kappa,+}^{\mathrm{sym}}$. We will follow the structure of the proof of Fact \ref{factqepfplus} but point out why it was crucial to have added the symbols $\kappa_{P,Q,\bar{h},\bar{r}}$ when working with the induced structure. We start by stating the essential ingredient of \cite{Hrushovski2021AxsTW}.

\begin{fact}(Part of the proof of Lemma 3.8 in \cite{Hrushovski2021AxsTW}.)\label{factextensioncharactertogaloisext}
Let $A\subseteq B\subseteq F\models\mathrm{PF}$ be subfields with $B$ a finite Galois extension of $A$. Let $\Psi,\Theta$ be two additive characters on $B$ that agree on $A$. If for all $n\in\mathbb{N}$ and $n$-tuples $\bar{a}$ in $A$ we have $\Psi_{\mathrm{sym}}^{n}(\bar{a})=\Theta_{\mathrm{sym}}^{n}(\bar{a})$, then there is some $\tau\in \mathrm{Gal}(B/A)$ such that $\Psi(b)=\Theta(\tau(b))$ for all $b\in B$.
\end{fact}

\begin{lemma}\label{lemmaextensiontoalgclosurepfplusind}(See Lemma 3.8 in \cite{Hrushovski2021AxsTW} for the statement for $\mathrm{PF}^{+}$.) Let $A$ and $A^{\prime}$ be isomorphic substructures of $F\models\mathrm{PF}^{+}_{\mathrm{ind}}$ in the language $\mathcal{L}_{\sigma-\mathrm{ind},\kappa,+}^{\mathrm{sym}}$. Then any isomorphism $\alpha$ between $A$ and $A^{\prime}$ extends to an isomorphism between the relative algebraic closures (in $F$) of $A$ and $A^{\prime}$.
    
\end{lemma}

\begin{proof}
    From quantifier elimination of $\mathrm{PF}_{\sigma-\mathrm{ind}}$ in $\mathcal{L}_{\sigma-\mathrm{ind},\kappa}$ (Lemma \ref{lemmainducedstructuredefinableclosure}), it follows that we obtain a field isomorphism extending $\alpha$ between the respective relative algebraic closures that preserves the $C_{n}$. We can thus reduce to the situation of working with $A$ and its relative algebraic closure $\Tilde{A}$ carrying two additive characters $\Psi_{1}$ and $\Psi_{2}$, that extend the one on $A$ and we have to find an $\mathcal{L}_{\sigma-\mathrm{ind}}$-structure automorphism $\tau$ of $\Tilde{A}$ over $A$ such that for all $c\in\Tilde{A}$ we have $\Psi_{1}(c)=\Psi_{2}(\tau(c))$. By Corollary \ref{corollaryGaloisoverdefinableclosure} we have $\mathrm{Gal}(\Tilde{A}/A)=\mathrm{Aut}_{\sigma-\mathrm{ind}}(\Tilde{A}/A)$.
    Thus, it suffices (by compactness of $\mathrm{Gal}(\Tilde{A}/A)$)) to show for any finite Galois extension $A\subset B\subset\Tilde{A}$ that there is some $\tau\in \mathrm{Gal}(B/A)$ satisfying $\Psi_{1}(c)=\Psi_{2}(\tau(c))$ for all $c\in B$. This is Fact \ref{factextensioncharactertogaloisext}. 
\end{proof}

The rest of the proof now follows as in
\cite{Hrushovski2021AxsTW}. Working with the $\mathcal{L}_{\sigma-\mathrm{ind}}$-structure instead of the field structure does not yield any complications, but to a certain extent is slightly easier. The latter is due to the fact that when working with the field structure at some point one has to choose a generator of the Galois group, which in our case has already been done.
We will now state the following embedding lemma which is the essential ingredient of the proof of quantifier elimination for $\mathrm{PF}^{+}_{\mathrm{ind}}$.

\begin{lemma}\label{lemmaembeddinglemma}
Let $F$ be an $\aleph_{1}$-saturated model of $\mathrm{PF}^{+}$ with Galois generator $\sigma\in \mathrm{Gal}(F^{\mathrm{alg}}/F)$, $K\subseteq F$ a countable substructure, relatively algebraically closed as a field in $F$. Let $E\supseteq K$ be an extension of $\mathcal{L}_{+}$-structures such that the following hold.
\begin{itemize}
    \item $E$ is a field of transcendence degree $1$ over $K$.
    \item $K$ is relatively algebraically closed in $E$.
    \item $\Psi:(E,+)\rightarrow S^{1}$ is a group homomorphism.
    \item There is a topological generator $\sigma^{\prime}$ of $\mathrm{Gal}(E^{\mathrm{alg}}/E)$ such that $\sigma^{\prime}\mapsto\sigma$ induces an isomorphism $\mathrm{Gal}(E^{\mathrm{alg}}/E)\rightarrow \mathrm{Gal}(F^{\mathrm{alg}}/F)$.
\end{itemize}
Then, we find a field embedding $\tau$ of $E^{\mathrm{alg}}$ over $K^{\mathrm{alg}}$ into $F^{\mathrm{alg}}$ such that $\tau(\sigma^{\prime}(a))=\sigma(\tau(a))$ for all $a\in E^{\mathrm{alg}}$, $\tau(E)$ is relatively algebraically closed in $F$ and $\tau$ preserves $\Psi$.
\end{lemma}

\begin{proof}
This follows exactly as in the proof of Proposition 3.9 in \cite{Hrushovski2021AxsTW} as in the proof generators of the Galois groups are chosen and fixed at the beginning. In Lemma 4.7 in \cite{ludwig2025pseudofinitefieldsadditivemultiplicative} we prove a more general embedding lemma in the same presentation as the above statement.
\end{proof}


\begin{fact}(Lemma 3.11 in \cite{Hrushovski2021AxsTW}.)\label{factexistenceoffullextensionoverQ}
Let $T$ be a completion of $\mathrm{PF}$. Then there exists $F\models T$ and a subfield $F_{1}$ of $F$ of transcendence degree $1$ over $\mathbb{Q}$ such that the restriction $\mathrm{Gal}(F^{\mathrm{alg}}/F)\rightarrow \mathrm{Gal}(F_{1}^{\mathrm{alg}}/F_{1})$ is an isomorphism.
    
\end{fact}

\begin{corollary}\label{corollaryqepfplusind}
    The theory $\mathrm{PF}^{+}_{\mathrm{ind}}$ has quantifier elimination in the language $\mathcal{L}_{\sigma-\mathrm{ind},\kappa,+}^{\mathrm{sym}}$.
\end{corollary}
\begin{proof}
As $\mathcal{L}_{\sigma-\mathrm{ind},\kappa,+}^{\mathrm{sym}}$ is countable, it suffices to show that we can extend every isomorphism between countable substructures $K_{1}$ and $K_{2}$ to a back-and-forth system between $\aleph_{1}$-saturated models $F_{1}$ and $F_{2}$. 
Without loss of generality, by Lemma \ref{lemmaextensiontoalgclosurepfplusind}, we can assume that $F_{1}$ and $F_{2}$ have a common countable substructure $K$ that is a relatively algebraically closed subfield of $F_{1}$ and $F_{2}$. Then, it suffices to show that we can embed every $E_{1}\subseteq F_{1}$, relatively algebraically closed, and of transcendence degree $1$ over $K$ into $F_{2}$ so that the image is relatively algebraically closed again.
If the map $\mathrm{Gal}(F_{2}^{\mathrm{alg}}/F_{2})\rightarrow \mathrm{Gal}(E_{1}^{\mathrm{alg}}/E_{1})$ obtained by identifying the Galois generators (which are specified as we work in $\mathcal{L}_{\sigma-\mathrm{ind}}$) is an isomorphism, this works by Lemma \ref{lemmaembeddinglemma}. (Note that Lemma \ref{lemmaembeddinglemma} indeed ensures that the $\mathcal{L}_{\sigma-\mathrm{ind},\kappa,+}^{\mathrm{sym}}$-structure is preserved.) Once we can make sure that $\mathrm{Gal}(F_{1}^{\mathrm{alg}}/F_{1})\rightarrow \mathrm{Gal}(K^{\mathrm{alg}}/K)$ and  $\mathrm{Gal}(F_{2}^{\mathrm{alg}}/F_{2})\rightarrow \mathrm{Gal}(K^{\mathrm{alg}}/K)$ are both isomorphisms, this condition on $E_{1}$ is automatically fulfilled (see e.g. the proof of Lemma 3.10 in \cite{Hrushovski2021AxsTW} for the argument). 
Hence, if we can show that we can always find an $K$-isomorphism between two relatively algebraically closed extensions $K_{1}/K$,  $K_{2}/K$ in $F_{1}$,$F_{2}$, respectively, such that $\mathrm{Gal}(F_{1}^{\mathrm{alg}}/F_{1})\rightarrow \mathrm{Gal}(K_{1}^{\mathrm{alg}}/K_{1})$ and  $\mathrm{Gal}(F_{2}^{\mathrm{alg}}/F_{2})\rightarrow \mathrm{Gal}(K_{2}^{\mathrm{alg}}/K_{2})$ are isomorphisms, we are done.
This now follows from Fact \ref{factexistenceoffullextensionoverQ} as in \cite{Hrushovski2021AxsTW} (last paragraph of proof of Lemma 3.9).
\end{proof}

\subsection{Quantifier elimination and consequences in $\mathrm{ACFA}^{+}$}\label{sectionQEacfapluschaptertwo}

\begin{definition}
      Let $\mathcal{L}_{\sigma,\kappa}$ be the language $\mathcal{L}_{\sigma}$ together with the symbols $C_{n}$ and $\kappa_{P,Q,\bar{h},\bar{r}}$ with their interpretations on (cartesian products of) $F$
      as given in Definition \ref{inducedstructurelanguage} and in Definition \ref{definitiondefclosedfunctionsforPFind} respectively.
      Let $\mathcal{L}_{\sigma,\kappa}^{\mathrm{qe}}$ be a definitional expansion of $\mathcal{L}_{\sigma,\kappa}$ such that $\mathrm{ACFA}^{+}$ has quantifier elimination in it\footnote{e.g., the one obtained by adding all the predicates of the form described in Corollary \ref{corollaryexplicitqeforacfa}} together with the functions $f_{\psi,\phi}$ as given by Lemma \ref{factdefbfunctionofcanocparamter} (here we use that $F$ is stably embedded) for all formulas $\phi$.
      Finally, let $\mathcal{L}_{\sigma,\kappa,+}^{\mathrm{sym}}$ be the language $\mathcal{L}_{\sigma,\kappa}^{\mathrm{qe}}$ together with the symbols $\Psi_{\mathrm{sym}}^{n}$ (including a symbol for $\Psi$). Then, let $\mathrm{ACFA}^{+}$ be the $\mathcal{L}_{\sigma,\kappa,+}^{\mathrm{sym}}$-theory obtained by adding the interpretations of the added symbols to the axioms in Definition \ref{definitionacfaplus}.
\end{definition}

\begin{theorem}\label{theoremQEacfaplus}

    The theory $\mathrm{ACFA}^{+}$ has quantifier elimination in the language $\mathcal{L}_{\sigma,\kappa,+}^{\mathrm{sym}}$. 
\end{theorem}

\begin{proof}We show that we can apply Corollary \ref{corollaryabstractQEconition}: First note that in $\mathcal{L}_{\sigma,\kappa}^{\mathrm{qe}}$ all the added symbols are $\emptyset$-definable in $\mathcal{L}_{\sigma}$ and thus $\mathrm{ACFA}$ has elimination of imaginaries ((1.10) in \cite{acfa}) in $\mathcal{L}_{\sigma,\kappa}^{\mathrm{qe}}$. By construction, it has quantifier elimination. The fixed field $F\subset K\models\mathrm{ACFA}$ is stably embedded by Fact \ref{factstableembeddednessacfaclassicallogic} (and remains so in a definitional expansion). In addition, $F$ as an $\mathcal{L}_{\sigma-\mathrm{ind},\kappa}$-structure carries the induced structure from $K$. Finally, $\mathcal{L}_{\sigma-\mathrm{ind},\kappa,+}^{\mathrm{sym}}$ is a relational expansion of $\mathcal{L}_{\sigma-\mathrm{ind},\kappa}$ and $\mathrm{PF}_{\sigma-\mathrm{ind}}^{+}$ has quantifier elimination in it. Hence, Corollary \ref{corollaryabstractQEconition} can be applied, and it follows that $\mathrm{ACFA}^{+}$ has quantifier elimination in $\mathcal{L}_{\sigma,\kappa,+}^{\mathrm{sym}}$.
\end{proof}

From this quantifier-elimination result we directly obtain the following two corollaries.

\begin{corollary}\label{characterisationeltequivacfaplus}

Let $(K_{1},\sigma_{1},\Psi_{1}),(K_{2},\sigma_{2},\Psi_{2})$ be two models of $\mathrm{ACFA}^{+}$
 and $E$ a common algebraically closed difference subfield on which $\sigma_{1}$ and $\sigma_{2}$ agree and such that $\Psi_{1}, \Psi_{2}$ agree on $\mathrm{Fix}(E)$, then
 \[(K_{1},\sigma_{1},\Psi_{1})\equiv_{E} (K_{2},\sigma_{2},\Psi_{2}).
 \]
 \end{corollary}
\begin{proof}
    Since $E$ is algebraically closed, it carries the same $\mathcal{L}_{\sigma,\kappa,+}^{\mathrm{sym}}$-structure if computed from $(\sigma_{1},\Psi_{1})$ or $(\sigma_{2},\Psi_{2})$.
\end{proof}

 \begin{corollary}\label{corollarymodelcompleteness}
     $\mathrm{ACFA}^{+}$ is model complete.
 \end{corollary}
 \begin{proof}
     Direct from Corollary \ref{characterisationeltequivacfaplus}.
 \end{proof}

\begin{theorem}\label{modelcompanion}
Every difference field $(K,\sigma)$ with a fixed field $F$ and a homomorphism $\Psi:(F,+)\rightarrow (S^{1},\cdot)$ embeds into a model of $\mathrm{ACFA}^{+}$. Thus, $\mathrm{ACFA}^{+}$ is the model companion of the theory of difference fields with an additive character on the fixed field.
\end{theorem}
\begin{proof}
For the first assertion, we note that we can find some $(F^{*},\Psi)\models\mathrm{PF}^{+}$ that extends the structure $(F,\Psi)$ by 3.13 of \cite{Hrushovski2021AxsTW}. We can assume (by possibly moving $K$) that $K$ and $F^{*}$ are linearly disjoint over $F$. Thus, by extending $\sigma$ (in the unique way, see Theorem 1.3 of \cite{acfa}) to the field amalgam of $K$ and $F^{*}$ we obtain a difference field $(\Tilde{K},\sigma)$ with fixed field $F^{*}$. Then by Theorem 5.2 in \cite{afshordel} we can embed $(\Tilde{K},\sigma)$ in a model $(L,\sigma)$ of $\mathrm{ACFA}$ with fixed field $F^{*}$. Then $(L,\sigma)$ together with the predicate $\Psi$ on $F^{*}$ yields a model of $\mathrm{ACFA}^{+}$. The second statement then already follows from model-completeness of $\mathrm{ACFA}^{+}$.
\end{proof}

From Corollary \ref{characterisationeltequivacfaplus} we obtain, as in the classical case for $\mathrm{ACFA}$, the following two consequences:
\begin{corollary}\label{substructureeltequiv}
Let $\mathcal{K}_{1}:=(K_{1},\sigma_{1},\Psi_{1})$ and $\mathcal{K}_{2}:=(K_{2},\sigma_{2},\Psi_{2})$ be models of $\mathrm{ACFA}^{+}$. Then $\mathcal{K}_{1}$ and $\mathcal{K}_{2}$ are elementarily equivalent if and only if
 \[(\mathbb{Q}^{\mathrm{alg}},\sigma_{1}\restriction_{\mathbb{Q}^{\mathrm{alg}}},\Psi_{1}\restriction_{\mathbb{Q}^{\mathrm{alg}}})\cong (\mathbb{Q}^{\mathrm{alg}},\sigma_{2}\restriction_{\mathbb{Q}^{\mathrm{alg}}},\Psi_{2}\restriction_{\mathbb{Q}^{\mathrm{alg}}}).
 \]
\end{corollary}
\begin{corollary}\label{typedescription}
Let $E$ be a substructure of a model $K$ of $\mathrm{ACFA}^{+}$ and $a,b$ tuples from $K$. Then $tp(a/E)=tp(b/E)$ if and only if there is an $E$-isomorphism $\alpha$ between $\mathrm{acl}_{\sigma}(E(a))$ and $\mathrm{acl}_{\sigma}(E(b))$ preserving the values of $\Psi$ and sending $a$ to $b$.
\end{corollary}
Finally, we see that the fixed field itself is stably embedded again.

\begin{lemma}\label{inducedstructureplusfixedfield}
    Let $K\models\mathrm{ACFA}^{+}$ and let $F$ be its fixed field. Then we obtain exactly the induced structure on $F$ by considering $F$ as an $\mathcal{L}_{+,\sigma-\mathrm{ind}}$-structure where the $C_{n}$ are chosen such that $(\bar{a},\bar{b})\in C_{n}$ if $\sigma(\alpha)=\Tilde{\sigma}(\alpha)$ for $\alpha$ and $\Tilde{\sigma}$ is defined as in Definition \ref{inducedstructurelanguage}. Further, $F$ is stably embedded and we even have that the restriction of every (parameter-)definable predicate to $F$ is definable with parameters from $F$ only using the language $\mathcal{L}_{+}$.
\end{lemma}
\begin{proof}
    We can apply Corollary \ref{corollarypropertiesofenrihcedstableembeddedstructures}.
\end{proof}

\begin{lemma}\label{qe-conservative}
The theory $\mathrm{ACFA}^{+}$ is a conservative expansion of $\mathrm{ACFA}$ with respect to definable sets. In other words, every set definable in some model $K\models\mathrm{ACFA}$ is already definable in $\mathcal{L}_{\sigma}$.
\end{lemma}
\begin{proof}
    This is a direct application of Lemma \ref{lemmaconservativeextonstablembset} using the corresponding result for $\mathrm{PF}^{+}$, that is, Corollary 4.5 in \cite{Hrushovski2021AxsTW}.
\end{proof}

\subsection{The asymptotic theory of $(\mathbb{F}_{q}^{\mathrm{alg}},\mathrm{Frob}_{q},\Psi_{q})$}\label{sectionlimittheory}
In this section, we will show that $\mathrm{ACFA}^{+}$ is the common theory of the algebraic closure of finite fields together with the Frobenius automorphism, the standard character (see Definition \ref{definitionstandardcharacter}) on the fixed field, and a set of axioms stating that the characteristic is $0$. In other words, we will see that Ax's theorem (Fact \ref{theoremaxtheorempsf}) generalises to our context.\\
Hrushovski already generalised Ax's work in two directions: He showed that $\mathrm{ACFA}$ is the common theory of the algebraic closure of finite fields with Frobenius in his highly important work \cite{elementaryfrobenius} by proving a generalisation of the Lang-Weil bounds. Further, he showed that $\mathrm{PF}^{+}$ is the common theory of finite fields with standard character in \cite{Hrushovski2021AxsTW}. Our result is merely a combination of those two. However, we will explain below that whereas one direction (ultraproducts of the above are models of $\mathrm{ACFA}^{+}$) directly follows from Hrushovski's results, the other direction (every completion of $\mathrm{ACFA}^{+}$ has such an ultraproduct as a model) is not completely immediate and we will have to extend Hrushvoski's argument in \cite{Hrushovski2021AxsTW} slightly. The necessary ingredients are however already present in \cite{Hrushovski2021AxsTW} and we only have to recombine them.

\begin{definition}\label{definitionstandardcharacter}
Given a prime $p$ the \textit{standard character} on $\mathbb{F}_{p}$ is defined by the following map:
\[\Psi_{p}:n+p\mathbb{Z} \;\rightarrow \exp(2\pi i \frac{n}{p}).\]
For an arbitrary finite field $\mathbb{F}_{q}$ with $q=p^{e}$ we define the standard character via
\[\Psi_{q}(x):=\Psi_{p}(\mathrm{Tr}(x))\]
where $\mathrm{Tr}:\mathbb{F}_{q}\rightarrow \mathbb{F}_{p}$ is the trace map.
\end{definition}
\begin{remark} (3.15 in \cite{Hrushovski2021AxsTW})
Any other additive character on $\mathbb{F}_{q}$ has the form $\Psi_{q}(cx)$ for a unique $c\in\mathbb{F}_{q}$.
\end{remark}
\begin{notation}
We denote by $\mathrm{Frob}_{q}$ the Frobenius-map on $\mathbb{F}_{q}^{\mathrm{alg}}$. In the following, when we write $\prod_{\mathcal{U}}\mathbb{F}_{p}$, then we will implicitly mean that $\mathcal{U}$ is a non-principal ultrafilter on the sets of primes and when we write $\prod_{\mathcal{U}}\mathbb{F}_{q}$ it will be a non-principal ultrafilter on the set of prime powers $\{p^e\,|\,p\;\text{prime},\;e\in\mathbb{N}\}$.  We will write $\mathrm{Th}((\mathbb{F}_{p})_{p})$ (resp. $\mathrm{Th}((\mathbb{F}_{q})_{q})$) for the common theory of all prime (resp. finite) fields.
\end{notation}
\begin{definition}
    The \textit{limit theory} $\mathrm{Th}_{0}(\mathcal{K})$ of a class $\mathcal{K}$ of (the enriched algebraic closure) of finite fields will be the common theory of the elements of that class together with a set of axioms stating that the characteristic is $0$. 
\end{definition}
The main theorem of this section will then be as follows.
\begin{theorem}\label{maintheoremlimittheory}
$\mathrm{ACFA}^{+}$ is the limit theory of the algebraic closure of finite fields with standard character and Frobenius, and of the algebraic closure of
prime fields with arbitrary character and Frobenius. In other words,
\[\mathrm{ACFA}^{+}=\mathrm{Th}_{0}((\mathbb{F}_{q}^{\mathrm{alg}},\mathrm{Frob}_{q},\Psi_{q})_{q})=\mathrm{Th}_{0}((\mathbb{F}_{p}^{\mathrm{alg}},\mathrm{Frob}_{p},\Psi_{p}(ax))_{p,a\in\mathbb{F}_{p}})\]
\end{theorem}
To start with, we show that one of the directions of the above theorem directly follows from two results of Hrushovski in \cite{elementaryfrobenius} and \cite{Hrushovski2021AxsTW} respectively:
\begin{theorem}\label{theoremdirectionhrushovskisresults}
Any characteristic $0$ ultraproduct $\prod_{\mathcal{U}}(\mathbb{F}_{q}^{\mathrm{alg}},\mathrm{Frob}_{q},\Psi_{q})$ is a model of $\mathrm{ACFA}^{+}$. This even holds if we allow for an arbitrary additive character $\Tilde{\Psi}_{q}$ on $\mathbb{F}_{q}$.
\end{theorem}
\begin{proof}
    First we use the result of Hrushovski (Theorem 1.4 in \cite{elementaryfrobenius}, see \cite{HrushLangWeilbyShuddhodanVayarslan} for a different proof of the main result of \cite{elementaryfrobenius}) that the $\mathcal{L}_{\sigma}-$structure $\prod_{\mathcal{U}}(\mathbb{F}_{q}^{\mathrm{alg}},\mathrm{Frob}_{q})$ is a model of $\mathrm{ACFA}$.
    Then, to conclude, we apply Lemma 3.6 of \cite{Hrushovski2021AxsTW} which states that $\prod_{\mathcal{U}}(\mathbb{F}_{q},\Tilde{\Psi}_{q})$ is a model of $\mathrm{PF}^{+}$.
\end{proof}
For the other direction, we show that for every $(K,\sigma,\Psi)\models\mathrm{ACFA}^{+}$ we find some $\prod_{\mathcal{U}}(\mathbb{F}_{q}^{\mathrm{alg}},\mathrm{Frob}_{q},\Psi_{q})$ elementarily equivalent to it. Let us remark here that while the corresponding statement for $\mathrm{PF}^{+}$ was proved in \cite{Hrushovski2021AxsTW} and the one for $\mathrm{ACFA}$ in \cite{Macintyre1997-MACGAO}, this does not immediately yield a proof of Theorem \ref{maintheoremlimittheory}. This is due to the fact that, in general, the $\mathcal{L}_{\mathrm{ring}}$-theory of $F:=\mathrm{Fix}(\sigma)\subset K\models\mathrm{ACFA}$ does not uniquely determine the $\mathcal{L}_{\sigma}$-theory of $(F^{\mathrm{alg}},\sigma)$ (this follows from Fact \ref{inducedstructureonfixedfield}). However, we will see that this turns out to be only a subtlety that reduces to an application of Lemma \ref{primecomputationlemma} when following the proof of Hrushovski in \cite{Hrushovski2021AxsTW} adapted to our context. Before turning to the proof let us consider the question whether Theorem \ref{maintheoremlimittheory} still holds if we replace \textit{finite} by \textit{prime}. For $\mathrm{ACFA}$ this is not problematic; it is the $\mathcal{L}_{\sigma}$-theory of the algebraic closure of prime fields with Frobenius. However, for $\mathrm{PF}^{+}$ this is already not the case as discussed in Chapter 5 of \cite{Hrushovski2021AxsTW}, and thus neither for $\mathrm{ACFA}^{+}$. 
\begin{notation}
    From now on we write $T_{p}$ for the limit theory of the structures $(\mathbb{F}_{p},\mathrm{Frob}_{p},\Psi_{p})$ and by $T_{q}$ the limit theory of the structures $(\mathbb{F}_{q},\mathrm{Frob}_{q},\Psi_{q})$.
\end{notation}

\begin{fact}[See 5.1 in \cite{Hrushovski2021AxsTW}]\label{primefieldsstandardobservation}
Let $(F,\sigma,\Psi)\models T_{p}$. Let $n\in\mathbb{N}$ and $\mu_{n}$ be the set of $n$-th roots of unity. Let $\alpha_{n}\in(\mathbb{Z}/(n))^{\times}$ be such that $\sigma^{-1}(\omega)=\omega^{\alpha_{n}}$ for all $w\in\mu_{n}$. Then $\Psi(\frac{-1}{n})=\exp(2\pi i\frac{\alpha_{n}}{n})$.

\end{fact}

\begin{remark}
    Finally, we could ask if $T_{p}$ consists of $\mathrm{ACFA}^{+}$ together with axioms stating that $\Psi(1)=1$ and $\Psi(\frac{-1}{n})=\exp(2\pi i\frac{\alpha_{n}}{n})$ for all $n\in\mathbb{N}$. There is not much hope to answer this question easily, since a similar question in the $\mathrm{PF}^{+}$-case already seems to be out of reach. This is what Chapter 5 in \cite{Hrushovski2021AxsTW} is about. See Conjecture 5.2 in \cite{Hrushovski2021AxsTW} for the statement.
\end{remark}

The outline of the proof of Theorem \ref{maintheoremlimittheory} is as follows: By Corollary \ref{substructureeltequiv} two models of $\mathrm{ACFA}^{+}$ are elementarily equivalent if their substructures on $\mathbb{Q}^{\mathrm{alg}}$ have the same isomorphism-type. By compactness, we only have to consider finitely many finite Galois extensions of $\mathbb{Q}$ and can then reduce to work with only one such extension. Fix such an extension $L$ of degree $n$. By an application of the \v{C}ebotarev density theorem we can find infinitely many primes $p$ such that the $\mathcal{L}_{\sigma}$- sentences describing the finite Galois extension over $\mathbb{Q}$ are fulfilled in the algebraic closure of those prime fields, together with the corresponding Frobenius automorphism. Concretely, this sentence describes for $f_{n}(X)\in\mathbb{Q}[X]$ a monic irreducible polynmial of degree $n$ with splitting field $L_{n}$ the action of $\sigma$ on its roots.\\
Then, to deal with the character, i.e., the $\mathcal{L}_{\sigma,+}$ sentences that describe it in $L$ we pass from $\mathbb{F}_{p}^{\mathrm{alg}}$ to $\mathbb{F}_{p^{e_{p}}}^{\mathrm{alg}}$ carefully choosing the $e_{p}$ using that $\Psi_{p^{e_{p}}}(x)=\Psi_{p}(e_{p}x)$ for $x\in\mathbb{F}_{p}$. So far, this is the same proof as given by Hruhsovski in Proposition 3.16 in \cite{Hrushovski2021AxsTW}. In our case however, when passing from $\mathbb{F}_{p}^{\mathrm{alg}}$ to $\mathbb{F}_{p^{e_{p}}}^{\mathrm{alg}}$ we have to make sure that we do not \textit{change} the generator of the Galois group the Frobenius represents in the following sense\footnote{The necessity of that is an instance of the fact that the induced structure on a fixed field of a model of $\mathrm{ACFA}$ is more than the (pure) field structure (Fact \ref{inducedstructureonfixedfield}).}: We have to make sure that 
$\mathrm{Frob}_{p^{e_{p}}}(\alpha)=h(\alpha)=\mathrm{Frob}_{p}(\alpha)$ holds for $\alpha$ being a root of $f_{n}$ and $h\in\mathbb{Q}[X]$. To show that this is possible will, however, reduce to an elementary number theoretic statement given in Lemma \ref{primecomputationlemma} that is already present in Lemma 3.17 in \cite{Hrushovski2021AxsTW}.
We start with some preparational work, following \cite{Hrushovski2021AxsTW}, that will later allow us to choose the $e_{p}$.

\begin{fact}\label{torsorsubgroup}
 Let $G\leq \mathbb{T}^{n}$ be a proper closed subgroup. Then there is some $\bar{m}=(m_{1},\dots,m_{n})\in\mathbb{Z}^{n}\backslash 0$ such that for any $\bar{g}=(g_{1},\dots,g_{n})\in G$ we have $\prod_{1\leq i \leq n}g_{i}^{m_{i}}=1$. And moreover, we can choose finitely many such tuples $(m_{1,j},\dots,m_{n,j})_{j}$ such that $G$ is given by all those $(g_{1},\dots,g_{n})$ such that $\prod_{1\leq i \leq n}g_{i}^{m_{i,j}}=1$ for all $1\leq j\leq k$.    
\end{fact}
\begin{proof}
    This can be shown using Pontryagin-duality: The character group, which is $Hom(\mathbb{T}^{n},S^{1})$ (in the category of topological groups) is isomorphic to $(\mathbb{Z}^{n},+)$. Let $\widehat{G}$ be the dual of $G$, i.e., $\widehat{G}$ is given by those $\phi\in Hom(\mathbb{T}^{n},S^{1})$ that annihilate $G$.
    Now (see Theorem 33, Chapter 5, in \cite{pontryagin1966topological}) we have that $G$ is its double dual, i.e., $G=\skew{2.3}\widehat{\widehat{G}}$ and the Fact directly follows as subgroups of $(\mathbb{Z}^{n},+)$ are finitely generated.
\end{proof}

\begin{fact}[Lemma 3.18 in \cite{Hrushovski2021AxsTW}]\label{fullimageofPsiFc}
Let $(F,\Psi)$ be an ultraproduct of finite fields with arbitrary non-trivial character. Let $n\in\mathbb{N}$, $\bar{c}=(c_{1},\dots,c_{n})\in F^{n}$ and write $F\bar{c}$ for $\{(fc_{1},\dots,fc_{n})\,|\,f\in F\}$.
Assume that for any $\bar{m}=(m_{1},\dots,m_{n})\in\mathbb{Z}^{n}\backslash 0$ we have $\sum_{1\leq i\leq n}m_{i}c_{i}\neq 0$, then $\Psi^{(n)}(F\bar{c})=\mathbb{T}^{n}$.
\end{fact}

\begin{fact}(See Lemma 3.17 in \cite{Hrushovski2021AxsTW})\label{finitefieldextensionascorrection}
Let $(p_{m})_{m\in\mathbb{N}}$ be a sequence of strictly increasing prime numbers. Let $\bar{c}_{m}=(c_{m_{1}},\dots,c_{m_{n}})\in\mathbb{F}_{p_{m}}^{n}$ be a sequence of $n$-tuples. Assume that for any given nonzero integer vector $(\alpha_{1},\dots,\alpha_{n})$ for all but finitely many $m$ it is not the case that $\sum_{i=1}^{n}\alpha_{i}c_{m_{i}}=0.$ Let $U$ be a non-empty open subset of $\mathbb{T}^{n}$.
Then for all but finitely many $m$ we find some $e_{m}\in\mathbb{N}$ such that $\Psi_{p_{m}}(e_{m}\bar{c}_{m})\in U$.
\end{fact}


Now we give a description of the compactness argument we use to reduce to finite Galois extensions. The idea is the same as in \cite{Hrushovski2021AxsTW}, adapted to our context and extended by details. From now on, we fix the following notation.
\begin{notation}\label{notationenrichedgaloissentences}
     Let $(K,\sigma,\Psi)\models\mathrm{ACFA}^{+}$ and write $\mathbb{Q}^{\mathrm{alg}}$ as an increasing union of finite Galois-extensions $(L_{n})_{n\in\mathbb{N}}$ of $\mathbb{Q}$. For each $n$ denote by $(L_{n},\sigma,\Psi)$ the induced structure on $L_{n}$. Let $f_{n}(X)\in\mathbb{Q}[X]$ be some monic irreducible polynomial of degree $n$ such that $L_{n}$ is its splitting field and $[L_{n}:\mathbb{Q}]=n$. Now we take some root $\beta$ of $f_{n}(X)$ and write the other roots of $f_{n}(X)$ as polynomials in $\beta$ over $\mathbb{Q}$. I.e., we have that $\{h_{1}(\beta),\dots,h_{n}(\beta)\}$ is the set of roots of $f_{n}$, where $h_{i}(X)\in\mathbb{Q}[X]$ and $h_{1}(X)=X$. We fix our choice of $\beta$ and of the $h_{i}(X)$. Note that $\text{deg} h_{i}\leq n-1$ for $n\geq 2$.
\end{notation}
\begin{fact}\label{factisomtypesentenceclassic}(Introduction of Section 1.6 in \cite{Macintyre1997-MACGAO}, see also 1.6 of \cite{acfa}.)
For $(L_{n},\sigma)$ seen as an $\mathcal{L}_{\sigma} $ -structure we can find an $\mathcal{L}_{\sigma}$- (and thereby also $\mathcal{L}_{\sigma,+}$)-sentence $\phi_{L_{n},c}$ that expresses $\sigma\in c$ where $c$ denotes the conjugacy class of $\sigma$ in $\mathrm{Gal}(L_{n}/\mathbb{Q})$. Concretely $\phi_{L_{n},c}$ is of the form
\[\exists\beta\left(f_{n}(\beta)=0\;\land\;\sigma(\beta)=h_{i}(\beta)\right)\]
for some $1\leq i\leq n$. 
\end{fact}

Now, our aim will be to consider $\Theta_{n}$, a countable set of $\mathcal{L}_{\sigma,+}$-sentences that describes $(L_{n},\sigma,\Psi)$ up to isomorphism.  We continue to work in the notation of Notation \ref{notationenrichedgaloissentences}.

\begin{notation}
    Let $\{\alpha_{1},\dots,\alpha_{m}\}$ denote a $\mathbb{Q}$-basis of $\mathrm{Fix}(\sigma)$ in $L_{n}$. We can write $\alpha_{j}=g_{j}(\beta)$ for some $g_{j}\in\mathbb{Q}[X]$ with $\text{deg}\, g_{j}\leq n-1$ for every $1\leq j\leq m$. We fix the $g_{1}(X),\dots,g_{m}(X)$ for the rest of the section.
\end{notation}

\begin{definition}\label{definitionthetasentence}
For any $s,l\in\mathbb{N}$ we denote by $\Theta_{n,s}^{l}$ the $\mathcal{L}_{\sigma,+}$-sentence
    \[\exists\beta\left(f_{n}(\beta)=0\;\land\;\sigma(\beta)=h_{i}(\beta)\;\land\;\bigwedge_{1\leq j\leq m}\left|\Psi\left(\frac{1}{s}g_{j}(\beta)\right)-r_{j,s}\right|\leq\frac{1}{l}\right)\;\;\;\;(\star)
    \]where the $r_{j,s}$ are given by the value $\Psi(\frac{1}{s}\alpha_{j})$ for $1\leq j\leq m$ in the theory of $(L_{n},\sigma,\Psi)$. Now we define the set of $\mathcal{L}_{+,\sigma}$-sentences $\Theta_{n}$ as $\Theta_{n}:=\{\Theta_{n,s}^{l}\,|\,s,l\in\mathbb{N}\}$.    
\end{definition}

\begin{remark}
    There is no problem in working with an existential quantifier instead of an infimum in the sentences $\Theta_{n,s}^{l}$. Indeed, modulo $\mathrm{ACFA}^{+}$, the sentences $\Theta_{n,s}^{l}$ are expressible in $\mathcal{L}_{\sigma,+}$ since the quantification is over a finite (strongly) definable set and thus infimum and existential quantifier coincide.
\end{remark}
\begin{lemma}\label{lemmaenrichedfinitegaloisisomorphismtype}
    The set of sentences $\Theta_{n}$ describes the $\mathcal{L}_{\sigma,+}$-isomorphism type of $(L_{n},\sigma,\Psi)$.
\end{lemma}
\begin{proof}
    Let $(L_{n}^{\prime},\sigma^{\prime},\Psi^{\prime})$ be such that $L_{n}^{\prime}$ is also given as the splitting field of $f_{n}(X)$ over $\mathbb{Q}$. Let $(L_{n}^{\prime},\sigma^{\prime},\Psi^{\prime})\models \Theta_{n}$. We want to show that
    \[(L_{n},\sigma,\Psi)\cong_{\mathcal{L}_{\sigma,+}}(L_{n}^{\prime},\sigma^{\prime},\Psi^{\prime}).\]
    Note that in general for $x\in L_{n}$ we have that $\Psi(\frac{1}{s^{\prime}}x)$ determines $\Psi(\frac{1}{s}x)$ where $s$ divides $s^{\prime}$. Then, by the pigeonhole principle there exists a root $\beta^{\prime}$ of $f_{n}(X)$ in $L_{n}^{\prime}$ that witnesses $\Theta_{n,s}^{l}$ for all $s\in\mathbb{N}$ and $l\in\mathbb{N}$ at the same time. Now there is an $\mathcal{L}_{\sigma}$-isomorphism $\tau$ between $(L_{n},\sigma)$ and $(L_{n}^{\prime},\sigma^{\prime})$ sending $\beta$ to $\beta^{\prime}$. Now $\tau$ also preserves $\Psi$:
    Let $\alpha_{j}^{\prime}\in L_{n}^{\prime}$ be given by $g_{j}(\beta^{\prime})$ for all $1\leq j\leq m$, then $\{\alpha_{1}^{\prime},\dots,\alpha_{m}^{\prime}\}$ is a $\mathbb{Q}-$basis of $\mathrm{Fix}(\sigma^{\prime})$ in $L_{n}^{\prime}$ and $\tau(\alpha_{j})=\alpha_{j}^{\prime}$ for all $1\leq j\leq m$. But then $\tau$ clearly preserves $\Psi$ as $\beta^{\prime}$ witnesses $\Theta_{n}$.
\end{proof}

\begin{lemma}\label{lemmareductiontotonesentence}
 Let $\Theta_{n_{1},s_{1}}^{l_{1}},\dots,\Theta_{n_{k},s_{k}}^{l_{k}}$ be as in Definition \ref{definitionthetasentence} with corresponding monic irreducible polynomials $f_{n_{1}}(X),\dots,f_{n_{k}}(X)\in\mathbb{Q}[X]$. We can find an $\mathcal{L}_{\sigma,+}$ - sentence $\Theta_{n,s}^{l}$ of the form $(\star)$ as in Definition \ref{definitionthetasentence} such that $\Theta_{n,s}^{l}$ implies all the $\Theta_{n_{1},s_{1}}^{l_{1}},\dots,\Theta_{n_{k},s_{k}}^{l_{k}}$.  
\end{lemma}
\begin{proof}
Let $L$ be the Galois-extension of $\mathbb{Q}$ obtained as the splitting field of the $f_{n_{1}}(X),\dots,f_{n_{k}}(X)$. Next, we simply treat all the $\Theta_{n_{1},s_{1}}^{l_{1}},\dots,\Theta_{n_{k},s_{k}}^{l_{k}}$ in the Galois extension $L$. Concretely, let $f(X)$ be a monic irreducible  polynomial of degree $n$ over $\mathbb{Q}$ such that $L$ is its splitting field and $[L:\mathbb{Q}]=n$. As before, we fix a root $\beta$ of $f(X)$ and list all its conjugates as $h_{1}(\beta),\dots,h_{n}(\beta)$ with $h_{i}(X)\in\mathbb{Q}[X],\;h_{1}(X)=X$. Now there is a unique $1\leq i\leq n$ such that $\sigma(\beta)=h_{i}(\beta)$. From this we can compute the action of $\sigma$ on all the roots of all the $f_{n_{t}}(X)$ for $1\leq t\leq k$. 
Next, we consider some $\mathbb{Q}$-basis $\{\gamma_{1},\dots,\gamma_{m}\}$ of $\mathrm{Fix}(\sigma)$. Then we can compute all the elements from the fixed fields of the splitting fields of the $f_{n_{t}}(X)$ that occur in $\Theta_{n_{1},s_{1}}^{l_{1}},\dots,\Theta_{n_{k},s_{k}}^{l_{k}}$ using $\mathbb{Q}$-linear combinations of the $\{\gamma_{1},\dots,\gamma_{m}\}$. By replacing $\gamma_{j}$ by $\frac{1}{v_{j}}\gamma_{j}$ for some $v_{j}\in\mathbb{N}$ if necessary, we can even make sure to work with $\mathbb{Z}$-linear combinations. Or, in other words, this way we can make sure that the value of the $\Psi(\gamma_{j})$ determines the value of $\Psi$ on all the elements from the fixed fields of the splitting fields of the $f_{n_{t}}(X)$ that occur in $\Theta_{n_{1},s_{1}}^{l_{1}},\dots,\Theta_{n_{k},s_{k}}^{l_{k}}$. Thus, if we choose the $r_{j,s}$ appropriately and $l$ sufficiently big we can make sure that $\left|\Psi\left(\frac{1}{s}\gamma_{j}(\beta)\right)-r_{j,s}\right|\leq\frac{1}{l}$ implies that $\Psi$ takes its values as described by $\Theta_{n_{1},s_{1}}^{l_{1}},\dots,\Theta_{n_{k},s_{k}}^{l_{k}}$.
Now write all the $\gamma_{j}$ as $g_{j}(\beta)$ for some polynomial $g_{j}(X)\in\mathbb{Q}[X]$. Consequently, it suffices to work with a sentence of the form $(\star)$ with the polynomials $f,h_{i},g_{1},\dots,g_{m}$.
\end{proof}
After having completed the reduction to one sentence of the form $(\star)$, 
before we show that for any sentence of the form $(\star)$ there can be found structures $(\mathbb{F}_{q}^{\mathrm{alg}},\mathrm{Frob}_{q},\Psi_{q})$ of arbitrarily large characteristic such that the sentence holds, we will give the calculation on prime numbers that we need in order to use Hrushovski's strategy from \cite{Hrushovski2021AxsTW} in our context.

\begin{lemma}\label{primecomputationlemma}
    Given some $e_{p}^{\prime}\in\mathbb{N}$ and some prime $p$ coprime to $n$ and coprime to $\phi(n)$ (Euler-totient function) we can find some $e_{p}$ of the form $e_{p}=e_{p}^{\prime}+mp$ where $m\in\mathbb{N}$ such that the following two statements hold for $q=p^{e_{p}}$ and $q^{\prime}=p^{e_{p}^{\prime}}$:
    \begin{itemize}
        \item $p^{e_{p}}\equiv p\mod{n}.$
        \item $\Psi_{q^{\prime}}(x)=\Psi_{q}(x)$ for $x\in\mathbb{F}_{q}$.
    \end{itemize}
   
\end{lemma}
\begin{proof}
    Let $G$ be the subgroup of $(\mathbb{Z}/(n))^{\times}$ generated by $p\mod{n}$, then $G$ is cyclic and $|G|=r$ divides $\phi(n)$, hence $r$ and $p$ are coprime. But from this we can directly deduce that $p^{p}\mod{n}$ also yields a generator of $G$. Now we can find some $m\in\mathbb{N}$ such that the first assertion holds, simply by choosing $m$ such that $p^{e_{p}^{\prime}}\cdot (p^{p})^{m}=p\mod{n}$. For the second assertion, note that
    \[\Psi_{q}(x)=\Psi_{p}(e_{p}x)=\Psi_{p}((e_{p}^{\prime}+mp)x)=\Psi_{p}(e_{p}^{\prime}x)=\Psi_{q^{\prime}}(x).\]
\end{proof}
\begin{lemma}\label{mainlemmafiniteproof}
Let $\Theta_{n,s}^{l}$ be an $\mathcal{L}_{\sigma,+}$-sentence of the form $(\star)$ as obtained in Lemma \ref{lemmareductiontotonesentence}, i.e., with polynomials $f,h_{i},g_{1},\dots,g_{m}$. Let $L$ be the Galois extension of degree $n$ over $\mathbb{Q}$ obtained as the splitting field of $f$. Assume that for some $\sigma\in \mathrm{Gal}(L/\mathbb{Q})$ and group homomorphism $\Psi:(\mathrm{Fix}(\sigma),+)\rightarrow (S^{1},\cdot)$ we have \[(L,\sigma,\Psi)\models \Theta_{n,s}^{l}.\]
Then there are infinitely many primes $p$ with corresponding $e_{p}\in\mathbb{N}$ such that for the character $\Tilde{\Psi}_{p}$ on $\mathbb{F}_{p}$ defined by $\Tilde{\Psi}_{p}(x):=\Psi_{p}(e_{p}x)$ we have 
    \[(\mathbb{F}_{p}^{\mathrm{alg}},\mathrm{Frob}_{p},\Tilde{\Psi}_{p})\models \Theta_{n,s}^{l}\]
    Moreover we can choose the $e_{p}\in\mathbb{N}$ such that for $q:=p^{e_{p}}$ we have
    \[(\mathbb{F}_{q}^{\mathrm{alg}},\mathrm{Frob}_{q},\Psi_{q})\models \Theta_{n,s}^{l}.\]
\end{lemma}
\begin{proof}
As before, we write $\beta$ for a root of $f(X)$ and let $1\leq i\leq n$ be such that $\sigma(\beta)=h_{i}(\beta)$ holds. Further, let $\gamma_{j}=g_{j}(\beta)$ for all $1\leq j\leq m$. Then $\{\gamma_{1},\dots,\gamma_{m}\}$ is a $\mathbb{Q}$-basis of $\mathrm{Fix}(\sigma)$ in $L$.\\
We apply the \v{C}ebotarev density Theorem as in 1.6 of \cite{Macintyre1997-MACGAO} (see also 1.14 of \cite{acfa}) to obtain a sequence of primes $(p_{k})_{k\in\mathbb{N}}$ such that the $(\mathbb{F}_{p_{k}}^{\mathrm{alg}},\mathrm{Frob}_{p_{k}})$ satisfy the $\mathcal{L}_{\sigma}$-sentence describing the isomorphism type of $(L,\sigma)$ as in Fact \ref{factisomtypesentenceclassic}. So, it follows that $\mathrm{Frob}_{p_{k}}(\beta^{\prime})=h_{i}(\beta^{\prime})$ for a root $\beta^{\prime}$ of $f(X)$ in $\mathbb{F}_{p_{k}}^{\mathrm{alg}}$. At this point, we already discard the finitely many primes dividing $s$ or a denominator of coefficients in $f(X),h_{i}(X)$ or some $g_{j}(X)$. Moreover, in any such $(\mathbb{F}_{p_{k}}^{\mathrm{alg}},\mathrm{Frob}_{p_{k}})$ there is a set $\{\gamma_{1}^{\prime},\dots,\gamma_{m}^{\prime}\}$ inside $\mathbb{F}_{p_{k}}$ that corresponds to the set $\{\gamma_{1},\dots,\gamma_{m}\}$. I.e., $\gamma_{j}^{\prime}$ is of the form $g_{j}(\beta^{\prime})$ in $\mathbb{F}_{p_{k}}$ where $\gamma_{j}=g_{j}(\beta)$ in $\mathrm{Fix}(\sigma)\subseteq L$.
Now we have to modify the character in a way such that the prescribed conditions hold.\\ Let $U$ be a nonempty open set of $\mathbb{T}^{m}$ such that $\Psi^{(m)}(\gamma_{1},\dots,\gamma_{m})\in U$ would imply that $\Theta_{n,s}^{l}$ holds.
At this point, we apply Fact \ref{finitefieldextensionascorrection} to find a sequence of natural numbers $(e_{k})_{k\in\mathbb{N}}$ such that $\Psi^{(m)}_{p_{k}}(e_{k}\gamma_{1}^{\prime},\dots,e_{k}\gamma_{m}^{\prime})\in U$. Now we are done proving the first part of the theorem (prime fields, arbitrary character).\\
To prove the second part, we will pass from the structure $(\mathbb{F}_{p_{k}}^{\mathrm{alg}},\mathrm{Frob}_{p_{k}},\Psi_{p_{k}})$ to $(\mathbb{F}_{q_{k}}^{\mathrm{alg}},\mathrm{Frob}_{q_{k}},\Psi_{q_{k}})$, where $q_{k}=p_{k}^{e_{k}}$. We have to take care that we do not change the Galois generator that is given by the Frobenius, i.e., that we do not change the $i$ such that $\mathrm{Frob}_{p_{k}}(\beta^{\prime})=h_{i}(\beta^{\prime})$.
Note, that we can assume that $p_{k}$ and $n$ (the degree of the Galois extension $L$) are coprime. By applying Lemma \ref{primecomputationlemma} we can replace $e_{k}$ by some $\Tilde{e}_{k}$ (of the form $\Tilde{e}_{k}=e_{k}+N_{k}\cdot p_{k}$ for $N_{k}\in\mathbb{N}$) such that the two following conditions hold \begin{itemize}
    \item $\Psi_{p_{k}^{\Tilde{e}_{k}}}(\gamma_{j}^{\prime})=\Psi_{p_{k}^{e_{k}}}(\gamma_{j}^{\prime})\;\;\forall\,1\leq j\leq m.$
    \item $p_{k}^{\Tilde{e}_{k}}\mod{n}=p_{k}\mod{n}.$
\end{itemize}
Then, we have for $\Tilde{q}_{k}:=p_{k}^{\Tilde{e}_{k}}$ that $(\mathbb{F}_{\Tilde{q}_{k}}^{\mathrm{alg}},\mathrm{Frob}_{\Tilde{q}_{k}})$ satisfies the same $\mathcal{L}_{\sigma}$-sentence describing the Galois-extension of degree $n$ over $\mathbb{Q}$, i.e., the isomorphism type of $(L,\sigma)$. Consequently, it follows that $\Theta_{n,s}^{l}$ holds in $(\mathbb{F}_{\Tilde{q}_{k}}^{\mathrm{alg}},\mathrm{Frob}_{\Tilde{q}_{k}},\Psi_{\Tilde{q}_{k}})$, which completes the proof.
\end{proof}

Now we can finally proceed to prove the main theorem of this section: 
\begin{proof} (of Theorem \ref{maintheoremlimittheory})
Fix some $(K,\sigma,\Psi)\models\mathrm{ACFA}^{+}$ and $(L_{n},\sigma,\Psi)$ as in Notation \ref{notationenrichedgaloissentences}. Let $(H,\Tilde{\sigma},\Tilde{\Psi})\models\mathrm{ACFA}^{+}$ and denote the induced structure of $(H,\Tilde{\sigma},\Tilde{\Psi})$ on $L_{n}$ by $(H_{n},\Tilde{\sigma},\Tilde{\Psi})$. Let $\Theta_{n}$ be as in Definition \ref{definitionthetasentence}.
Then it follows by Lemma \ref{lemmaenrichedfinitegaloisisomorphismtype} that
\[(H,\Tilde{\sigma},\Tilde{\Psi})\models \Theta_{n}\;\iff\;(H_{n},\Tilde{\sigma},\Tilde{\Psi})\cong(L_{n},\sigma,\Psi).\]
Now simply by applying compactness and Lemma \ref{lemmareductiontotonesentence} it follows that there is some ultraproduct of the form $\prod_{\mathcal{U}}(\mathbb{F}_{q}^{\mathrm{alg}},\mathrm{Frob}_{q},\Psi_{q})$ with  \[(K,\sigma,\Psi)\equiv\prod_{\mathcal{U}}(\mathbb{F}_{q}^{\mathrm{alg}},\mathrm{Frob}_{q},\Psi_{q})\] if and only if for every sentence $\Theta_{n,s}^{l}$ of the form $(\star)$ as in Definition \ref{definitionthetasentence} we have that
\[(\mathbb{F}_{q}^{\mathrm{alg}},\mathrm{Frob}_{q},\Psi_{q})\models\Theta_{n,s}^{l}\;\]for infinitely many $q$. This works similarly for prime fields with arbitrary character. Now, the statement follows directly from Lemma \ref{mainlemmafiniteproof} and from Theorem \ref{theoremdirectionhrushovskisresults}.
\end{proof}

\section{3-amalgamation over $\sigma$-AS-closed sets}\label{section3amalgam}
We will now deal with the question when $3$-amalgamation holds over an algebraically closed set $E\subseteq\mathcal{M}\models \mathrm{ACFA}^{+}$. Here we use algebraically closed in the sense of the theory $\mathrm{ACFA}$. When we pass to the continuous logic algebraic closure, including CL-imaginaries, we will make this explicit by writing $\mathrm{acl}_{\mathrm{CL}}^{\mathrm{eq}}$. Let us first recall that in the theories $\mathrm{ACFA}$ and $\mathrm{PF}^{+}$ $3$-amalgamation (and even $n$-amalgamation for every $n\in\mathbb{N}$) holds over all algebraically closed sets. However, we will see that in $\mathrm{ACFA}^{+}$ this is not the case in general. Nevertheless, we will show how to characterise those algebraically closed $E$ such that $3$-amalgamation holds over $E$. This will be the main result that we will later apply to obtain a description of the connected component (which can be nontrivial) of the Kim-Pillay group for the completions of $\mathrm{ACFA}^{+}$ in Section \ref{sectionautomgrouptorsors} and a natural extension of $\mathrm{ACFA}^{+}$ that (weakly) eliminates CL-imaginaries in Section \ref{sectioneliminationofimaginaries}.\\ Moreover, we will see, as 3-amalgamation does hold over models, that $\mathrm{ACFA}^{+}$ is still a simple theory. In the context of classical discrete logic, it is a longstanding open question if any simple theory has $3$-amalgamation over all $\mathrm{acl}^{\mathrm{eq}}$-closed sets and equivalently if the Kim-Pillay group is totally disconnected. However, note that since $\mathrm{ACFA}^{+}$ is a continuous logic theory, it is not an answer to this question. Also, it was already observed by Hrushovski in \cite{Hrushovski1998SimplicityAT} that outside of the classical first-order context (in his case in the context of Robinson theories) simplicity does not necessarily imply that the Kim-Pillay group is totally disconnected.

\begin{notation}
    From now on we will denote by $\mathrm{acl}_{\sigma}(A)$ the field-theoretic algebraic closure of $\mathrm{cl}_{\sigma}(A)$, the difference field generated by $A$, where $A$ is a subset of some model of $\mathrm{ACFA}$ (or $\mathrm{ACFA}^{+}$). This coincides with the model-theoretic algebraic closure in $\mathrm{ACFA}$. (See (1.7) in \cite{acfa}.)
    We will write $p_{\mathcal{L}_{\sigma}}(\bar{x})$ for $p(\bar{x})\restriction_{\mathcal{L}_{\sigma}}$ where $p(\bar{x})$ is some $\mathcal{L}_{\sigma}^{+}$-type.
\end{notation}

\begin{definition}\label{definitionindependenceacfa}
  For subsets $A,B,C\subseteq M$ with $C\subseteq A,B$ we say that $A\ind_{C} B$ ($A$ is independent from $B$ over $C$) if and only if $\mathrm{acl}_{\sigma}(A)$ is algebraically independent from $\mathrm{acl}_{\sigma}(B)$ over $\mathrm{acl}_{\sigma}(C)$.
\end{definition}
\begin{fact}((1.9) in \cite{acfa}.)
    The $\mathcal{L}_{\sigma}$-theory $\mathrm{ACFA}$ is simple and forking independence coincides with the above notion of independence.
\end{fact}

We start with a small but important lemma that will later allow us to extend $\Psi$ freely whenever it does not interfere with any $\mathbb{Q}$-linear dependencies. 

\begin{lemma}\label{collaryextensionlinearindtypes}
 Let $K\models \mathrm{ACFA}^{+}$ and $E=\mathrm{acl}_{\sigma}(E)\subseteq K$. Let $\bar{c},\bar{d}$ be tuples in $K\backslash E$ and, moreover, $\bar{c}\in F^{n}$. Assume that $\bar{c}$ is not contained in any rational hyperplane over $F\cap \mathrm{acl}_{\sigma}(E\bar{d})$. Let $p(\bar{x},\bar{y}):=\mathrm{tp}_{\mathcal{L}_{\sigma}}(\bar{c}\bar{d}/E)\cup \mathrm{tp}(\bar{d}/E)$.
 For any $\bar{k}\in\mathbb{N}^{n}$ and any $\bar{r}\in\mathbb{T}^{n}$, the partial $\mathcal{L}_{\sigma,+}$-type $p(\bar{x},\bar{y})$ is consistent with $\Psi^{(n)}(\frac{1}{k_{1}}x_{1},\dots,\frac{1}{k_{n}}x_{n})=\bar{r}$.
\end{lemma}
\begin{proof}
    Take some $\mathcal{L}_{\sigma}$-structure $H\supset E$ such that $\mathrm{acl}_{\sigma}(H)=H$ and $H$ contains a realisation of $\mathrm{tp}_{\mathcal{L}_{\sigma}}(\bar{c}\bar{d}/E)$. Then, since $\bar{c}$ is $\mathbb{Q}$-linear independent from $\mathrm{acl}_{\sigma}(E\bar{d})$, we can extend $\Psi$ to an homomorphism on $H\cap F$ respecting $\mathrm{tp}(\bar{d}/E)$ as well as $\Psi^{(n)}(\frac{1}{k_{1}}x_{1},\dots,\frac{1}{k_{n}}x_{n})=\bar{r}$. By model completeness (Corollary \ref{corollarymodelcompleteness}) of $\mathrm{ACFA}^{+}$ this structure embeds into a model of $\mathrm{ACFA}^{+}$. This completes the proof.
\end{proof}

We use Theorem 1.51 in \cite{simplicityCATs} which generalises to the context of CATs (which comprises continuous logic) the main theorem by Kim and Pillay in
\cite{Kim199simpletheories}, the characterisation of simplicity by the presence of a good notion of independence (which most importantly satisfies the Independence Theorem over models).
It follows directly from the corresponding properties for $\mathrm{ACFA}$ that properties 1.-6. of Theorem 1.51 in \cite{simplicityCATs} for the above notion of independence still hold in $\mathrm{ACFA}^{+}$. We will now observe that property 7 (Extension) holds and then later prove the Independence Theorem over models for $\mathrm{ACFA}^{+}$ at the end of this section (Corollary \ref{corollaryindependencetheoremsigmaASstrongtypes}) which using the above mentioned Theorem 1.51 in \cite{simplicityCATs} suffices to prove simplicity.

\begin{lemma}\label{lemmaextnesionindependenceinacfaplus}We work in the notation of Definition \ref{definitionindependenceacfa}. If $A\ind_{C}B$ and $D\supseteq B$, then we find $A^{\prime}$ with $\mathrm{tp}(A^{\prime}/C)=\mathrm{tp}(A/C)$ and $A^{\prime}\ind_{C}D$.
\end{lemma}
\begin{proof}
   We use Extension in $\mathrm{ACFA}$ for the corresponding $\mathcal{L}_{\sigma}$-reducts and observe that we can extend $\Psi$ to $A^{\prime}$ consistent with $\mathrm{tp}(A/C)$ by Lemma \ref{collaryextensionlinearindtypes}. 
\end{proof}
Let us briefly recall what is meant by 3-amalgamation over an (algebraically closed) set.

\begin{definition}
     We work in some theory $T$ with $A\subseteq \mathcal{M}\models T$ and a notion of independence satisfying properties 1.-7. of Theorem 1.51 in \cite{simplicityCATs}. A 3-amalgamation problem over $A$ is a system of complete $A$-types $\mathcal{V}=\{p_{1}(\bar{x}_{1}),p_{2}(\bar{x}_{2}),p_{3}(\bar{x}_{3}),p_{12}(\bar{x}_{12}),p_{13}(\bar{x}_{13}),p_{23}(\bar{x}_{23})\}$ such that
    \begin{itemize}
        \item for any $\bar{a}_{w}\models p_{w}$ the infinite tuple $\bar{a}_{w}$ enumerates its algebraic closure. 
        \item For any $w\subseteq w^{\prime}$ we have $\bar{x}_{w}\subseteq \bar{x}_{w^{\prime}}$ and $p_{w^{\prime}}\restriction_{w}=p_{w}$.
        \item For any $\bar{a}_{ij}\models p_{ij}$ and corresponding subtuples $\bar{a}_{i}$ and $\bar{a}_{j}$, we have $\bar{a}_{i}\ind_{\mathcal{M}}\bar{a}_{j}$.
    \end{itemize}
    We say that $\mathcal{V}$ has a solution, if there is some $p_{123}(\bar{x}_{123})$ such that the system $(p_{w})_{w\subseteq \{1,2,3\}}$ still satisfies the above properties (in particular for $\bar{a}_{123}\models p_{123}$, we require $\bar{a}_{ij}\ind_{\bar{a}_{j}}\bar{a}_{jk})$ for $\{i,j,k\}=\{1,2,3\}$. We say that $T$ has 3-amalgamation over $A$, if every 3-amalgamation problem over $A$ has a solution.
\end{definition}

We do not consider questions on higher amalgamation (i.e., $n$-amalgamation for $n\geq 4$) in $\mathrm{ACFA}^{+}$ as this will be the content of a subsequent paper but we briefly recall the two results below which in particular entail that 3-amalgamation holds over all algebraically closed sets in $\mathrm{ACFA}$ and $\mathrm{PF}^{+}$.

\begin{fact}\label{factnamalgamationacfa}(See (1.9) in \cite{acfa})
  In $\mathrm{ACFA}$ $n$-amalgamation holds for all $n$ over all algebraically closed sets.  
\end{fact}

\begin{fact}\label{factnamalgamationpfplus}(Proposition 3.20 in \cite{Hrushovski2021AxsTW}.)
    In $\mathrm{PF}^{+}$ $n$-amalgamation holds for all $n$ over all algebraically closed sets.
\end{fact}

\begin{definition}
We write $\mathcal{P}(\mathbf{n})$ for $\mathcal{P}(\{1,\dots,n\})$.
We work in $\mathrm{ACF}, \mathrm{ACFA}$ or $\mathrm{ACFA}^{+}$. We call a system of tuples $(\bar{a}_{w})_{w\in\mathcal{P}(\mathbf{n})}$ (with $\bar{a}_{\emptyset}=A$) an \textit{independent n-system over $A$} if we have for any $w,w^{\prime}\in \mathcal{P}(\mathbf{n})$ that $\bar{a}_{w}$ is independent from $ \bar{a}_{w^{\prime}}$ over $\bar{a}_{w\cap w^{\prime}}$ and, moreover, that $\bar{a}_{w}$ enumerates $\mathrm{acl}(E,\bar{a}_{i}\,|\,i\in w)$.
\end{definition}

We will start by proving the existence of algebraically closed sets $A\subseteq M$ and a $3$-amalgamation problem over $A$ that does not have a solution. We will introduce the notion of when an algebraically closed set is \textit{$\sigma$-AS-closed}. Later we will see that those sets are indeed exactly the algebraically closed sets over which $3$-amalgamation holds in general.

\begin{definition}\label{torsordefinition}
Given $a\in K$ we denote by $\mathfrak{T}_{a}$ the additive $F$-Torsor defined by $\sigma(x)-x=a$.
\end{definition}
\begin{definition}\label{definitiontorsorclosed}
We define a set $A$ as \textit{$\sigma$-AS-closed} (for Artin-Schreier) if $A=\mathrm{acl}_{\sigma}(A)$ and, moreover, if for every $a\in A$ there exists $b\in A$ such that $b\in \mathfrak{T}_{a}$.
\end{definition}
The following lemma provides an easy way to construct examples of algebraically closed sets that are not $\sigma$-AS-closed. It will also be a key ingredient in a subsequent paper on higher amalgamation in $\mathrm{ACFA}^{+}$.
\begin{lemma}\label{lemmatorsornotrealisedinalgclosure}
Let $(K,\sigma)$ be a difference field of characteristic $0$. Assume that $\mathfrak{T}_{a}$ is not realised in $K$ where $a\in K$, then $\mathfrak{T}_{a}$ is also not realised in $K^{\mathrm{alg}}$.    
\end{lemma}
\begin{proof}
    Assume the contrary, that is, $\sigma(b)-b=a$ for some $b\in K^{\mathrm{alg}}\backslash K$. Let $P(X)=X^{n}+c_{1}X^{n-1}+\dots+c_{n}$ be the minimal polynomial of $b$ over $K$. Then the minimal polynomial of $\sigma(b)$ is given by $P_{\sigma}(X)=X^{n}+\sigma(c_{1})X^{n-1}+\dots+\sigma(c_{n})$. Next, since $\sigma(b)=a+b$ we have $P_{\sigma}(a+b)=0$. We write $P_{\sigma}(a+X)$ as a polynomial in $X$, denoted by $\Tilde{P}(X)=X^{n}+c_{1}^{\prime}X^{n-1}+\dots+c_{n}^{\prime}$.
    Then, we have $c_{1}^{\prime}=na+\sigma(c_{1})$. Since $\Tilde{P}(b)=0$ and $\text{deg}_{X}\Tilde{P}= n$ we have by uniqueness of the minimal polynomial that $c_{1}^{\prime}=c_{1}$. But from $c_{1}=na+\sigma(c_{1})$ it follows that
    \[\sigma\left(-\frac{c_{1}}{n}\right)+\frac{c_{1}}{n}=a\]
    which contradicts that $\mathfrak{T}_{a}$ is not realised in $K$.
\end{proof}

\begin{corollary}
    
For every completion $T$ of $\mathrm{ACFA}^{+}$ and every $\mathcal{M}\models T$ the set $\mathbb{Q}^{\mathrm{alg}}$ is not $\sigma$-AS-closed.
\end{corollary}

\begin{proof}
    This follows directly from Lemma \ref{lemmatorsornotrealisedinalgclosure} as $\sigma$ acts as identity on $\mathbb{Q}$ and consequently, $\mathfrak{T}_{a}$ is not realised in $\mathbb{Q}^{\mathrm{alg}}$ for any $a\in \mathbb{Q}\backslash\{0\}$.
\end{proof}

\begin{definition}\label{definitionequivrelationEa}
    Fix $K\models \mathrm{ACFA}^{+}$ and $a\in K$. Let $E_{a}$ be the equivalence relation on the $(F,+)$-Torsor $\mathfrak{T}_{a}$ given by
    \[x\,E_{a}\,y\;\iff\;\Psi(x-y)=1.\]
    Depending on the context we extend $E_{a}$ by a single class outside of $\mathfrak{T}_{a}$. Then we call any $\mathbf{b}_{a}\in \mathfrak{T}_{a}/E_{a}$ (or rather the corresponding CL-imaginary as explained in Remark \ref{imaginaryasCLimaginary} below) a $\sigma$-AS-\textit{imaginary}.
\end{definition}

\begin{remark}\label{imaginaryasCLimaginary}Work in the notation of Definition \ref{definitionequivrelationEa}.
Note that the $\mathfrak{T}_{a}/E_{a}$ can be recovered as the CL-imaginaries $\mathfrak{T}_{a}/\rho_{a}$ where $\rho_{a}(x,y)$ is the pseudometric on $\mathfrak{T}_{a}$ given by
\[\rho_{a}(x,y):=|1-\Psi(x-y)|.\]
\end{remark}

In Section 6 in \cite{Hrushovski2021AxsTW} Hrushovski describes several potentially interesting lines of research extending his work on $\mathrm{PF}^{+}$. In (6.3) he mentions that the $\mathfrak{T}_{a}/E_{a}$ could be used to show that the Kim-Pillay group is not necessarily totally disconnected for completions of $\mathrm{ACFA}^{+}$. We want to make this observation explicit and the main step will be to show the following Lemma \ref{lemmanottorsorclosedimpliesno3amalgamation} which relates 3-amalgamation (over $E$) to $\sigma$-AS-closedness (of $E$). Afterwards, in Theorem \ref{theorem3amalgamationcharacterisation} we show that this is indeed an equivalence. Hrushovski already states in (6.3) in \cite{Hrushovski2021AxsTW} that the key to determine the Kim-Pillay groups should be to show 3-amalgamation when the $\mathfrak{T}_{a}/E_{a}$ are taken into account for. He gives a sketch for the proof whenever $p_{1},p_{2}$ and $p_{3}$ are the generic type of some given $\mathfrak{T}_{a}$ and states that the general case would need further analysis. In the proof of our result (Theorem \ref{theorem3amalgamationcharacterisation}) which covers this general case we essentially reduce the problem to the situation considered by Hrushovski. We start by giving a lemma which follows from the stability of $\mathrm{ACF}$.

\begin{lemma}\label{lemmadecomposabilityofgeneraladditiveequations}
Fix some independent 3-system $(\bar{a}_{w})_{w\in\mathcal{P}(\mathbf{3})}$ over $A$. For every equation of the form $b_{12}+b_{13}+b_{23}=0$ with $b_{ij}\in\bar{a}_{ij}$, there are $c_{i}\in\bar{a}_{i}$ for $i=1,2,3$ such that \begin{itemize}
    \item $b_{12}=c_{1}+c_{2}$,
    \item $b_{13}=c_{3}-c_{1}$,
    \item $b_{23}=-c_{2}-c_{3}$.
\end{itemize}

\end{lemma}

\begin{proof} 
 Let $\Tilde{p}_{ij}:=\mathrm{tp}^{\mathrm{qf}}_{\mathcal{L}_{\mathrm{ring}}}(\bar{a}_{ij}/A)$ be the quantifier-free $\mathcal{L}_{\mathrm{ring}}$-type (over $A$) of $\bar{a}_{ij}$. Then, $\mathrm{tp}^{\mathrm{qf}}_{\mathcal{L}_{\mathrm{ring}}}(\bar{a}_{12}/\bar{a}_{3})$ is the co-heir of $\Tilde{p}_{12}$ and thus for every $1\leq i<j\leq 3$, we can write $b_{ij}=d_{ij}^{i}+d_{ij}^{j}$ where $d_{ij}^{i}\in \bar{a}_{i}$ and $d_{ij}^{j}\in \bar{a}_{j}$. (See e.g. 1.9 in \cite{acfa} or Theorem 7.2 in \cite{ludwig2025pseudofinitefieldsadditivemultiplicative} for a similar argument.) We regroup the elements and write $\delta_{1}=d_{12}^{1}+d_{13}^{1}\;\in\bar{a}_{1}$ and similarly for $\delta_{2}\in\bar{a}_{2}$ and $\delta_{3}\in\bar{a}_{3}$. It then follows that $\delta_{1}+\delta_{2}+\delta_{3}=0$ and thus by independence $\delta_{i}\in E$ for all $1\leq i\leq 3$.
Now, we simply set 
\begin{itemize}
    \item $c_{1}=d_{12}^{1}=\delta_{1}-d_{13}^{1}$
    \item $c_{2}=d_{12}^{2}=\delta_{2}-d_{23}^{2}$
    \item $c_{3}=d_{13}^{3}+\delta_{1}=\delta_{3}+\delta_{1}-d_{23}^{3}=-\delta_{2}-d_{23}^{3}$
\end{itemize}
For this constellation, it follows that $b_{12}=c_{1}+c_{2}$, $b_{13}=-c_{1}+c_{3}$ and $b_{23}=-c_{2}-c_{3}$ which completes the proof.
\end{proof}

\begin{lemma}\label{lemmanottorsorclosedimpliesno3amalgamation}
    Assume that $E=\mathrm{acl}_{\sigma}(E)$ is not $\sigma$-AS-closed. Then, 3-amalgamation does not hold over $E$.
\end{lemma}
\begin{proof}
      Take $e\in E$ such that $\mathfrak{T}_{e}$ is not realised in $E$. Let $\alpha_{1},\alpha_{2},\alpha_{3}$ be independent realisations (using 3-amalgamation in $\mathrm{ACFA}$ they exist such that $\alpha_{i}\ind_{\alpha_{j}}\alpha_{k}$ for $\{i,j,k\}=\{1,2,3\}$) of $\mathfrak{T}_{e}$ and $\bar{a}_{i}$ an enumeration of $\mathrm{acl}_{\sigma}(E\alpha_{i})$ for $1\leq i\leq 3$. Let $p_{i}:=\mathrm{tp}(\bar{a}_{i}/E)$. Further, let $b_{12}=\alpha_{1}-\alpha_{2},\;b_{13}=\alpha_{3}-\alpha_{1}$ and $b_{23}=\alpha_{2}-\alpha_{3}$. Let us argue why all of the $b_{ij}\in F$ are not in the $\mathbb{Q}$-vector space generated by $F\cap \bar{a}_{i}$ and $F\cap \bar{a}_{j}$:\\
      Assuming otherwise, we find $d_{ij}^{i}\in F\cap \bar{a}_{i}$ and $d_{ij}^{j}\in F\cap \bar{a}_{j}$ such that $b_{ij}=d_{ij}^{i}+d_{ij}^{j}$ (where $i<j$ as before). Now, with the exact same trick of regrouping the summands as in the proof of Lemma \ref{lemmadecomposabilityofgeneraladditiveequations} we can reduce to elements $c_{i}\in F\cap \bar{a}_{i}$ for $i=1,2,3$ such that $b_{12}=c_{1}-c_{2}$, $b_{13}=c_{3}-c_{1}$ and $b_{23}=c_{2}-c_{3}$. 
      We recall that we defined $b_{12}$ as $b_{12}=\alpha_{1}-\alpha_{2}$ and we obtain
      \[\alpha_{1}-\alpha_{2}=b_{12}=c_{1}-c_{2}\] which yields $\alpha_{1}-c_{1}=\alpha_{2}-c_{2}\in \bar{a}_{1}\cap \bar{a}_{2}=E$ by independence. But, on the other hand $\sigma(\alpha_{1}-c_{1})-(\alpha_{1}-c_{1})=e$ since $c_{1}\in F$ and thus we found a contradiction to the assumption that $\mathfrak{T}_{e}$ was not realised in $E$.
      Finally, since the $b_{ij}\in F$ are not in the $\mathbb{Q}$-vector space generated by $F\cap \bar{a}_{i}$ and $F\cap \bar{a}_{j}$, the following is a consequence of Corollary \ref{collaryextensionlinearindtypes}:\\
      For any $r_{12},r_{13},r_{23}\in S^{1}$ we have that $\mathrm{tp}(\bar{a}_{i}/E)\cup \mathrm{tp}(\bar{a}_{j}/E)$ is consistent with $\Psi(b_{ij})=r_{ij}$. Thus, in particular, if we choose $r_{12},r_{13},r_{23}\in S^{1}$ such that $r_{12}\cdot r_{13}\cdot r_{23}\neq 1$ and take the types $p_{ij}$ as completions of $\mathrm{tp}(\bar{a}_{i}/E)\cup \mathrm{tp}(\bar{a}_{j}/E)$ together with a sentence stating $\Psi(b_{ij})=r_{ij}$ for all the $r_{ij}$, then the $3$-amalgamation problem $(p_{w})_{w\in\mathcal{P}^{-}(\mathbf{3})}$ does not have a solution, as there is no way to extend $\Psi$ to a homomorphism on an amalgam $\bar{a}_{123}$ consistent with the types $p_{12},p_{13},p_{23}$.\end{proof}

Now, we will focus on the other direction and show that $3$-amalgamation indeed holds over $\sigma$-AS-closed sets. 

\begin{theorem}\label{theorem3amalgamationcharacterisation}
    Over $E=\mathrm{acl}_{\sigma}(E)$ 3-amalgamation holds if and only if $E$ is $\sigma$-AS-closed.
\end{theorem}

\begin{proof}
One direction is already given by Lemma \ref{lemmanottorsorclosedimpliesno3amalgamation}. For the other, assume that $E$ is $\sigma$-AS-closed. Let $(p_{w})_{w\in \mathcal{P}^{-}(\mathbf{3})}$ be a 3-amalgamation problem over $E$. Then, the system $(p_{w}\restriction_{\mathcal{L}_{\sigma}})_{w\in \mathcal{P}^{-}(\mathbf{3})}$ is a 3-amalgamation problem in the $\mathcal{L}_{\sigma}$-theory $\mathrm{ACFA}$ and by Fact \ref{factnamalgamationacfa} it has a solution $\Tilde{p}_{\mathbf{3}}$. We will now extend $\Tilde{p}_{\mathbf{3}}$ to a complete type $p_{\mathbf{3}}$ in $\mathrm{ACFA}^{+}$ which completes the system $(p_{w})_{w\in \mathcal{P}^{-}(\mathbf{3})}$. Take a realisation $\bar{a}_{\mathbf{3}}$ of $\Tilde{p}_{\mathbf{3}}$ and let $(\bar{a}_{w})_{w\in\mathcal{P}(\mathbf{3})}$ be the corresponding independent $3$-system (in $\mathcal{L}_{\sigma}$) over $E$ (where $\bar{a}_{w}\models p_{w}\restriction_{\mathcal{L}_{\sigma}}$). Next, we equip the set $A:=(\bar{a}_{12}\cup\bar{a}_{13}\cup\bar{a}_{23})\cap F$ with a map $\Psi: A\rightarrow S^{1}$ such that $\bar{a}_{w}$ together with $\Psi\restriction_{\bar{a}_{w}\cap F}$ realises $p_{w}$ for all $w\in\mathcal{P}^{-}(\mathbf{3})$. Let $V$ denote the $\mathbb{Q}$-vector space generated by $A$. Then, it suffices to show that $\Psi$ extends to a homomorphism on $V$ because once this is done we can freely choose $\Psi$ outside of $V$ by Lemma \ref{collaryextensionlinearindtypes}.\\
By compactness $\Psi$ extends to a homomorphism on $V$ if and only if for every equation $b_{12}+b_{13}+b_{23}=0$ where $b_{ij}\in\bar{a}_{ij}\cap F$ we have that $\Psi(b_{12})\cdot \Psi(b_{13})\cdot\Psi(b_{23})=1$. From Lemma \ref{lemmadecomposabilityofgeneraladditiveequations} we obtain $c_{i}\in \bar{a}_{i}$ for $i=1,2,3$ such that $b_{12}=c_{1}+c_{2}$, $b_{13}=-c_{1}+c_{3}$ and $b_{23}=-c_{2}-c_{3}$. Since $\sigma(b_{ij})=b_{ij}$, it follows that $\sigma(c_{i})-c_{i}=\sigma(c_{j})-c_{j}$ for any $1\leq i,j\leq 3$.
From the fact that $\bar{a}_{i}$ is independent from $\bar{a}_{j}$ (for $i\neq j$) over $E$ it follows that $\sigma(c_{i})-c_{i}=\sigma(c_{j})-c_{j}=e\in E$.
Finally, we use the assumption that $E$ is $\sigma$-AS-closed to find some $\Tilde{e}\in \mathfrak{T}_{e}\cap E$. We set $d_{i}:=c_{i}-\Tilde{e}$ for any $1\leq i\leq 3$. It follows that $\sigma(d_{i})-d_{i}=0$, or, in other words $d_{i}\in F$. Moreover, for any $1\leq i,j\leq 3$ with $i\neq j$ we have $c_{i}-c_{j}=d_{i}-d_{j}$. As a consequence, we obtain $b_{12}=d_{1}+d_{2}$, $b_{13}=-d_{1}+d_{3}$ and $b_{23}=-d_{2}-d_{3}$.
Now $\Psi(b_{12})=\Psi(d_{1})\cdot\Psi(d_{2})$ as $\Psi$ is a homomorphism on $\bar{a}_{12}\cap F$. This works similarly for $b_{13}$ and $b_{23}$. Putting everything together, we get \[\Psi(b_{12})\cdot \Psi(b_{13})\cdot\Psi(b_{23})=\Psi(d_{1})\cdot\Psi(-d_{2})\cdot\Psi(-d_{1})\cdot\Psi(d_{3})\cdot\Psi(-d_{2})\cdot\Psi(-d_{3})=1,\]
which is exactly what we had to show.
\end{proof}

We will now establish an analogous formulation of Theorem \ref{theorem3amalgamationcharacterisation} in the terminology of strong types. This will later be useful when we determine (the connected component of) the Kim-Pillay group.

\begin{definition}\label{defintionTORA}
    Let $E=\mathrm{acl}_{\sigma}(E)\subseteq M$. We define $\wp_{\sigma}(E)$ as the set of elements $e\in E$ such that $\mathfrak{T}_{e}$ is realised in $E$, i.e., \[\wp_{\sigma}(E):=\{e\in E\;|\;\exists b\in E\;\sigma(b)-b=e\}.\]
\end{definition}



\begin{definition}\label{definitiontorsorstrongtype} 
Let $A\subseteq K\models \mathrm{ACFA}^{+}$ be a small set. Recall the definition of a $\mathcal{D}$-strong type in Definition \ref{definitionDstrongtype}. Here we take $\mathcal{D}$ to be the set of all $A$-definable equivalence relations (with finitely many classes) together with the $E_{a}$ for $a\in \mathrm{acl}_{\sigma}(A)$ as in the above Definition \ref{definitionequivrelationEa}. Such a strong type will be called a $\sigma$-AS-strong type. We write $\bar{b}\equiv_{A}^{\wp_{\sigma}}\bar{c}$ to denote that $\bar{b}$ and $\bar{c}$ have the same  $\sigma$-AS-strong type over $A$. 
\end{definition}


\begin{lemma}\label{remarkfrom3amalgtoInd}
The proof of Theorem \ref{theorem3amalgamationcharacterisation} works in the exact same way, if we start with some $E=\mathrm{acl}_{\sigma}(E)$, not necessarily $\sigma$-AS-closed, but assume the following condition on $(p_{w})_{w\in \mathcal{P}^{-}(\mathbf{3})}$:\\ 
Whenever $\bar{f}_{12}\models p_{12}$ and $\bar{h}_{13}\models p_{13}$, then $\bar{f}_{1}\equiv_{E}^{\wp_{\sigma}}\bar{h}_{1}$ where $\bar{f}_{1}$ and $\bar{h}_{1}$ are the subtuples corresponding to realisations of $p_{1}$ and similarly for $p_{2}$ and $p_{3}$.

\end{lemma}
\begin{proof}
    We assume this and work in the exact same notation as in the proof of Theorem \ref{theorem3amalgamationcharacterisation}. Again, we only have to show that $\Psi(b_{12})\cdot \Psi(b_{13})\cdot\Psi(b_{23})=1$. As before we obtain $c_{i}\in \bar{a}_{i}$ for $i=1,2,3$ such that $b_{12}=c_{1}+c_{2}$, $b_{13}=-c_{1}+c_{3}$ and $b_{23}=-c_{2}-c_{3}$. Recall that $c_{1},c_{2},c_{3}\in \mathfrak{T}_{e}$. Now the condition in the statement of the lemma yields that $c_{1}$ is in the same $E_{e}$-class whether computed within $\bar{a}_{12}$ or $\bar{a}_{13}$ and similarly for $c_{2}$ and $c_{3}$. But then $\Psi(c_{1}+c_{2})\Psi(-c_{1}+c_{3})\Psi(-c_{2}-c_{3})=1$ already follows.
\end{proof}


\begin{corollary}\label{corollaryindependencetheoremsigmaASstrongtypes}
    The Independence Theorem holds for $\sigma$-AS-strong types.
\end{corollary}
\begin{proof}
 Recall that for every small set $A$ we have to show that whenever we have given the tuples $\bar{a}_{2}\ind_{A} \bar{a}_{3},\, \bar{c}\ind_{A}\bar{a}_{2},\,\bar{d}\ind_{A}\bar{a}_{3}$ and $\bar{c}\equiv_{A}^{\wp_{\sigma}}\bar{d}$, then we can find some $\bar{a}_{1}\ind_{A}\bar{a}_{2}\bar{a}_{3}$ such that $\bar{a}_{1}\equiv_{A\bar{a}_{2}}^{\wp_{\sigma}}\bar{c}$ and $\bar{a}_{1}\equiv_{A\bar{a}_{3}}^{\wp_{\sigma}}\bar{d}$.
We may assume that the above tuples enumerate their respective algebraic closures. Next, we consider the 3-amalgamation problem $(p_{w})_{w\in \mathcal{P}^{-}(\mathbf{3})}$ which is given by the types $\mathrm{tp}(\mathrm{acl}_{\sigma}(\bar{c}\bar{a}_{2})/A),\mathrm{tp}(\mathrm{acl}_{\sigma}(\bar{d}\bar{a}_{3})/A),\mathrm{tp}(\mathrm{acl}_{\sigma}(\bar{a}_{2}\bar{a}_{3})/A)$ together with the corresponding restrictions (using that $\mathrm{tp}(\bar{c}/A)=\mathrm{tp}(\bar{d}/A)$).
By Lemma \ref{remarkfrom3amalgtoInd} we can apply Theorem \ref{theorem3amalgamationcharacterisation} to obtain a solution $p_{\mathbf{3}}(\bar{x})$. If we now take a realisation $\bar{b}_{\mathbf{3}}\models p_{\mathbf{3}}$ and $\tau\in \mathrm{Aut}(M/A)$ such that $\tau(\bar{b}_{23})=\bar{a}_{23}$, then we can set $\bar{a}_{1}:=\tau(\bar{b}_{1})$ and it follows that $\bar{a}_{1}\ind_{A}\bar{a}_{2}\bar{a}_{3}$, $\bar{a}_{1}\equiv_{A\bar{a}_{2}}\bar{c}$ and $\bar{a}_{1}\equiv_{A\bar{a}_{3}}\bar{d}$.\\
Let us show that we can then even assume that $\bar{a}_{1}\equiv_{A\bar{a}_{2}}^{\wp_{\sigma}}\bar{c}$ and $\bar{a}_{1}\equiv_{A\bar{a}_{3}}^{\wp_{\sigma}}\bar{d}$. Note that if $A\bar{a}_{2}$ and $A\bar{a}_{3}$ are both $\sigma$-AS-closed, it already follows from $\bar{a}_{1}\equiv_{A\bar{a}_{2}}\bar{c}$ that  $\bar{a}_{1}\equiv_{A\bar{a}_{2}}^{\wp_{\sigma}}\bar{c}$ holds (and similarly for $\bar{a}_{1}\equiv_{A\bar{a}_{3}}^{\wp_{\sigma}}\bar{d}$). However, we can force this by applying Lemma \ref{lemmaextnesionindependenceinacfaplus} several times. Concretely, we replace $\bar{a}_{2}$ by $\Tilde{a}_{2}\supseteq A\bar{a}_{2}$ where $\Tilde{a}_{2}\ind_{\bar{a}_{2}}\bar{a}_{3}\bar{c}$ and $\Tilde{a}_{2}$ is $\sigma$-AS-closed. Similarly, we construct $\Tilde{a}_{3}\supseteq A\bar{a}_{3}$ $\sigma$-AS-closed such that, moreover,
$\Tilde{a}_{3}\ind_{\bar{a}_{3}}\Tilde{a}_{2}$ holds.
It remains to show that we can find some $\Tilde{d}$ with $\Tilde{d}\ind_{A}\Tilde{a}_{3}$ such that $\Tilde{d}\equiv_{A}^{\wp_{\sigma}}\bar{d}$. Let $\Tilde{b}\supseteq A\bar{a}_{3}$ be $\sigma$-AS-closed with $\bar{d}\ind_{\bar{a}_{3}}\Tilde{b}$ and $\Tilde{d}$ be such that $\Tilde{d}\ind_{\Tilde{b}}\Tilde{a}_{3}$ and $\Tilde{d}\equiv_{\Tilde{b}}d$. Then, $\Tilde{d}\equiv_{A}^{\wp_{\sigma}}\bar{d}$ holds
which shows that $\Tilde{d}$ with the required properties can be found.
Now, the above argument still works for the quadrupel $(\Tilde{a}_{2},\Tilde{a}_{3},\bar{c},\Tilde{d})$ and the $\bar{a}_{1}$ we obtain yields a solution to our initial problem.
\end{proof}

\begin{definition}
    For any $A\subseteq M$ we denote by $\mathrm{acl}_{\wp_{\sigma}}^{\mathrm{eq}}(A)$ the set $\mathrm{acl}_{\sigma}(A)\cup\{\mathbf{b}_{a}\in \mathfrak{T}_{a}/E_{a}\,|\,a\in \mathrm{acl}_{\sigma}(A)\}$.
\end{definition}

\begin{corollary}\label{corollarydescriptionofbddinACFAplus}
For any $A\subseteq M$ the set $\mathrm{bdd}(A)=\mathrm{acl}_{\mathrm{CL}}^{\mathrm{eq}}(A)$ is interdefinable with $\mathrm{acl}_{\wp_{\sigma}}^{\mathrm{eq}}(A)$.
\end{corollary}
\begin{proof}
    Using Theorem \ref{theorewhendstrongtypeequivlascarstrongtype} it follows from Corollary \ref{corollaryindependencetheoremsigmaASstrongtypes} that having the same $KP$-strong-type coincides with having the same $\sigma$-AS-strong type. Since, moreover, $\mathrm{ACFA}$ eliminates classical imaginaries (and in $\mathrm{ACFA}^{+}$ there are no new classical ones by Lemma \ref{qe-conservative}) it follows for $A=\mathrm{acl}_{\sigma}(A)$ that $\mathrm{bdd}(A)=\mathrm{acl}_{\mathrm{CL}}^{\mathrm{eq}}(A)$ is interdefinable with the union of $A$ and the sets of $\sigma$-AS-imaginaries $\mathfrak{T}_{a}/E_{a}$ for $a\in A$. 
\end{proof}

\begin{corollary}\label{corollaryacfaplusissimple}The theory $\mathrm{ACFA}^{+}$ is simple.
    
\end{corollary}
\begin{proof}
We apply Theorem 1.51 in \cite{simplicityCATs}. The first properties hold since they already hold in $\mathrm{ACFA}$ and the notion of independence coincides with the one of $\mathrm{ACFA}$. For Extension we can use Lemma \ref{collaryextensionlinearindtypes} and finally since $\equiv^{\mathrm{KP}}$ is finer then $\equiv^{\wp_{\sigma}}$ it follows directly from Theorem \ref{theorem3amalgamationcharacterisation} that the Independence Theorem holds.
\end{proof}

 \section{Description of the Kim-Pillay group}\label{sectionautomgrouptorsors}

Let $\mathcal{M}\models \mathrm{ACFA}^{+}$ be a monster and $A=\mathrm{acl}_{\sigma}(A)\subseteq M$ some small set. In this section, we will describe the group $\mathrm{Aut}(\mathrm{acl}^{\mathrm{eq}}_{\wp_{\sigma}}(A)/A)$. By Corollary \ref{corollarydescriptionconnectedcompasAutomgroup} this group is the connected component $H_{\mathrm{KP}}^{A}$ of the Kim-Pillay group of $\mathrm{Th}(\mathcal{M}_{A})$. In (6.3) in \cite{Hrushovski2021AxsTW} Hrushovski conjectures that for any completion of $\mathrm{ACFA}^{+}$ the group $H_{\mathrm{KP}}$ is abelian. This will immediately follow from our result, as we will determine $H_{\mathrm{KP}}^{A}$ up to isomorphisms (of topological groups) in Theorem \ref{theoremidentycompnenKPgrouptopequiv}. Apart from the fact that the possible non-triviality of $H_{\mathrm{KP}}$ in a natural theory like $\mathrm{AFCA}^{+}$ as well as how it naturally arises from the $\sigma$-AS equations considered in the last section are interesting phenomena in itself, Hrushovski gives another reason to investigate $H_{\mathrm{KP}}$ in 6.3 (iii) in \cite{Hrushovski2021AxsTW}. He states that the probable abelian nature of $H_{\mathrm{KP}}$ provides a reason that $\mathrm{ACFA}^{+}$ does not yet serve as a framework for a potential archimedean analogue of Grothendieck's $l$-adic local systems. While one might hope that $\mathrm{ACFA}^{+}$ will eventually serve as a preparation of grounds for such considerations, we will not provide any insights in this direction throughout this article. We will now start by considering the automorphism group of a single fixed $\mathfrak{T}_{a}/E_{a}$ before progressively turning to $\mathrm{Aut}(\mathrm{acl}^{\mathrm{eq}}_{\wp_{\sigma}}(A)/A)$.

\begin{definition}\label{definitionautomorphismtotorsor}
    Let $\mathcal{M}\models \mathrm{ACFA}^{+}$ be a monster and $A\subseteq M$ some small set. Let $a\in A$ and $e\in \mathfrak{T}_{a}$. We denote by $\mathrm{rot}_{a}:\mathrm{Aut}(\mathcal{M}/A)\rightarrow(S^{1},\cdot)$ the map that sends $\tau\in \mathrm{Aut}(\mathcal{M}/A)$ to $\tau_{a}\in S^{1}$ where $\tau_{a}:=\Psi(\tau(e)-e)$. We will call $\tau_{a}$ the \textit{rotation} (on $\mathfrak{T}_{a}/E_{a}$) induced by $\tau$. We call the image of $\mathrm{rot}_{a}$ the group of rotations of $\mathfrak{T}_{a}/E_{a}$ over $A$. 
\end{definition}
\begin{lemma}\label{lemmarotationsarhomom}
    The map $\mathrm{rot}_{a}$ from Definition \ref{definitionautomorphismtotorsor} is a well-defined homomorphism and does not depend on the choice of $e\in \mathfrak{T}_{a}$. Moreover, the automorphism group of $\sigma$-AS-imaginaries $\mathrm{Aut}_{A}(\mathfrak{T}_{a}/E_{a})$ (for $a\in A$ fixed) can be identified with the group of rotations of $\mathfrak{T}_{a}/E_{a}$ (over $A$).
\end{lemma}
\begin{proof}
    First, let us note that $\mathrm{rot}_{a}$ is indeed a well-defined map, since for $e_{1},e_{2}\in \mathfrak{T}_{a}$ we have
    \[E_{a}(e_{1},e_{2})\iff\Psi(e_{1}-e_{2})=1\iff\Psi(\tau(e_{1})-\tau(e_{2}))=1.\]
    It is a homomorphism: Let $\tau_{1},\tau_{2}\in \mathrm{Aut}(\mathcal{M}/A)$ and $e\in \mathfrak{T}_{a}$ be given. Let $r_{i}:=\Psi(\tau_{i}(e)-e)$ for $i=1,2$. But then  \[\Psi(\tau_{2}\circ\tau_{1}(e)-e)=\Psi(\tau_{2}(\tau_{1}(e))-\tau_{1}(e))\cdot\Psi(\tau_{1}(e)-e)=r_{2}\cdot r_{1}.\]
\end{proof}

\begin{lemma}\label{lemmarealisedtorsoryieldsidentity}
Let $A=\mathrm{acl}_{\sigma}(A)\subseteq M,\; a\in A$.
If $a\in \wp_{\sigma}(A)$, then $\mathrm{Aut}_{A}(\mathfrak{T}_{a}/E_{a})=\{id\}$.

\end{lemma}
\begin{proof}
    Let $\tau\in \mathrm{Aut}(\mathcal{M}/A)$ and  $\Tilde{a}\in \mathfrak{T}_{a}\cap A$. Then $\tau(\Tilde{a})=\Tilde{a}$ and consequently $\tau_{a}=1$ for any $\tau\in \mathrm{Aut}(\mathcal{M}/A)$.
\end{proof}

\begin{lemma}\label{lemmarotationonsumoftorsors}Let $a_{1},\dots, a_{n}\in A$, $(m_{1},\dots,m_{n})\in\mathbb{Z}^{n}$ and $\tau\in \mathrm{Aut}(\mathcal{M}/A)$ be given. Let $b:=\sum_{1\leq i\leq n}m_{i}a_{i}$. We can describe the rotation induced from $\tau$ on $\mathfrak{T}_{b}/E_{b}$ as follows:
\[\tau_{b}=\prod_{1\leq i\leq n}\tau_{a_{i}}^{m_{i}}.\]
\end{lemma}
\begin{proof}
    Let $e_{1},\dots, e_{n}$ be elements of $\mathfrak{T}_{a_{1}},\dots,\mathfrak{T}_{a_{n}}$ respectively. The element $\Tilde{e}:=\sum_{1\leq i\leq m}m_{i}e_{i}$ is an element from $\mathfrak{T}_{b}$ and we have
    \[\Psi(\tau(\Tilde{e})-\Tilde{e})=\Psi\left(\sum_{1\leq i\leq n}m_{i}(\tau(e_{i})-e_{i})\right)=\prod_{1\leq i\leq n}\tau_{a_{i}}^{m_{i}}.\]
    
    \end{proof}

\begin{lemma}\label{mainstatementtorsorharacterisation}
    Let $\mathcal{M}\models \mathrm{ACFA}^{+}$ be a monster and $A=\mathrm{acl}_{\sigma}(A)$ some small set. Let $\bar{a}=(a_{1},\dots,a_{n})\in A^{n}$ be $\mathbb{Q}$-linearly independent over $\wp_{\sigma}(A)$ (i.e., no non-trivial $\mathbb{Q}$-linear combination of the $a_{i}$ lies in $\wp_{\sigma}(A)$). Let $\bar{r}=(r_{1},\dots,r_{n})\in\mathbb{T}^{n}$ be given, then there is an automorphism $\tau\in \mathrm{Aut}(\mathcal{M}/A)$ such that for all $1\leq i\leq n$ we have $\tau_{a_{i}}=r_{i}$.
\end{lemma}
\begin{proof}
Let $\bar{d}\in M^{n}$, such that $d_{i}\in \mathfrak{T}_{a_{i}}$ for any $1\leq i\leq n$ and let $\bar{e}\models \mathrm{tp}(\bar{d}/A)$ be independent from $\bar{d}$ over $A$. We write $c_{i}=d_{i}-e_{i}$ and $D:=\mathrm{acl}_{\sigma}(\bar{d})\cap F$ as well as $E:=\mathrm{acl}_{\sigma}(\bar{e})\cap F$. Let $p(\bar{x},\bar{y}):=\mathrm{tp}(\bar{d}\bar{e}/A)$ and $p_{\mathcal{L}_{\sigma}}(\bar{x},\bar{y}):=\mathrm{tp}_{\mathcal{L}_{\sigma}}(\bar{d}\bar{e}/A)$. If $p_{\mathcal{L}_{\sigma}}(\bar{x},\bar{y})$ does not imply any $\mathbb{Q}$-linear dependence between $x_{1}-y_{1},\dots,x_{n}-y_{n}$ over the $\mathbb{Q}-$vector space generated by $D$ and $E$, then by Lemma \ref{collaryextensionlinearindtypes} $p(\bar{x},\bar{y})$ is consistent with $\bigwedge_{1\leq i\leq n}\Psi(x_{i}-y_{i})=r_{i}$. But from this the statement of the theorem already follows:\\ 
Let $(\bar{b},\bar{b^{\prime}})$ be such a realisation of $\mathrm{tp}(\bar{d}/A)\cup \mathrm{tp}(\bar{e}/A)\cup\{\bigwedge_{1\leq i\leq n}\Psi(x_{i}-y_{i})=r_{i}\}$. Then, we have $\mathrm{tp}(\bar{b}/A)= \mathrm{tp}(\bar{b^{\prime}}/A)$ and thus find some $\tau\in \mathrm{Aut}(\mathcal{M}/A)$ that sends $\bar{b}$ to $\bar{b^{\prime}}$.
Thus, it suffices to show that there is no such linear dependence among the $x_{1}-y_{1},\dots,x_{n}-y_{n}$ implied by $p_{\mathcal{L_{\sigma}}}(\bar{x},\bar{y})$:\\
We assume the opposite, i.e., we assume that there are $q_{1},\dots,q_{n}\in\mathbb{Q}$, not all $0$, such that for any realisation $\bar{d}\bar{e}\models p_{\mathcal{L_{\sigma}}}$ (with the notation as before) there are $\Tilde{d}\in D$ and $\Tilde{e}\in E$ such that
\[(\star)_{1}\;\;\;\sum_{1\leq i\leq n}q_{i}c_{i}=\Tilde{d}-\Tilde{e}.\]
From $(\star)_{1}$ we obtain directly by substituting $c_{i}$ by $d_{i}-e_{i}$ for every $1\leq i\leq n$ that
\[(\star)_{2}\;\;\;
\left(\sum_{1\leq i\leq n}q_{i}d_{i}\right)-\Tilde{d}=\left(\sum_{1\leq i\leq n}q_{i}e_{i}\right)-\Tilde{e}.\]
We write $d^{\prime}:=\sum_{1\leq i\leq n}q_{i}d_{i}$.
By the independence of $\bar{d}$ and $\bar{e}$ it follows from $(\star)_{2}$ that
$d^{\prime}-\Tilde{d}=g\in A$. But then since $\Tilde{d}\in F$ and $d_{i}\in \mathfrak{T}_{a_{i}}$ we have
\[(\star)_{3}\;\;\;\sigma(g)-g=\sigma(d^{\prime})-d^{\prime}=\sum_{1\leq i\leq n}q_{i}(\sigma(d_{i})-d_{i})=\sum_{1\leq i\leq n}q_{i}a_{i}.\]
We set $\Tilde{a}:=\sum_{1\leq i\leq n}q_{i}a_{i}$. Then $(\star)_{3}$ yields that $\Tilde{a}\in \wp_{\sigma}(A)$ contradicting the assumption that $a_{1},\dots, a_{n}$ are $\mathbb{Q}$-linearly independent over $\wp_{\sigma}(A)$.
\end{proof}

\begin{remark}Let $A=\mathrm{acl}_{\sigma}(A)$, $a\in A\backslash \wp_{\sigma}(A)$ and $\tau\in \mathrm{Aut}(\mathcal{M}/A)$ be given. Now $\tau_{a}=\tau_{\sigma(a)}$ (and thus $\tau_{a}=\tau_{\sigma^{n}(a)}$ for all $n\in\mathbb{N}$) since for any $e\in \mathfrak{T}_{a}$, we have $\sigma(e)\in \mathfrak{T}_{\sigma(a)}$ and moreover
\[\Psi(\tau(e)-e)=\Psi(\sigma(\tau(e)-e))=\Psi(\tau(\sigma(e))-\sigma(e)).\]
Now it might not be entirely obvious a priori why this case does not contradict Lemma \ref{mainstatementtorsorharacterisation}. However, if $b=\sigma^{n}(a)$, then
\[a+\sum_{1\leq i\leq n}\sigma^{i}(a)-\sigma^{i-1}(a)=b\] and clearly $\sigma^{i}(a)-\sigma^{i-1}(a)\in  \wp_{\sigma}(A)$ for all $1\leq i\leq n$. Thus, $b$ and $a$ are not $\mathbb{Q}$-linearly independent over $ \wp_{\sigma}(A)$.

\end{remark}

\begin{definition}\label{definitiongeneraltorsorautomgroup}
    Work with the notation of Definition \ref{definitionautomorphismtotorsor}. Let $B\subseteq A$ some subset. We denote by $\mathrm{rot}_{B}:\mathrm{Aut}(\mathcal{M}/A)\rightarrow\prod_{B}(S^{1},\cdot)$ the homomorphism that is given by the map $\mathrm{rot}_{a}$ for $a\in B$ in each component. We denote by $G_{B}$ the group that is given by the image of $\mathrm{rot}_{B}$. Note that it depends on $A$, which nonetheless we suppress in the notation as it will not change for the rest of this section.
\end{definition}

\begin{notation}
    We denote by $\varprojlim \mathbb{R}/n\mathbb{Z}$ the inverse limit over the index set $\mathbb{N}$ with transition maps given by the projections
    $\mathbb{R}/m\mathbb{Z}\twoheadrightarrow\mathbb{R}/l\mathbb{Z}$ for $l|m$.
\end{notation}

\begin{remark}\label{remarkcommutingdiagram}
    For every $m\in\mathbb{N}$ we have a natural identification $\xi_{m}$ between $\mathbb{R}/m\mathbb{Z}$ and $S^{1}$ defined by 
    \[\xi_{m}(x):=\exp\left(2\pi i\frac{x}{m}.\right)\]
    With this identification, the following diagram commutes for $l|m$.
    \begin{center}
        \begin{tikzcd}
\mathbb{R}/m\mathbb{Z} \arrow[rr, two heads] \arrow[dd, "\xi_{m}"'] &  & \mathbb{R}/l\mathbb{Z} \arrow[dd, "\xi_{l}"] \\
                                         &  &                   \\
 S^{1} \arrow[rr, "x\,\rightarrow\, x^{m/l}"]                        &  & S^{1}                
\end{tikzcd}
    \end{center}
\end{remark}

\begin{proposition}\label{automgroupofonetorsor}

    Let $\mathcal{M}\models \mathrm{ACFA}^{+}$ be a monster model, $A=\mathrm{acl}_{\sigma}(A)\subseteq M$ small and $a\in A\backslash \wp_{\sigma}(A)$ be given. We set $\mathbb{Q}a:=\{qa\,|\,q\in\mathbb{Q}\}$. Then $G_{\mathbb{Q}a}\cong \varprojlim\mathbb{R}/n\mathbb{Z}$.
\end{proposition}

\begin{proof}
We write $\frac{1}{\mathbb{N}}a:=\{\frac{1}{n}a\,|\,n\in\mathbb{N}\}$. From Lemma \ref{mainstatementtorsorharacterisation} (since $a\in A\backslash \wp_{\sigma}(A)$) and saturation it follows that $G_{\frac{1}{n}a}\cong S^{1}$ for every $n\in\mathbb{N}$.
We will first show that $G_{\frac{1}{\mathbb{N}}a}\cong \varprojlim\mathbb{R}/n\mathbb{Z}$. Fix $n\in\mathbb{N}_{>1}$ and $\tau\in \mathrm{Aut}(\mathcal{M}/A)$. We observe that $\tau_{\frac{1}{n}a}$ is indeed not uniquely determined by $\tau_{a}$. Precisely, there are $n$-many choices for $\tau_{\frac{1}{n}a}$ consistent with $\tau_{a}$ since for $e\in \mathfrak{T}_{a}$ we have
\[\left(\tau_{\frac{1}{n}a}\right)^{n}=\Psi\left(\frac{1}{n}(\tau(e)-e)\right)^{n}=\Psi(\tau(e)-e)=\tau_{a},\]
so if $\tau_{a}=\exp(2\pi ir_{a})$, then $\tau_{\frac{1}{n}a}$ is of the form $\zeta_{n}^{k}\cdot \exp(2\pi i\frac{r_{a}}{n})$ where $\zeta_{n}=\exp(2\pi i\frac{k}{n})$ for $1\leq k\leq n$. By Lemma \ref{mainstatementtorsorharacterisation} every such choice can arise, i.e., for every $1\leq k\leq n$ there is $\Tilde{\tau}\in \mathrm{Aut}(\mathcal{M}/A)$ with $\Tilde{\tau}_{a}=\tau_{a}$ and $\Tilde{\tau}_{\frac{1}{n}a}=\zeta_{n}^{k}\cdot \exp(2\pi i\frac{r_{a}}{n})$. Let $m\in\mathbb{N}$ be such that $n|m$. Then, for $\tau, e$ as above, we have $\left(\tau_{\frac{1}{m}a}\right)^{\frac{m}{n}}=\tau_{\frac{1}{n}a}$
and using Remark \ref{remarkcommutingdiagram} we obtain that the following diagram commutes for $n|m$.
\begin{center}
    \begin{tikzcd}
\mathbb{R}/m\mathbb{Z} \arrow[rrrr, two heads] \arrow[dd, "\xi_{m}"'] &  &  & &\mathbb{R}/n\mathbb{Z} \arrow[dd, "\xi_{n}"] \\
                                         &  &  &                 \\
 G_{\frac{1}{m}a} \arrow[rrrr, "\tau_{(1/m)a}\,\rightarrow\, \left(\tau_{(1/m)a}\right)^{m/n}"]                      &  &  &  &G_{\frac{1}{n}a}                
\end{tikzcd}
\end{center}
Thus we obtain that $G_{\frac{1}{\mathbb{N}}a}\cong\varprojlim\mathbb{R}/n\mathbb{Z}$ (again using that by Lemma \ref{mainstatementtorsorharacterisation} every possible sequence $\left(\tau_{\frac{1}{n}a}\right)_{n\in\mathbb{N}}$ arises). To complete the proof, it remains to show that the projection $G_{\mathbb{Q}a}\twoheadrightarrow G_{\frac{1}{\mathbb{N}}a}$ is injective. But this is immediate by Lemma \ref{lemmarotationonsumoftorsors}.
\end{proof}

\begin{theorem}\label{corollarydescriptionofautomgroupofalltorsors}
Let $W$ be a $\mathbb{Q}$-basis of a vector space complement of $\wp_{\sigma}(A)$ in $A=\mathrm{acl}_{\sigma}(A)$. Then we have
\[G_{A}\cong\prod_{W}\left(\varprojlim\mathbb{R}/n\mathbb{Z}\right).\]
Moreover, $G_{A}$ is the connected component of the Kim-Pillay group $\mathrm{Gal}_{\mathrm{KP}}^{A}$.
\end{theorem}
\begin{proof}
    This follows directly from Lemma \ref{mainstatementtorsorharacterisation} and Proposition \ref{automgroupofonetorsor}:
    For any $w\in W$, we have $G_{\frac{1}{\mathbb{N}}w}\cong \varprojlim\mathbb{R}/n\mathbb{Z}$ by Proposition \ref{automgroupofonetorsor}. Furthermore, by compactness it follows from Lemma \ref{mainstatementtorsorharacterisation} and Lemma \ref{lemmarotationonsumoftorsors} that $G_{\langle W\rangle_{\mathbb{Q}}}\cong \prod_{W}\left(\varprojlim\mathbb{R}/n\mathbb{Z}\right)$ holds (with $\langle W\rangle_{\mathbb{Q}}$ the $\mathbb{Q}$-vector space generated by $W$). Finally, by Lemma \ref{lemmarealisedtorsoryieldsidentity} we have $G_{\wp_{\sigma}(A)}=1$ and consequently $G_{A}\cong\prod_{W}\left(\varprojlim\mathbb{R}/n\mathbb{Z}\right)$. The moreover part follows directly from Corollary \ref{corollarydescriptionofbddinACFAplus} and Corollary \ref{corollarydescriptionconnectedcompasAutomgroup}.
\end{proof}

After having obtained this purely group-theoretic description of $\mathrm{Gal}_{\mathrm{KP}}^{A}$, we will now turn to its topology. 
\begin{definition}\label{definitionzetabmap}
Let $a_{1},\dots,a_{n}\in A\backslash\wp_{\sigma}(A)$ and $b_{i}\in \mathfrak{T}_{a_{i}}$ for all $1\leq i\leq n$ be such that, moreover, whenever $a_{j}=\sum_{i\neq j} q_{i}a_{i}$ for $q_{i}\in\mathbb{Q}$, then $b_{j}=\sum_{i\neq j} q_{i}b_{i}$ also has to hold. We write $\zeta_{\bar{b}}:\mathfrak{T}_{a_{1}}/E_{a_{1}}\times\cdots\times \mathfrak{T}_{a_{n}}/E_{a_{n}}\rightarrow \mathbb{T}^{n}$ for the map whose $i$-th component is induced by the map sending $c\in \mathfrak{T}_{a_{i}}$ to $\Psi(c-b_{i})\in S^{1}$.    
\end{definition}

\begin{lemma}\label{lemmahomeomorphicS1andTorsorim}
For all $1\leq i\leq n$, let $\mathfrak{T}_{a_{i}}/E_{a_{i}}$ be equipped with the logic topology and $\mathfrak{T}_{a_{1}}/E_{a_{1}}\times\cdots\times \mathfrak{T}_{a_{n}}/E_{a_{n}}$ by the product topology. Then, the map $\zeta_{\bar{b}}$ from Definition \ref{definitionzetabmap} is a $G_{A}$-equivariant homeomorphism.
\end{lemma}
\begin{proof}
    From Lemma \ref{lemmarotationsarhomom} and Lemma \ref{lemmarotationonsumoftorsors} it follows that $\zeta_{\bar{b}}$ is a $G_{A}$-equivariant bijection. Both spaces are compact Hausdorff. Thus, it suffices to show that $\zeta_{\bar{b}}$ is continuous. For any closed intervals $I_{i}$ in $S^{1}$ the set defined by $\bigwedge_{1\leq i\leq n}\Psi(x_{i}-b_{i})\in I_{i}$ is type-definable over the parameters $\bar{b}$ and thus closed in the logic topology. (Recall that when defining the logic topology, a small set of parameters from $M\backslash A$ is permitted.) Products of closed intervals form a subbasis of closed sets in $\mathbb{T}^{n}$ and hence $\zeta_{\bar{b}}$ is continuous.
\end{proof}

\begin{theorem}\label{theoremidentycompnenKPgrouptopequiv}
    As a topological group, the connected component $H_{\mathrm{KP}}^{A}$ of the Kim-Pillay group $\mathrm{Gal}_{\mathrm{KP}}^{A}$ is homeomorphic to $G_{A}\cong\prod_{W}\left(\varprojlim\mathbb{R}/n\mathbb{Z}\right)$, where the latter is equipped with the topology induced by the euclidean one.
\end{theorem}
\begin{proof}
By Corollary \ref{corollarydescriptionofbddinACFAplus} $\mathrm{acl}_{\mathrm{CL}}^{\mathrm{eq}}(A)$ is interdefinable with $\mathrm{acl}_{\wp_{\sigma}}^{\mathrm{eq}}(A)$. Let $J$ be the directed system consisting of finite products of the latter, that is, products of the form $\mathfrak{T}_{a_{1}}/E_{a_{1}}\times\cdots\times \mathfrak{T}_{a_{n}}/E_{a_{n}}$. (Note that we can extend $E_{i}$ to the whole sort by one class outside of $\mathfrak{T}_{a_{i}}$, if necessary.)
By Lemma \ref{lemmaKPgruppeinverserlimesmitPKstarkemsystem} $\mathrm{Gal}_{\mathrm{KP}}^{A}\cong\varprojlim_{J}\mathrm{Aut}(\mathfrak{T}_{a_{1}}/E_{a_{1}}\times\cdots\times \mathfrak{T}_{a_{n}}/E_{a_{n}})$ as topological groups where the latter is equipped with the topology induced from the inverse limit and $\mathrm{Aut}(\mathfrak{T}_{a_{1}}/E_{a_{1}}\times\cdots\times \mathfrak{T}_{a_{n}}/E_{a_{n}})$ with the compact-open topology for its action on $\mathfrak{T}_{a_{1}}/E_{a_{1}}\times\cdots\times \mathfrak{T}_{a_{n}}/E_{a_{n}}$.
By Lemma \ref{lemmahomeomorphicS1andTorsorim} there is an $\mathrm{Aut}(\mathfrak{T}_{a_{1}}/E_{a_{1}}\times\cdots\times \mathfrak{T}_{a_{n}}/E_{a_{n}})$-equivariant homeomorphism $\mathfrak{T}_{a_{1}}/E_{a_{1}}\times\cdots\times \mathfrak{T}_{a_{n}}/E_{a_{n}}\cong \mathbb{T}^{n}$.
As $\mathrm{Aut}(\mathfrak{T}_{a_{1}}/E_{a_{1}}\times\cdots\times \mathfrak{T}_{a_{n}}/E_{a_{n}})$ acts via component-wise rotation on $\mathbb{T}^{n}$, the compact-open topology coincides with the euclidean topology, and thus it follows that the isomorphism between $\mathrm{Gal}_{\mathrm{KP}}^{A}$ and  $\prod_{W}\left(\varprojlim\mathbb{R}/n\mathbb{Z}\right)$ from Theorem \ref{corollarydescriptionofautomgroupofalltorsors} is indeed an isomorphism of topological groups.
\end{proof}

\section{Elimination of CL-imaginaries up to $\sigma$-AS-imaginaries}\label{sectioneliminationofimaginaries}
As we have seen, the only obstruction to $3$-amalgamation comes from the $\sigma$-AS-imaginaries. It is natural to ask whether this gives us a description of a class of CL-imaginaries which it suffices to add to obtain elimination of all CL-imaginaries. We will show that weak elimination of CL-imaginaries indeed follows almost directly. Whereas we will not get full elimination of CL-imaginaries, we will later see that we can slightly improve and obtain elimination of CL-imaginaries up to finite sets of finite tuples (instead of up to compact sets).

\begin{definition}\label{definitionelominationofimaginaries} (Definition 4.2 in \cite{thorn_fork_cont}.) Fix a continuous logic theory $T$ and a saturated model $\mathcal{M}\models T$. We say that $T$
has weak elimination of imaginaries, if for every $a\in \mathcal{M}^{\mathrm{eq}}$ there is a compact set $B$ from $M$ such that $a\in \mathrm{dcl}_{\mathrm{CL}}(B)$ and $B\subseteq \mathrm{acl}_{\mathrm{CL}}(a)$. We say that $T$
has elimination of imaginaries up to finite sets of finite tuples, if in the above we can moreover always assume that $B$ is a finite set of finite tuples.
    
\end{definition}


For convenience, we will add only one sort that captures the whole of $\mathrm{acl}_{\wp_{\sigma}}^{\mathrm{eq}}(\emptyset)$.

\begin{definition}\label{definitionlanguageextensionfortorsors}
    Let $\mathcal{L}_{\sigma,+}^{\mathrm{EI}}$ be the expansion of the language $\mathcal{L}_{\sigma}^{+}$ such that $\mathcal{L}_{\sigma,+}^{\mathrm{EI}}$-structures consist of two sorts $(K,H)$ where $K$ is an $\mathcal{L}_{\sigma}^{+}$-structure, $H$ an $\mathcal{L}_{\emptyset}$-structure (but equipped with a non-trivial metric $\rho_{H}$) and a map $\pi:K\rightarrow H$. 
\end{definition}
We will now define the $\mathcal{L}_{\sigma,+}^{\mathrm{EI}}$-theory $\mathrm{ACFA}^{+}_{\mathrm{EI}}$ such that it extends $\mathrm{ACFA}^{+}$ by axioms ensuring that $\rho_{H}$ is induced by the pseudometric on $\mathfrak{T}_{a}$ for every $a\in K$.
\begin{definition}\label{definitionEIextension}
    We define $(K,H)$ to be a model of $\mathrm{ACFA}^{+}_{\mathrm{EI}}$, if
    \begin{itemize}
        \item $K\models \mathrm{ACFA}^{+}$.
        \item $\pi(K)$ is dense in $H$.
        \item $\rho(\pi(x),\pi(y)):=\begin{cases}
        |1-\Psi(x-y)|,\;\;\text{if}\;\;\sigma(x)-x=a=\sigma(y)-y\\
        4,\;\;\;\;\;\;\;\;\;\;\;\;\;\;\;\text{otherwise}.
    \end{cases}$
       
    \end{itemize}
\end{definition}

\begin{notation}
    From now on we work in a fixed (sufficiently saturated) model $\mathcal{M}\models \mathrm{ACFA}^{+}_{\mathrm{EI}}$. We write $\mathcal{M}=(K,H)$ where by abuse of notation $K,H$ will stand for $K(\mathcal{M}), H(\mathcal{M})$ respectively.
\end{notation}

\begin{lemma}\label{lemmabddaclCLequalsaclinexpansion}
    Let $E=\mathrm{acl}_{\sigma}(E)\subseteq K$. Then, $\mathrm{bdd}(E)=\mathrm{acl}_{\mathrm{CL}}^{\mathrm{eq}}(E)$ is interdefinable with $\pi(\wp_{\sigma}^{-1}(E))$ where $\wp_{\sigma}^{-1}(E)=\{x\in K\,|\,\sigma(x)-x\in E\}$. 
\end{lemma}
\begin{proof}
    This is immediate by Remark \ref{imaginaryasCLimaginary}.
\end{proof}
\begin{lemma}\label{lemmaaclwthimaginaries}
    Let $(A,B)\subseteq (K,H)$. The CL-algebraic closure of $(A,B)$ in $\mathcal{M}$ is given by $\mathrm{acl}_{\mathrm{CL}}(A,B)=(\Tilde{A},\pi(\wp_{\sigma}^{-1}(\Tilde{A})))$ for $\Tilde{A}=\mathrm{acl}_{\sigma}(A,\wp_{\sigma}(\pi^{-1}(B)))$.
\end{lemma}
\begin{proof}
For any $c\in K$ the corresponding space of $\sigma$-AS imaginaries $\mathfrak{T}_{c}/E_{c}=\pi(\wp_{\sigma}^{-1}(c))$ is compact, definable from $c$ and also $c$ is definable from $\mathfrak{T}_{c}/E_{c}$. The statement then follows using the fact that any automorphism $\tau\in \mathrm{Aut}(\mathcal{M})$ fixes $c$ if and only if it fixes $\mathfrak{T}_{c}/E_{c}$ set-wise. 
\end{proof}

\begin{definition}
    Given $C\subseteq A,B$ all subsets of some $\mathcal{M}\models \mathrm{ACFA}^{+}_{\mathrm{EI}}$, we say that $A\ind_{C}B$ if the $K$-part of $\mathrm{acl}_{\mathrm{CL}}(A)$ is (field-theoretically) algebraically independent from the $K$-part of $\mathrm{acl}_{\mathrm{CL}}(B)$ over the $K$-part of $\mathrm{acl}_{\mathrm{CL}}(C)$.
\end{definition}


\begin{corollary}\label{lemma3amlagamacfactor}
    In $\mathrm{ACFA}^{+}_{\mathrm{EI}}$ the Independence Theorem holds over all algebraically closed sets. 
\end{corollary}
\begin{proof}
 We have to show, given $A=\mathrm{acl}(A)\subseteq M$ and $\bar{a}_{2}\ind_{A} \bar{a}_{3},\, \bar{c}\ind_{A}\bar{a}_{2}$, $\bar{d}\ind_{A}\bar{a}_{3}$ and $\bar{c}\equiv_{A}\bar{d}$, that there are $\bar{a}_{1}\ind_{A}\bar{a}_{2}\bar{a}_{3}$ such that $\bar{a}_{1}\equiv_{A\bar{a}_{2}}\bar{c}$ and $\bar{a}_{1}\equiv_{A\bar{a}_{3}}\bar{d}$. Without loss of generality we can assume to work with tuples that enumerate their corresponding algebraic closure. As in the proof of Corollary \ref{corollaryindependencetheoremsigmaASstrongtypes} we can even assume that the $K$-part of $\bar{c},\bar{d},\bar{a}_{2}$ and $\bar{a}_{3}$ is $\sigma$-AS-closed. But in that case the $H$-part of the tuple is already in the definable closure of its $K$-part and the problem reduces to working with $\bar{c},\bar{d},\bar{a}_{2}$ and $\bar{a}_{3}$ in $K$.
For any tuples $\bar{a},\bar{b}$ from $K$ and $A\subseteq M$ we have by Lemma \ref{lemmabddaclCLequalsaclinexpansion} that $\mathrm{tp}(\bar{a}/\mathrm{acl}_{\mathrm{CL}}(A))=\mathrm{tp}(\bar{b}/\mathrm{acl}_{\mathrm{CL}}(A))$ in $\mathcal{M}$ holds if and only if for $\Tilde{A}$, the $K$-part of $\mathrm{acl}_{\mathrm{CL}}(A)$, we have
   $\mathrm{tp}(\bar{a}/\mathrm{acl}_{\mathrm{CL}}^{\mathrm{eq}}(\Tilde{A}))=\mathrm{tp}(\bar{b}/\mathrm{acl}_{\mathrm{CL}}^{\mathrm{eq}}(\Tilde{A}))$ holds (where we are working in the $eq$-expansions of $K\models \mathrm{ACFA}^{+}$). Now, the latter is equivalent to $\bar{a}\equiv_{\Tilde{A}}^{\wp_{\sigma}}\bar{b}$ and we thus can apply Corollary \ref{corollaryindependencetheoremsigmaASstrongtypes}.
\end{proof}

\begin{theorem}\label{theoremweakEIACFAplusEI}
  The theory $\mathrm{ACFA}^{+}_{\mathrm{EI}}$ has weak elimination of imaginaries.
\end{theorem}
\begin{proof}
    We follow the proof of elimination of imaginaries in $\mathrm{ACFA}$ (1.10 in \cite{acfa}). We work in $\mathcal{M}^{\mathrm{eq}}$ (in the sense of continuous logic) and $\mathrm{dcl}_{\mathrm{CL}}, \mathrm{acl}_{\mathrm{CL}}$ are considered to be taken in $\mathcal{M}^{\mathrm{eq}}$ for the rest of this proof. Let $e$ be an imaginary element. We have to show that there is some compact set $B$ in $M$ such that any automorphism fixes $e$ if and only if it fixes $B$. Let $E=\mathrm{acl}_{\mathrm{CL}}(e)\cap M$. Take a tuple $\bar{a}\in M^{n}$ and a definable function $f$ such that $f(\bar{a})=e$. Let $P$ denote the set of realisations of $\mathrm{tp}(\bar{a}/E)$. If we can show that $f$ is constant on $P$, then, it follows that $e\in \mathrm{dcl}_{\mathrm{CL}}(E)$, i.e., weak elimination of imaginaries holds.\\
    To show that $f$ is constant on $P$ first we show that there is $\bar{c}\in P$ independent from $\bar{a}$ over $E$ with $f(\bar{c})=e$. In \cite{Conant2022SeparationFI} it is shown (Corollary 4.9) using an adaption of Neumann's lemma to the metric setting that we can find a conjugate $\bar{d}$ of $\bar{a}$ over $E\cup \{e\}$ such that $\mathrm{acl}_{\mathrm{CL}}(E\bar{a})\cap \mathrm{acl}_{\mathrm{CL}}(E\bar{a})\cap M=E$. Now we can proceed with the exact same proof as in 1.10 in \cite{acfa} to obtain the $\bar{c}$ as desired since our notion of independence stems from the one in $\mathrm{ACFA}$ and thus the argument via ranks works exactly in the same way.
    Next, we will use that the Independence Theorem holds over algebraically close sets in $\mathrm{ACFA}^{+}_{\mathrm{EI}}$ (Corollary \ref{lemma3amlagamacfactor}). Assume $f$ is not constant, then we find $\bar{d}\in P$ with $f(\bar{d})\neq f(\bar{a})$. We consider the Independence theorem for $\bar{a}_{0},\bar{a}_{1},\bar{b}_{0},\bar{b}_{1}$ such that $\mathrm{tp}(\bar{a}_{0}\bar{a}_{1})=\mathrm{tp}(\bar{a}\bar{c}/E)$, $\mathrm{tp}(\bar{b}_{0}/E\bar{a}_{0})=\mathrm{tp}(\bar{c}/E\bar{a})$ and $\mathrm{tp}(\bar{b}_{1}/E\bar{a}_{1})=\mathrm{tp}(\bar{d}/E\bar{a})$. Then, this configuration would be a counterexample to the Independence Theorem over $E$ contradicting Corollary \ref{lemma3amlagamacfactor}.
\end{proof}

\begin{definition}\label{definitioncodecompactsets}
    Let $T$ be a theory with monster model $\mathcal{N}\models T$. We say that $T$ \textit{codes compact sets up to finite sets of finite tuples} if for every compact set $A$ of finite tuples of the same length there is a finite set of finite tuples $\{a_{1},\dots,a_{n}\}$ such that for every automorphism $f\in \mathrm{Aut}(\mathcal{N})$ we have $f(\{a_{1},\dots,a_{n}\})=\{a_{1},\dots,a_{n}\}$ if and only $f(A)=A$. If we can always assume that $n=1$, i.e., $A$ is coded by a finite tuple, then we say that $T$ \textit{strongly codes compact sets}.
\end{definition}


\begin{notation}\label{notationproductsoftorsorgroups}
    Let $\bar{a}=a_{1}\dots a_{n}$ be a finite tuple of elements from $K$ and $G\leq \mathrm{Aut}(\mathcal{M})$ a group stabilising $\{a_{1},\dots, a_{n}\}$ setwise. We denote by $G_{\bar{a}}$ the subgroup of finite index of $G$ that fixes $\bar{a}$ pointwise. Then we let $\Tilde{G}_{\bar{a}}$ be the subgroup of $\mathrm{Aut}(\mathfrak{T}_{a_{1}}/E_{a_{1}})\times\dots\times \mathrm{Aut}(\mathfrak{T}_{a_{n}}/E_{a_{n}})\cong\mathbb{T}^{n}$ induced by $G_{\bar{a}}$.
\end{notation}

\begin{lemma}\label{lemmacodeofcompactstabilizedset}
   We use the notation as in \ref{notationproductsoftorsorgroups}.
    Let $D\subseteq \mathfrak{T}_{a_{1}}/E_{a_{1}}\times\dots\times \mathfrak{T}_{a_{n}}/E_{a_{n}}$ be compact and denote by $H_{\{D\}}$ the subgroup of $\Tilde{G}_{\bar{a}}$ that fixes $D$ setwise. Then we find finitely many $b_{j}\in \mathfrak{T}_{\Tilde{a}_{j}}/E_{\Tilde{a}_{j}}$ where $\Tilde{a}_{j}$ is a linear combination with integer coefficients of $a_{1},\dots,a_{n}$ such that $H_{\{D\}}=H_{\bar{b}}$ (for $H_{\bar{b}}$ the subgroup of $\Tilde{G}_{\bar{a}}$ that fixes $\bar{b}=b_{1}\dots b_{k}$ point-wise.
    
\end{lemma}
\begin{proof}
    Since $H_{\{D\}}\leq\Tilde{G}_{\bar{a}}\leq \mathrm{Aut}(\mathfrak{T}_{a_{1}}/E_{a_{1}})\times\dots\times \mathrm{Aut}(\mathfrak{T}_{a_{n}}/E_{a_{n}})\cong\mathbb{T}^{n}$ we consider $H_{\{D\}}$ as a group of rotations on the components $\mathfrak{T}_{a_{i}}/E_{a_{i}}$. 
    As $D$ is compact (and $\mathcal{M}$ sufficiently saturated) $H_{\{D\}}$ is closed. We can apply Fact \ref{torsorsubgroup} to obtain for $1\leq j\leq k$ some tuples $(m_{1,j},\dots,m_{n,j})\in \mathbb{Z}^{n}\backslash 0$ such that $H_{\{D\}}$ is given by those $(g_{1},\dots,g_{n})$ such that for all $1\leq j\leq k$
    \[(\star)\;\;\;\prod_{1\leq i \leq n}g_{i}^{m_{i,j}}=1.\] We suppose $k$ is chosen being minimal as such. Now choose some $(d_{1},\dots,d_{n})\in D$ (it actually suffices to choose $(d_{1},\dots,d_{n})\in \mathfrak{T}_{a_{1}}/E_{a_{1}}\times\dots\times \mathfrak{T}_{a_{n}}/E_{a_{n}}$ arbitrarily) and $(e_{1},\dots, e_{n})\in K^{n}$ such that for all $1\leq i\leq n$ we have that the class of $e_{i}$ equals the $\sigma$-AS-imaginary $d_{i}$. Let $\Tilde{a}_{j}=\sum_{1\leq i\leq n}m_{i,j}a_{i}$. For $1\leq j\leq k$, we take $b_{j}$ to be the class of $\sum_{1\leq i\leq n}m_{i,j}e_{i}$ in $\mathfrak{T}_{\Tilde{a}_{j}}/E_{\Tilde{a}_{j}}$. Let $\bar{b}=(b_{1},\dots ,b_{k})$. To complete the proof, we will show that $H_{\{D\}}= H_{\bar{b}}$ holds.\\ Every element $\tau\in H_{\{D\}}$ satisfies $(\star)$ for all $1\leq j\leq k$. But as $\tau=(\tau_{1},\dots,\tau_{n})$ is simply a tuple of rotations of  $\mathfrak{T}_{a_{1}}/E_{a_{1}},\dots,\mathfrak{T}_{a_{n}}/E_{a_{n}}$ respectively, it yields a rotation of $\mathfrak{T}_{\Tilde{a}_{j}}/E_{\Tilde{a}_{j}}$ for any $1\leq j\leq k$ as described in Lemma \ref{lemmarotationonsumoftorsors}. Then, by $(\star)$ this rotation is just the identity, i.e., every element in $\mathfrak{T}_{\Tilde{a}_{j}}/E_{\Tilde{a}_{j}}$ is fixed for any $1\leq j\leq k$. The other direction just follows as every $\tau\in \Tilde{G}_{\bar{a}}$ that satisfies the equations in $(\star)$ already is an element from $H_{\{D\}}$ by Fact \ref{torsorsubgroup}.

\end{proof}

\begin{remark}\label{remarkfreechoiceinlinearcombinations}
    Note that the construction in the proof does not depend on the choice of the $b_{j}$ from $\mathfrak{T}_{a_{j}}/E_{a_{j}}$ for $1\leq j\leq n$.
\end{remark}


\begin{lemma}\label{lemmacodingofcompactsetsacfator}
    The theory $\mathrm{ACFA}_{\mathrm{EI}}^{+}$ codes compact sets up to finite sets of finite tuples. 
\end{lemma}
\begin{proof}
    Let $A\subseteq M^{n}$ be compact and $G_{A}\subseteq \mathrm{Aut}(\mathcal{M})$ the subgroup consisting of those $\tau\in \mathrm{Aut}(\mathcal{M})$ such that $\tau(A)=A$. We can assume that $A$ is a set of tuples $(\bar{\alpha},\bar{\beta})$ where $\bar{\alpha}$ is a tuple from $K^{l}$ and $\bar{\beta}$ is a tuple from $H^{r}$. For any $1\leq i\leq r$, we define $\Tilde{a}_{i}:=\wp_{\sigma}(\pi^{-1}(\beta_{i}))$. Note that $\Tilde{\alpha}\in \mathrm{dcl}(\bar{\beta})$. For any $\bar{a}=(\bar{\alpha},\bar{\beta})\in A$ we write $\delta_{\bar{a}}$ for the corresponding tuple $(\bar{\alpha},\Tilde{a},\bar{\beta})$ (which lives in the product sort $K^{l+r}\times H^{r}$).
    We set $\Tilde{A}:=\{\delta_{\bar{a}}\,|\,\bar{a}\in A\}$. It follows for any $\tau\in \mathrm{Aut}(\mathcal{M})$ that $\tau(A)=A$ holds if and only if $\tau(\Tilde{A})=\Tilde{A}$ holds (using that above $\Tilde{\alpha}\in \mathrm{dcl}(\bar{\beta})$). Let $\Tilde{A}_{K}$ be the projection of $\Tilde{A}$ to $K^{l+r}$. Since $A$ is compact, $\Tilde{A}_{K}$ is finite. 
    Let $\Tilde{G}$ be the subgroup of $G_{A}$ that fixes $\Tilde{A}_{K}$ point-wise. For any $\bar{d}=(\bar{\alpha},\Tilde{\alpha})\in\Tilde{A}_{K}$ let $B_{\bar{d}}:=\{\bar{\beta}\,|\,(\bar{\alpha},\Tilde{\alpha},\bar{\beta})\in\Tilde{A}\}$.\\
    We can then apply Lemma \ref{lemmacodeofcompactstabilizedset} to the group $\Tilde{G}$ and the subgroup $\Tilde{H}\leq \prod_{1\leq i\leq r}\mathrm{Aut}(\mathfrak{T}_{\Tilde{\alpha}_{i}}/E_{\Tilde{\alpha}_{i}})$ induced by $\Tilde{G}$. We    obtain a tuple $\bar{b}_{\bar{d}}$ such that $\Tilde{H}_{\bar{b}_{\bar{d}}}=\Tilde{H}_{B_{\bar{d}}}$
    where the latter denote the subgroups of $\Tilde{H}$ fixing $\bar{b}_{\bar{d}}$ and $B_{\bar{d}}$ respectively. Let $\bar{b}:=(\bar{b}_{\bar{d}_{1}},\dots,\bar{b}_{\bar{d}_{s}} )$ for $\bar{d}_{1},\dots,\bar{d}_{s}$ an enumeration of $\Tilde{A}_{K}$. Then $\Tilde{H}_{\bar{b}}=\bigcap_{\bar{d}\in\Tilde{A}_{K}}\Tilde{H}_{B_{\bar{d}}}$.
    Now $G/\Tilde{G}$ is finite. We consider its action on $\Tilde{A}_{K}$. We choose a representative of every $G/\Tilde{G}$-orbit in $\Tilde{A}_{K}$ and take $\bar{c}$ to be a tuple consisting of all those representatives together with $\bar{b}$. Let $C$ be the set containing all the (finitely many) $G$-conjugates of $\bar{c}$, then $C$ codes $A_{M}$ which finishes the proof.
\end{proof}

\begin{corollary}\label{corollaryEIuptofinitesetsoffinitetuples}
    The theory $\mathrm{ACFA}^{+}_{\mathrm{EI}}$ has elimination of imaginaries up to finite sets of finite tuples.
\end{corollary}
\begin{proof}
    This follows directly from Theorem \ref{theoremweakEIACFAplusEI} together with Lemma \ref{lemmacodingofcompactsetsacfator}.
\end{proof}

\section{Further remarks and questions}\label{sectionfurtherremarks}
\subsection{Higher amalgamation}
Given the interesting behaviour of 3-amalgamation in $\mathrm{ACFA}^{+}$, it is natural to ask whether higher dimensional phenomena also occur. In the reducts $\mathrm{ACFA}$ and $\mathrm{PF}^{+}$ the proof for 3-amalgamation can be generalised straight-forwardly to higher dimensions to show that $n$-amalgamation holds for all $n\in\mathbb{N}$ over all algebraically closed sets in both cases.
In a forthcoming paper, we will show that, however, higher amalgamation in $\mathrm{ACFA}^{+}$ behaves more complicatedly. We will give two results. On one hand, we show that 4-amalgamation is strictly stronger than 3-amalgamation by constructing a $\sigma$-AS-closed set over which 4-amalgamation does not hold. 
On the other hand, we will prove, using a stability-theoretic argument, that over models $n$-amalgamation still holds for all $n\in\mathbb{N}$. This differs from the rank 1 context where it is easy to deduce from a non-trivial obstruction to 3-amalgamation that 4-amalgamation does not holds over models. As both proofs turn out to be surprisingly technical and involved we decided to present them in a subsequent paper.

\subsection{$\mathrm{ACFA}$ with additive and multiplicative character on the fixed field}
In \cite{ludwig2025pseudofinitefieldsadditivemultiplicative} we extend Hrushovski's work on the theory $\mathrm{PF}^{+}$ in \cite{Hrushovski2021AxsTW} to the theory $\mathrm{PF}^{+,\times}$ of pseudofinite fields with (sufficiently generic) additive and multiplicative character. This also entails the purely multiplicative case $\mathrm{PF}^{\times}$. It is natural to expect that the results of this paper have a multiplicative analogue: We consider for the moment the theory $\mathrm{ACFA}^{\times}$ which consists of $\mathrm{ACFA}$ together with axioms stating that for the fixed field with a multiplicative character $(F,\chi)$ we have $(F,\chi)\models \mathrm{PF}^{\times}$. Quantifier elimination should work exactly as for $\mathrm{ACFA}^{+}$ and it can be expected that equations of the form $\sigma(x)/x=a$ for $a\in F$ play the same role as the $\sigma$-AS-equations $\sigma(x)-x=a$ for $\mathrm{ACFA}^{+}$ and that thus a similar description of the Kim-Pillay groups can be obtained.
If we now consider both cases together and want to determine the Kim-Pillay group of $\mathrm{ACFA}^{+,\times}$ it is likely that it is simply given by the direct sum of the Kim-Pillay groups of the additive and multiplicative reducts (taking into account the way $\mathrm{PF}^{+,\times}$ is axiomatised in \cite{ludwig2025pseudofinitefieldsadditivemultiplicative}).\\
In terms of higher amalgamation it is however less clear what to expect. In particular, given that the above mentioned counterexample to 4-amalgamation over $\sigma$-AS-closed sets involves the action of the multiplicative group on the additive one, it would be already interesting to know whether 4-amalgamation does hold over the multiplicative analogue of a $\sigma$-AS-closed sets in $\mathrm{ACFA}^{\times}$.

\section{Appendix I: Model-theoretic Galois groups in continuous logic}\label{sectionmodeltheoreticgaloisgroups}

In this section we will study the Kim-Pillay group and its connected component of the identity for a continuous logic theory $T$. Some references on model-theoretic Galois groups in the classical context are among others \cite{CategoryofModelsLascar},\cite{Hrushovski1998SimplicityAT} ,\cite{hyperautomgroups} and \cite{galoisgroupsfirstordertheories}. Hyperdefinable sets in or with connections to continuous logic appear for example in \cite{fernandez2024canonicalquotientsmodeltheory}, \cite{piecewisehyperdefgroups} and \cite{hrushovski2022lascargroup}. The results in this section are certainly known to experts but we did not find a presentation entirely containing the statements needed for our purposes. We thus give a self-contained account, mainly generalising ideas from \cite{Hrushovski1998SimplicityAT} to the continuous context.\\

Again for this section, we fix $\mathcal{M}$ to be a monster model of $T$. We will denote by $X$ a type-definable (over $\emptyset$) set in $M^{m}$. 

\begin{definition}\label{definitionrelativedefinablesubset}(See Definition 2.9 of \cite{strongminsetshanson})
    We say that some subset $D\subseteq X^{n}$ is \textit{definable relative to $X$} if there is a predicate $P:M^{mn}\rightarrow[0,1]$ definable in $\mathcal{M}$ such that $P(\bar{x})=0$ for all $\bar{x}\in D$ and $\forall \bar{x}\in X^{n}$ we have $\mathrm{dist}(\bar{x},D)\leq P(\bar{x})$.
\end{definition}

This notion of a relative definable set already appears in the preprint \cite{strongminsetshanson} (see Definition 2.9). Note that, as mentioned there, it is in general not equivalent that $D$ is definable relative to $X$ to saying that there is some definable predicate $P$ that coincides with the distance predicate of $D$ on $X$. 


        


\begin{definition}
An equivalence relation $E$ on $X$ is called \textit{definable} (or \textit{strongly definable}) if $\{(\bar{x},\bar{y})\in X^{2}\;|\;\bar{x}\,E\,\bar{y}\}$ is a definable set relative to $X$. Further, $E$ is called \textit{weakly definable} if there is a definable predicate $\phi$ such that for all $\bar{x},\bar{y}\in X$ $\phi(\bar{x},\bar{y})=0$ holds if and only if $\bar{x}\,E\,\bar{y}$. Finally, $E$ is called \textit{type-definable} if there is some partial type $p$ such that for all $\bar{x},\bar{y}\in X$ the relation $\bar{x}\,E\,\bar{y}$ holds if and only if $p(\bar{x},\bar{y})$ holds.
\end{definition}

\begin{remark}\label{typedefinablesets}
A weakly definable equivalence relation is a definable equivalence relation in the language used by Chavarria and Pillay in \cite{Gomez2021OnPE}. A type-definable equivalence relation corresponds to the definition of type-definable sets as (possibly infinite) intersections of zerosets.  
\end{remark}

\begin{example}\label{torsorexample}
Let $G$ be a strongly definable group in a model $M\models T$ (i.e., the graph of multiplication is a strongly definable set).
Let $\Psi: G\rightarrow H$ be a definable predicate to some compact abelian group $H$ that is a group homomorphism. Denote by $1_{H}$ the identity element of $H$. Let $X$ be a strongly definable $G$-Torsor in $M^{n}$ for some $n$. The following is a weakly definable equivalence relation on $X$:
\[x\,E\,y\;\iff\;\Psi(xy^{-1})=1_{H}.\]
\end{example}

\begin{lemma}\label{lemmatypedeffinclassesisdefinable}
    A weakly definable equivalence relation $E(\bar{x},\bar{y})$ on $X$ with finitely many classes is definable on $X$. Moreover it is \textit{discrete} on $X$, i.e., there is some $\epsilon>0$ such that the different classes all have distance greater than $\epsilon$ from each other.
\end{lemma}

\begin{proof}
    Let $\phi(\bar{x},\bar{y})$ be such that $E(\bar{x},\bar{y})$ is defined by the zeroset $D$ of $\phi(\bar{x},\bar{y})$ on $X$. We assume that $D$ is not strongly definable on $X$. We find a sequence $(\bar{a}_{n},\bar{b}_{n})$ in $X$ and $\epsilon>0$ such that $\phi(\bar{a}_{n},\bar{b}_{n})\leq 1/n$ but $\mathrm{dist}((\bar{a}_{n},\bar{b}_{n}),E)\geq\epsilon$ for all $n\in\mathbb{N}$. Since there are only finitely many classes, we can pass to a subsequence such that all the $\bar{a}_{n}$ are in the same class $E_{1}$, whereas all the $\bar{b}_{n}$ are in the class $E_{2}$. By saturation we find some $\bar{c},\bar{d}\in D\cap X$ such that $\bar{c}\in E_{1},\,\bar{d}\in E_{2},\, \mathrm{dist}((\bar{c},\bar{d}),E)\geq\epsilon$ and $\phi(\bar{c},\bar{d})=0$ hold, contradicting that $E$ is an equivalence relation on $X$.\\
    The moreover part is proven by a similar argument: Assume that there is no such $\epsilon>0$. We find some class $E_{0}$ and a sequence $(\bar{a}_{i})_{i\in\mathbb{N}}$ in $X$ such that $\bar{a}_{i}\not\in E_{0}$ but $\mathrm{dist}(\bar{a}_{i},E_{0})\leq 1/n$. But then we can assume that $(\bar{a}_{i})_{i\in\mathbb{N}}$ lies in one class $E_{1}$ and by completeness (of the metric space) and continuity of $\phi(\bar{x},\bar{y})$ we obtain a tuple $\bar{c}$ with distance $0$ to $E_{0}$ such that $\phi(\bar{a}_{1},\bar{c})=0$, a contradiction.
\end{proof}

\subsection{Hyperimaginaries/CL-imaginaries}\label{sectionhyperclimaginaries}

\begin{definition}
A hyperimaginary (or $A$-hyperimaginary) is the equivalence class $a_{E}:=\bar{a}/E$ of a (possibly infinite) tuple $\bar{a}$ under a $\emptyset$-type-definable (or $A$-type-definable) equivalence relation $E$. We say that an automorphism $f\in \mathrm{Aut}(M)$ fixes $a_{E}$ if $E(f(\bar{a}),\bar{a})$ holds.
\end{definition}

\begin{definition}
    Let $A$ be a set of hyperimaginaries. We say that a hyperimaginary $a_{E}$ is \textit{definable} over $A$, if it is fixed by all automorphims that fix $A$ pointwise. By the \textit{definable closure} $\mathrm{dcl}^{\mathrm{heq}}(A)$ we denote the class of all hyperimaginaries definable over $A$. We say that two hyperimaginaries $a_{E}$ and $b_{F}$ are \textit{interdefinable} if $\mathrm{Aut}(M/a_{E})=\mathrm{Aut}(M/b_{F})$ (or equivalently if they have the same definable closure).
    The \textit{bounded closure} $\mathrm{bdd}(A)$ is given by all those hyperimaginaries having bounded orbit under $\mathrm{Aut}(M/A)$. Moreover, we write $a_{E}\equiv_{b_{F}}c_{E}$ to express that there is some $f\in \mathrm{Aut}(M/b_{F})$ such that $f(a_{E})=c_{E}$. These notions extend to tuples of hyperimaginaries in the usual way.
\end{definition}

There is a natural way to work with bounded hyperimaginaries in continuous logic allowing for a similar $T^{\mathrm{eq}}$-construction as for imaginaries in classical logic.
The following two facts were already proven by Ben Yaacov in the more general context of Cats in \cite{uncountablecats}. We directly cite the corresponding statement in continuous logic from \cite{Conant2022SeparationFI}.

\begin{fact}\label{factcountequivrelenough}(Fact 4.4 in \cite{Conant2022SeparationFI}.)
    Let $E(\bar{x},\bar{y})$ be a type-definable equivalence relation, and let $(E_{i}(\bar{x},\bar{y}))_{i<\lambda}$ enumerate all countably type-definable equivalence relations (i.e., defined by a countable type) such that $E(\bar{x},\bar{y})$ implies $E_{i}(\bar{x},\bar{y})$. Then for any real tuples $\bar{a}$ and $\bar{b}$ from $M$, $E(\bar{a}, \bar{b})$ holds if and only if $E_{i}(\bar{a},\bar{b})$ holds for all $i<\lambda$.
\end{fact}

\begin{fact}\label{factdefinablepseudometric}(Proposition 4.8 in \cite{Conant2022SeparationFI}.)
Let $E(\bar{x},\bar{y})$ be a countably type-definable equivalence relation. Then there is a definable pseudometric $\rho$ on $M^{|\bar{x}|}$ such that for any $\bar{a},\bar{b}\in M^{|\bar{x}|}$ we have $\rho(\bar{a},\bar{b})=0$ if and only if $E(\bar{a},\bar{b})$ holds.
    
\end{fact}

We only give a brief description of the construction of \textit{CL-imaginaries} and restrict ourselves to finitary imaginaries. For more details, we refer to Section 11 in \cite{mtfms}, Section 5 in \cite{locstabilitycontinuouslogic} and Chapter 3 of \cite{hansonphdthesis}.

\begin{definition}\label{definitionCLimaginaries}
  Let $\rho(\bar{x})$ define a pseudometric in $M^{n}$. Then we consider an element of the metric space given by the completion of $M^{n}/\rho$ as a \textit{CL-imaginary}. We can add the completion of $M^{n}/\rho$ with its metric (similarly to the $T^{\mathrm{eq}}$-construction in classical logic) as a new sort in a richer language as follows. Let $\mathcal{L}_{\rho}$ be the language that extends $\mathcal{L}$ by a new sort that will be interpreted by the completion of $M^{n}/\rho$ and an $n$-ary function symbol $\pi_{\rho}$ that will be interpreted by the canonical quotient mapping from $M^{n}$ to the completion of $M^{n}/\rho$. We denote by $\mathrm{acl}_{\mathrm{CL}}^{\mathrm{eq}}$ the CL-algebraic closure in the corresponding expansion $\mathcal{M}^{\mathrm{eq}}$ of $\mathcal{M}\models T$. (See Section 5 in \cite{locstabilitycontinuouslogic} for a more detailed account).
\end{definition}

\begin{corollary}\label{corollaryclimaginarieshyperimaginaries}
    Every hyperimaginary is interdefinable with a sequence of CL-imaginaries.
\end{corollary}
\begin{proof}
    This follows directly using Fact \ref{factcountequivrelenough} and Fact \ref{factdefinablepseudometric}.
\end{proof}

\begin{fact}\label{factKPequivrelationsequalsboundedequivrelation}(Corollary 1.47 of \cite{simplicityCATs}.)
    For any tuples $\bar{a}$ and $\bar{b}$ of the same length, $\bar{a}\equiv_{A}^{\mathrm{KP}} \bar{b}$ holds if and only if $\bar{a}\equiv_{\mathrm{bdd}(A)}\bar{b}$ holds.
\end{fact}

\begin{fact}\label{factKPgroupasautomgroup}
 For any $A\subseteq M$, the following isomorphies hold. \[\mathrm{Gal}_{\mathrm{KP}}^{A}(M)\cong \mathrm{Aut}_{\mathcal{M}}(\mathrm{bdd}(A)/A)\cong \mathrm{Aut}_{\mathcal{M}}(\mathrm{acl}_{\mathrm{CL}}^{\mathrm{eq}}(A)/A).\]
\end{fact}
\begin{proof}
    The first isomorphy directly follows from the definition, and the second from Corollary \ref{corollaryclimaginarieshyperimaginaries}.
\end{proof}

\subsection{The Kim-Pillay group as a topological group}
\begin{definition}
Let $X$ be an $\emptyset$-type-definable set and $E$ an $\emptyset$-type-definable equivalence relation with a bounded number of classes on $X$. We write $\projection_{E}:X\rightarrow X_{E}:=X/E$ for the projection map and $\pi_{E}(\bar{x},\bar{y})$ for the type defining $E$. We define a subset of $X/E$ to be closed if the preimage (under $\projection_{E}$) is a type-definable set (allowing a small set of parameters). The corresponding topology is called the \textit{logic topology}.
Let $E_{\mathrm{KP}}^{X}$ denote the finest bounded $\emptyset$-type-definable equivalence relation on $X$. We write $X_{\mathrm{KP}}$ for $X/E_{\mathrm{KP}}^{X}$, $\pi_{\mathrm{KP}}$ for $\pi_{E_{\mathrm{KP}}^{X}}$, and $\projection_{\mathrm{KP}}:X\rightarrow X_{\mathrm{KP}}$ for the projection map to the \textit{Kim-Pillay space} $X_{\mathrm{KP}}$. As before $X_{\mathrm{KP}}$ is equipped with the logic topology.
\end{definition}


\begin{lemma}\label{equivalence_seperation}
Let $E$ be a type-definable equivalence relation defined by the partial type $\pi(\bar{x},\bar{y})$. Let $\bar{a},\bar{b}$ be given such that $\lnot E(\bar{a},\bar{b})$ holds. Then there is some $\epsilon_{E}>0$ depending on $\bar{a}$ and $\bar{b}$ as well as a predicate $\psi_{E}(\bar{x},\bar{y})\in \pi$ (with $\psi_{E}(\bar{x},\bar{y})=\psi_{E}(\bar{y},\bar{x})$) such that the sets defined by $\psi_{E}(\bar{a},\bar{x})<\epsilon_{E}$ and $\psi_{E}(\bar{b},\bar{x})<\epsilon_{E}$ have empty intersection. 
\end{lemma}
\begin{proof}
Otherwise, there is by compactness for every $\delta>0$ some $\bar{c}_{\delta}$ such that for all $\phi(\bar{x},\bar{y})\in\pi(\bar{x},\bar{y})$ we have \[\max\{\phi(\bar{a},\bar{c}_{\delta}),\phi(\bar{c}_{\delta},\bar{a}),\phi(\bar{b},\bar{c}_{\delta}),\phi(\bar{c}_{\delta},\bar{b})\}\leq\delta.\]
But then, it follows that in a sufficiently saturated model there is some $\bar{c}$ such that $E(\bar{a},\bar{c})$ and $E(\bar{b},\bar{c})$ hold, which yields a contradiction.
\end{proof}

\begin{lemma}
With respect to the logic topology, the space $X_{\mathrm{KP}}$ is compact and Hausdorff.
\end{lemma}
\begin{proof}
We adapt the proof of Proposition 2.1 in \cite{Hrushovski1998SimplicityAT} to our context using the above.
Compactness of the Kim-Pillay space follows from saturation:
Since $X_{\mathrm{KP}}$ is of bounded size the family of preimages of a family of closed sets with nonempty finite intersections is a family of size $\leq 2^{2^{|T|}}$ of type-definable sets with non-empty finite intersections that has the finite intersection property by saturation. Consequently, the space $X_{\mathrm{KP}}$ has the finite intersection property. The Hausdorff property follows directly from Lemma \ref{equivalence_seperation}.
\end{proof}


\begin{lemma}\label{projoftypedefbsetsareclosed}
    The projections of type-definable sets $Y\subseteq X$ under the map $X\rightarrow X_{\mathrm{KP}}$ are closed.
\end{lemma}
\begin{proof}
    Let $B:=\projection_{\mathrm{KP}}(Y)$ and $\rho(\bar{x})$ the type that defines $Y$. Then the set defined by \[\bigcap_{\psi\in\rho}\;\bigcap_{\gamma\in\pi_{\mathrm{KP}}}\{\bar{x}\;|\;\inf_{\bar{y}}\max\{\rho(\bar{y}),\gamma(\bar{x},\bar{y})\}=0\}
    \]is type-definable and by compactness it equals $\projection_{\mathrm{KP}}^{-1}(B)$.
\end{proof}

\begin{lemma}\label{remarkbasistopology}
Let $A\subseteq X_{\mathrm{KP}}$ be open. Let $\rho(\bar{x})$ be the type that defines the set $Y\subseteq X$ such that $\projection_{\mathrm{KP}}(Y)=X_{\mathrm{KP}}\backslash A$. Then, we have  
\[A=\bigcup_{\psi\in\rho}\bigcup_{\epsilon>0}\{a\in X_{\mathrm{KP}}\,|\,\text{for all } \bar{c}\in\projection_{\mathrm{KP}}^{-1}(a)\,\text{one has}\,\models \psi(\bar{c})>\epsilon\}.\] 
\end{lemma}
\begin{proof}
The statement follows from the following reasoning. For every $a$ in the complement of the projection of the set defined by $\psi(\bar{x})=0$ there is, by compactness, some $\epsilon_{a}$ such that $\forall \bar{c}\in \projection_{\mathrm{KP}}^{-1}(a)$ the inequality $\psi(\bar{c})>\epsilon_{a}$ holds.
\end{proof}

\begin{definition}(Translation of Definition 3.6 of \cite{hyperautomgroups} to the continuous context.)
Let $\pi_{\mathrm{KP}}(\bar{x},\bar{y})$ be the type that defines the equivalence relation $E_{X}^{\mathrm{KP}}$ on the type-definable set $X$ (say defined by the type $p(\bar{x})$). Let $\phi(\bar{x},\bar{y})=0$ be in $\pi_{\mathrm{KP}}(\bar{x},\bar{y})$. Let $a\in X_{\mathrm{KP}}$. We define the set $U_{\phi,a,\epsilon}\subseteq X_{\mathrm{KP}}$ as follows \[U_{\phi,a,\epsilon}:=\{b\in X_{\mathrm{KP}}\,|\,\forall \bar{c}\in\projection_{\mathrm{KP}}^{-1}(b),\forall\bar{d}\in\projection_{\mathrm{KP}}^{-1}(a)\,\text{ one has } \models \phi(\bar{c},\bar{d})<\epsilon\}.\]
    
\end{definition}

\begin{notation}(Definition 6.4 in \cite{mtfms})
We define a binary function $\dotminus:\mathbb{R}_{\geq 0}\times\mathbb{R}_{\geq 0}\rightarrow \mathbb{R}_{\geq 0}$ by
\[x\dotminus y=\begin{cases}
(x-y) & \text{if }\; x \geq y\\
0 & \text{otherwise.}
\end{cases}\]
\end{notation}

\begin{lemma}\label{lemmabasisoftopologymorespecific}
    (This is an adaptation of Lemma 3.7 in \cite{hyperautomgroups} to the continuous context.) The sets of the form $U_{\phi,a,\epsilon}$ form a basis of the topology, i.e., we have that for every open set $U$ with $a\in U$, there are some $\phi$ and $\epsilon>0$ such that $a\in U_{\phi,a,\epsilon}\subseteq U$.
\end{lemma}
\begin{proof}
    Clearly $a\in U_{\phi,a,\epsilon}$. Let $\bar{d}\in\projection_{\mathrm{KP}}^{-1}(a)$. Now, $U_{\phi,a,\epsilon}$ is open by Lemma \ref{projoftypedefbsetsareclosed} as it is the complement of the projection of the type-definable set given by 
    \[\bigcap_{\gamma(\bar{x},\bar{y})\in\pi_{\mathrm{KP}}(\bar{x},\bar{y})}\;\{\bar{x}\;|\;\inf_{\bar{y}}\;\max\{\gamma(\bar{y},\bar{d}),\;\epsilon\dotminus \phi(\bar{x},\bar{y})\}=0\}.\]
    Now let $U$ be a neighbourhood of $a$. It remains to prove that $U_{\phi,a,\epsilon}\subseteq U$ for some $\phi$ and $\epsilon>0$. By Lemma \ref{remarkbasistopology} we can assume that $U=\{b\in X_{\mathrm{KP}}\,|\,\forall\,\bar{c}\in\projection_{\mathrm{KP}}^{-1}(b)\;\xi(\bar{c})>\delta\}$ for some formula $\xi(\bar{x})$ and $\delta>0$. Since $a\in U$, it follows that $\pi_{\mathrm{KP}}(\bar{d},\bar{x})=0$ together with $\xi(\bar{x})\leq\delta$ and the type $p(\bar{x})$ that defines $X$ is inconsistent. Hence, by compactness we find some $\phi(\bar{x},\bar{y})\in\pi_{\mathrm{KP}}(\bar{x},\bar{y})$ and $\epsilon>0$ such that the type $p(\bar{x})$ together with $\phi(\bar{x},\bar{d})<\epsilon$ imply $\xi(\bar{x})>\delta$. Thus, we have $U_{\phi,a,\epsilon}\subseteq U$.
\end{proof}

\begin{definition}\label{definitionlocalKPgroup}
    Let $\mathrm{Stab}(X_{\mathrm{KP}})$ be defined as the subgroup of $\mathrm{Aut}(M)$ given by those automorphisms stabilizing the classes of $X_{\mathrm{KP}}$. We define $G_{\mathrm{KP}}^{X}$ as the group obtained by taking the restriction to $X$ of the elements in the quotient $\mathrm{Aut}(M)/\mathrm{Stab}(X_{\mathrm{KP}})$. Similarly we define $G_{E}^{X}$ for $X_{E}$ where $E$ is an $\emptyset$-type-definable equivalence relation on $X$ with boundedly many classes (but not necessarily the finest). 
\end{definition}
\begin{lemma}\label{lemmalocalKPgroupastopologicalgroup}
    The group $G_{\mathrm{KP}}^{X}$ equipped with the topology induced by the Tychonoff-topology on $X_{\mathrm{KP}}^{X_{\mathrm{KP}}}$ is a compact (Hausdorff) topological group and the action on $X_{\mathrm{KP}}$ is continuous.
\end{lemma}

\begin{proof}
     To show compactness, we follow the proof for the classical context, which is given in Proposition 2.1(4) and Lemma 2.3 of \cite{Hrushovski1998SimplicityAT}. First, we show that the image of $\mathrm{Aut}(M)$ under the projection map $\mathrm{Aut}(M)\rightarrow G_{\mathrm{KP}}^{X}$ in the Tychonoff-product $X_{\mathrm{KP}}^{X_{\mathrm{KP}}}$ is closed. (This is Proposition 2.1(4) in \cite{Hrushovski1998SimplicityAT}.) Given a net of elements $f_{i}$ from the image of $\mathrm{Aut}(M)$ converging point-wise to some $f\in X_{\mathrm{KP}}^{X_{\mathrm{KP}}}$ we want to show that we can find some $g\in \mathrm{Aut}(M)$ that projects to $f$. Let us first note that $f$ preserves the projections of $\emptyset$-definable zerosets:
     Let $D$ be an $\emptyset$-definable zeroset. Then, $D$ is clearly preserved by all the $f_{i}$. But then it also follows for every $d\in D$ that $f(d)\in D$ since $f$ is the limit of pointwise convergence of the $f_{i}$ and $D$ is closed.\\
    To show that there is some $g\in \mathrm{Aut}(M)$ projecting onto $f$, by compactness it suffices to find for every finite set $\{\bar{a}_{1},\dots,\bar{a}_{n}\}\subseteq X$ some $\Tilde{g}\in \mathrm{Aut}(M)$ such that $\Tilde{g}(\bar{a}_{i})/E_{\mathrm{KP}}=f(\bar{a}_{i}/E_{\mathrm{KP}})$. Assume the contrary, i.e., for some given $\bar{a}_{1},\dots,\bar{a}_{n}$ we cannot find some $\bar{b}_{1},\dots, \bar{b}_{n}$ with $\mathrm{tp}(\bar{a}_{1},\dots,\bar{a}_{n})=\mathrm{tp}(\bar{b}_{1},\dots,\bar{b}_{n})$ and $\bar{b}_{i}/E_{\mathrm{KP}}=f(\bar{a}_{i}/E_{\mathrm{KP}})$ for all $1\leq i\leq n$. But then by compactness, we find some zeroset $C\subseteq X^{n}$ with $\bar{a}=(\bar{a}_{1},\dots,\bar{a}_{n})\in C$ such that there is no tuple $\bar{b}=(\bar{b}_{1},\dots, \bar{b}_{n})\in C$ such that $\bar{b}_{i}/E_{\mathrm{KP}}=f(\bar{a}_{i}/E_{\mathrm{KP}})$ for all $1\leq i\leq n$. However, since $f$ preserves projections of zerosets we can lift $f(\bar{a}/E_{\mathrm{KP}})$ to some tuple in $C$ which yields a contradiction. The Hausdorff property of $G_{\mathrm{KP}}^{X}$ is obvious since $X_{\mathrm{KP}}$ and, consequently, $X_{\mathrm{KP}}^{X_{\mathrm{KP}}}$ are Hausdorff. We proceed to proving that the group action $G_{\mathrm{KP}}^{X}\times X_{\mathrm{KP}}\rightarrow X_{\mathrm{KP}}$ is continuous: This will be an adaptation of the second part of the proof of Lemma 3.11 in \cite{hyperautomgroups} to the continuous context using Lemma \ref{lemmabasisoftopologymorespecific}.\\
     Let $a\in X_{\mathrm{KP}}$ and $U$ an open neighbourhood of $a$. Let $b\in X_{\mathrm{KP}}$ be such that $g.b=a$ for some $g\in G_{\mathrm{KP}}^{X}$. We want to find an open neighbourhood of $(g,b)\in G_{\mathrm{KP}}^{X}\times X_{\mathrm{KP}}$ that is sent into $U$ under the action. By Lemma \ref{lemmabasisoftopologymorespecific} we can take $U$ to be of the form $U_{\phi,a,\epsilon}$. Now we find by compactness some $\psi(\bar{x},\bar{y})\in\pi_{\mathrm{KP}}(\bar{x},\bar{y})$ and $\delta>0$ such that $\psi(\bar{x},\bar{z})<\delta$ and $\psi(\bar{y},\bar{z})<\delta$ imply that $\phi(\bar{x},\bar{y})<\epsilon$. We then consider the set $V=U_{\psi,b,\delta}$ and $W=\{g\in G_{\mathrm{KP}}^{X}\,|\,g.b\in U_{\psi,a,\delta}\}$ and note that the latter is an open neighbourhood of $g$ in $G_{\mathrm{KP}}^{X}$.
     Given $(b^{\prime},g^{\prime})\in V\times W$, we have $g^{\prime}.b^{\prime}\in g^{\prime}(V)$. Further, for all $\bar{e}\in \projection_{\mathrm{KP}}^{-1}(g^{\prime}.b^{\prime})$ and for all $\bar{h}\in \projection_{\mathrm{KP}}^{-1}(g^{\prime}.b)$ the inequality $\psi(\bar{e},\bar{h})<\delta$ holds by definition of $V$. Similarly, for all $\bar{h}\in \projection_{\mathrm{KP}}^{-1}(g^{\prime}.b)$ and for all $\bar{f}\in \projection_{\mathrm{KP}}^{-1}(a)$ the inequality $\psi(\bar{h},\bar{f})<\delta$ holds by definition of $W$. Consequently, it follows that $\phi(\bar{e},\bar{f})<\epsilon$ and hence $(g^{\prime}.b^{\prime})\in U_{\phi,a,\epsilon}$, which completes the proof.
     Showing that the multiplication is continuous works in a similar way. Thus, we only give a sketch of the proof without details. First, we note that we can take the open neighbourhoods of elements $g,g^{\prime},g\cdot g^{\prime}$, to be (finite intersections of) sets of the form as $W$ above. Then, assuming that the neighbourhood of $g\cdot g^{\prime}$ is given by $U_{\phi,a,\epsilon}$, we only have to choose $\psi$ and $\delta$ similarly to above and consider the sets $U_{\psi,g^{-1}(a),\delta}$ and $U_{\psi,(gg^{\prime})^{-1}(a),\delta}$ to conclude with a similar argumentation as above.
\end{proof}
\begin{fact}\label{lemmarestrictioneqrelation}(Corollary 1.44 in \cite{simplicityCATs}.)
Let $X\subseteq Y$ be type-definable (over $A$). Then $E_{\mathrm{KP}}^{X}$ is given by $E_{\mathrm{KP}}^{Y}\restriction_{X}$ and, moreover, coincides with being in the same orbit under $\mathrm{Aut}_{\mathrm{KP}}^{A}(M)$.
\end{fact}

\begin{lemma}\label{lemmakimpillaygroupasinverselimit}(4.21 of \cite{hyperautomgroups} in the classical context.) Let $X_{\mathrm{KP}}$ be a Kim-Pillay space as above. The group $\mathrm{Aut}(M)$ acts on $X_{\mathrm{KP}}$, the induced homomorphism $\gamma_{X}:\mathrm{Aut}(M)\rightarrow G_{\mathrm{KP}}^{X}$ is surjective and the intersection of the kernels of all those homomorphisms is given by $\mathrm{Aut}_{\mathrm{KP}}(M)$. The same holds if $X$ varies only over all (finite products of) the sorts. Consequently, $\mathrm{Gal}_{\mathrm{KP}}(M)$ is given as the inverse limit of the groups $G_{\mathrm{KP}}^{S}$ for the $S$ varying over all (finite products of) the sorts.\footnote{Similarly, we can express $\mathrm{Gal}_{\mathrm{KP}}(M)$ as the inverse limit of the groups $G_{\mathrm{KP}}^{X}$ for $X$ type-definable with the directed system coming from projections of sorts and (contravariantly) from inclusions of type-definable sets in the same sort.}   
\end{lemma}

\begin{proof}
    This follows directly using Fact \ref{lemmarestrictioneqrelation}.
\end{proof}

\begin{definition}\label{definitiontopologyonKPgroup}
    Using Lemma \ref{lemmakimpillaygroupasinverselimit} we equip the Kim-Pillay group $\mathrm{Gal}_{\mathrm{KP}}(M)$ with the topology induced from the $G_{\mathrm{KP}}^{X}$'s.
\end{definition}

\begin{corollary}
    With the topology of Definition \ref{definitiontopologyonKPgroup} the Kim-Pillay group $\mathrm{Gal}_{\mathrm{KP}}(M)$ is a compact Hausdorff topological group.
\end{corollary}
\begin{proof}
    This follows directly from Lemma \ref{lemmakimpillaygroupasinverselimit} and Lemma \ref{lemmalocalKPgroupastopologicalgroup} since being Hausdorff and compact is preserved under inverse limits. (See Theorem 3.8, Chapter 8 in \cite{EilenbergSteenrod+1952}.)
\end{proof}
\subsection{The identity component of the Kim-Pillay group}
The aim of this section will be to investigate the connected component of the identity in the Kim-Pillay group and its relation to another natural model-theoretic Galois group.




\begin{definition}
Let $A\subseteq \mathcal{M}\models T$. We say that two (possibly infinite) tuples $\bar{a}$ and $\bar{b}$ have the same \textit{strong type} or \textit{Shelah strong type} (over $A$), written $\bar{a}\equiv^{\mathrm{Sh}}\bar{b}$ (or $\bar{a}\equiv_{A}^{\mathrm{Sh}}\bar{b}$) if $\bar{a}$ and $\bar{b}$ have the same type and are equivalent in any definable (over $A$) equivalence relation with finitely many classes.
We define the group $\mathrm{Aut}_{\mathrm{Sh}}(M)$ (resp. $\mathrm{Aut}_{\mathrm{Sh}}^{A}(M)$) as the subgroup of the automorphism group consisting of those $f\in \mathrm{Aut}(M)$ that fix all equivalence classes of all strongly definable (over $A$) equivalence relations with finitely many classes. The \textit{Shelah-Galois group} is the quotient $\mathrm{Aut}(M)/\mathrm{Aut}_{\mathrm{Sh}}(M)$ (or $\mathrm{Aut}(M/A)/\mathrm{Aut}_{\mathrm{Sh}}^{A}(M)$) and will be denoted by $\mathrm{Gal}_{\mathrm{Sh}}(M)$ (or $\mathrm{Gal}_{\mathrm{Sh}}^{A}(M)$).\\ Equivalently to the Kim-Pillay space, we define on a sort $X$ (or even type-definable set) the \textit{Shelah space} (or space of strong types) $X_{\mathrm{Sh}}$ as $X/E_{\mathrm{Sh}}$ where $E_{\mathrm{Sh}}$ is the intersection of all definable equivalence relations with finitely many classes.
Again, we equip $X_{\mathrm{Sh}}$ with the \textit{logic topology} as for $X_{\mathrm{KP}}$.
\end{definition}

\begin{remark}
    Note that in our definition $E_{\mathrm{Sh}}(\bar{a},\bar{b})$ does not imply $\mathrm{tp}(\bar{a})=\mathrm{tp}(\bar{b})$. This, of course, differs from the purely classical context. If we chose to define $E_{\mathrm{Sh}}$ as being in the same $\mathrm{Aut}_{\mathrm{Sh}}(M)$-orbit, which arguably might seem more natural, $X_{\mathrm{Sh}}$ would no longer necessarily be totally disconnected. This property is what we will want to use. See, however, Remark \ref{remarkdiscussiononShelahspace} for a further discussion of the consequences of our choice of the definition.
\end{remark}

\begin{definition}
An \textit{imaginary} (or sometimes \textit{classical imaginary}) over $A$ is the equivalence class $a_{E}:=\bar{a}/E$ of a tuple $\bar{a}$ under a definable\footnote{Recall that by our convention this is meant in this strong sense, i.e., $E$ defines a definable set.} (over $A$) equivalence relation $E$. We say that an automorphism $f\in \mathrm{Aut}(M)$ fixes $a_{E}$ if $E(f(\bar{a}),\bar{a})$ holds. We denote by $\mathrm{acl}_{\mathrm{DL}}^{\mathrm{eq}}(A)$ the set of all $a_{E}$ where $E$ ranges over all $A$-definable equivalence relations with finitely many classes.
\end{definition}

Directly from the definition, we now obtain the following.
\begin{fact}\label{factShelahgroupasautomgroup}
  For any $A\subseteq M$, we have $\mathrm{Gal}_{\mathrm{Sh}}^{A}(M)\cong \mathrm{Aut}_{\mathcal{M}}(\mathrm{acl}_{\mathrm{DL}}^{\mathrm{eq}}(A)/A)$. 
\end{fact}

\begin{definition}
    Let $\epsilon>0$ and $\phi(\bar{x},\bar{y})$ be a formula in two tuples of variables $\bar{x},\bar{y}$ of the same length. We say that $\phi(\bar{x},\bar{y})$ is $\epsilon$\textit{-thick} if there is no model $\mathcal{M}\models T$ and no infinite $\epsilon$\textit{-antichain} in $\mathcal{M}$, i.e., if there is no infinite sequence $(\bar{a}_{i})_{i\in\omega}$ such that $\max\{\phi(\bar{a}_{i},\bar{a}_{j}),\phi(\bar{a}_{j},\bar{a}_{i})\}>\epsilon$ for all $i<j\in\omega$. Similarly, we say that $\phi(\bar{x},\bar{y})$ is $\epsilon$-thick on $X$, if there is no $\epsilon$-antichain for $\phi(\bar{x},\bar{y})$ in $X$.
\end{definition}
\begin{convention}
    Let $E(\bar{x},\bar{y})$ be a bounded type-definable equivalence relation on $X$ defined by a type $\pi_{E}(\bar{x},\bar{y})$. For $\phi(\bar{x},\bar{y}),\psi(\bar{x},\bar{y})\in\pi_{E}(\bar{x},\bar{y})$ we say that $\psi(\bar{x},\bar{y})$ \textit{refines} $\phi(\bar{x},\bar{y})$ on $Y\supseteq X$, if $\psi(\bar{x},\bar{y})\geq \phi(\bar{x},\bar{y})$ on $Y^{2}$.
\end{convention}
The following is an adaptation of the methods in Section 2.1 of \cite{Hrushovski1998SimplicityAT} to the continuous context.

\begin{lemma}\label{lemmaformulasinKPtypearethick}
Let $E(\bar{x},\bar{y})$ be a bounded type-definable equivalence relation (defined by the type $\pi_{E}(\bar{x},\bar{y})$) on a type-definable set $X$. For every $\phi(\bar{x},\bar{y})\in \pi_{E}(\bar{x},\bar{y})$ and every $\epsilon>0$ we find a zeroset $Y\supseteq X$ and a finite conjunction\footnote{By a conjunction of $\phi_{1}(\bar{x}),\dots,\phi(\bar{x})$ here we simply mean the formula $\psi(\bar{x}):=\max\{\phi_{1}(\bar{x}),\dots,\phi(\bar{x})\}$.} $\psi(\bar{x},\bar{y})$ of formulas of $\pi_{E}(\bar{x},\bar{y})$ such that $\psi(\bar{x},\bar{y})$ is $\epsilon$-thick on $Y$ and $\psi(\bar{x},\bar{y})$
refines $\phi(\bar{x},\bar{y})$. 
\end{lemma}
\begin{proof}
Assume not. Then, by compactness it already follows that $E(\bar{x},\bar{y})$ is not bounded.
\end{proof}

We now adapt a classical strategy working with thick formulas to our context. (See for example Chapter 9 in \cite{casanovas_2011} for similar arguments in the classical context.)
\begin{lemma}\label{lemmatransitiveclosureofthickdefinable}
   Let $Y$ be a zeroset (defined over $\emptyset$). Let $\phi(\bar{x},\bar{y})$ be an $\epsilon$-thick formula on $Y$ (defined over $\emptyset$) such that $\sup_{\bar{x}}\phi(\bar{x},\bar{x})\leq\epsilon$ holds on $Y$. Then the transitive closure on $Y$ of the relation $D(\bar{x},\bar{y})$ given by $D(\bar{x},\bar{y}):\iff\max\{\phi(\bar{x},\bar{y}),\phi(\bar{y},\bar{x})\}\leq\epsilon$, is a definable equivalence relation relative to $Y$ (defined over $\emptyset$) with finitely many classes.
\end{lemma}
\begin{proof}
    Let $n\in\mathbb{N}$ be the maximal length of an $\epsilon$-antichain for $\phi(\bar{x},\bar{y})$ on $Y$. We write $Y(\bar{x})$ for a formula with zeroset $Y$. Let $R(\bar{x},\bar{y})$ be the formula given by $\max\{\phi(\bar{x},\bar{y}),\phi(\bar{y},\bar{x})\}\dotminus\epsilon$. Let $\rho(\bar{x},\bar{y})$ be the formula defined as follows.
    \[\rho(\bar{x},\bar{y})=\inf_{\bar{z}_{1},\dots,\bar{z}_{2n}}\;\max_{1\leq i\leq 2n-1}\{R(\bar{x},\bar{z}_{1}),R(\bar{z}_{i},\bar{z}_{i+1}),R(\bar{z}_{2n},\bar{y}),Y(\bar{z}_{i}),Y(\bar{z}_{2n})\}.\]
    We will show that the zeroset defined by $\rho(\bar{x},\bar{y})$ yields an equivalence relation with finitely many classes on $Y$. Then, by Lemma \ref{lemmatypedeffinclassesisdefinable} this equivalence relation is definable relative to $Y$. The relation given by the zeroset of $\rho(\bar{x},\bar{y})$ on $Y$ is clearly symmetric and reflexive, as is $R(\bar{x},\bar{y})$. Since it is included in the transitive closure of $D(\bar{x},\bar{y})$ on $Y$, it suffices to show that it is transitive to complete the proof. Assume the contrary, i.e., we find $\bar{a},\bar{b},\bar{c}\in Y$ with $\rho(\bar{a},\bar{b})=\rho(\bar{b},\bar{c})=0$ but $\rho(\bar{a},\bar{c})>0$.
    Let $(\alpha_{i})_{0\leq i\leq m-1}$ be a sequence in $Y$ of minimal length such that $\alpha_{0}=\bar{a}, \alpha_{m}=\bar{c}$ and $\alpha_{i}=\bar{b}$ for some $0<i<m-1$ and $D(\alpha_{i},\alpha_{i+1})$ holds for all $0\leq i\leq m-2$. Such a sequence clearly exists since $\rho(\bar{a},\bar{b})=\rho(\bar{b},\bar{c})=0$ (and by saturation it can be found in $Y$). Moreover, from $\rho(\bar{a},\bar{c})>0$ it follows that $m\geq 2n+3$.\\
    W.l.o.g assume $m=2n+3$. We define the sequence $(\beta_{i})_{0\leq i\leq n}$ as $\beta_{i}=\alpha_{2i}$ for $0\leq i\leq n-1$ and $\beta_{n}=\alpha_{2n+3}$. Then $(\beta_{i})_{0\leq i\leq n}$ is an $\epsilon$-antichain for $\phi(\bar{x},\bar{y})$ because of the following reasoning. If it was not an $\epsilon$-antichain for $\phi(\bar{x},\bar{y})$, then we would find some $\beta_{j},\beta_{k}$ with $k<j$ and $\max\{\phi(\beta_{j},\beta_{k}),\phi(\beta_{k},\beta_{j}))\}\leq\epsilon$. But then the sequence that we obtain from $(\alpha_{i})_{0\leq i\leq m-1}$ by deleting every element between $\beta_{k}$ and $\beta_{j}$ contradicts the minimality of $m$. Thus, $(\beta_{i})_{0\leq i\leq n}$ is an $\epsilon$-antichain for $\phi(\bar{x},\bar{y})$ on $Y$. Since $(\beta_{i})_{0\leq i\leq n}$ has $n+1$ elements, we get a contradiction to the assumption that the maximal length of an $\epsilon$-antichain for $\phi(\bar{x},\bar{y})$ was $n$.
\end{proof}

The next two lemmas appear in \cite{Hrushovski1998SimplicityAT} (see Lemma 2.1 and Lemma 2.2) in the classical context. We adapt the proofs to the continuous one.
\begin{lemma}\label{lemmaconnectedcomponentofKPspace}
  Let $X$ be type-definable. There is a canonical map $\zeta: X_{\mathrm{KP}}\rightarrow X_{\mathrm{Sh}}$ that is continuous and surjective and such that the fibres are exactly given by the connected components of $X_{\mathrm{KP}}$.
\end{lemma}

\begin{proof}
    Surjectivity is clear as $E_{\mathrm{KP}}$ refines $E_{\mathrm{Sh}}$. The logic topology on $X_{\mathrm{Sh}}$ is exactly the topology induced by the map $\zeta: X_{\mathrm{KP}}\rightarrow X_{\mathrm{Sh}}$, hence $\zeta$ is continuous. Moreover $X_{\mathrm{Sh}}$ is a profinite space, and thus it is totally disconnected. By continuity, images of connected subsets are connected, thus by total disconnectedness given by a point in $X_{\mathrm{Sh}}$. Now it remains to show that if $X_{\mathrm{KP}}=X_{\mathrm{KP},1}\dot{\cup}X_{\mathrm{KP},2}$ for $X_{\mathrm{KP},1},X_{\mathrm{KP},2}$ closed, then the images of $X_{\mathrm{KP},1}$ and $X_{\mathrm{KP},2}$ under the projection are disjoint. Furthermore, let $X_{1},X_{2}$ be the preimages of $X_{\mathrm{KP},1},X_{\mathrm{KP},2}$ under $\mathrm{proj}_{\mathrm{KP}}$.
    The first part of the argument will simply follow from compactness. However, since there are many parameters to take care of, we subdivide the argument into the following steps:
    \begin{enumerate}[(A)]
        \item There exist $\gamma\geq 0,\,\delta>\gamma$, a formula $\phi(\bar{x})$ (possibly defined with parameters) and a zeroset $Z\supseteq X$ (defined over $\emptyset$) such that $Z=Z_{1}\dot{\cup}Z_{2}$ for zerosets $Z_{1},Z_{2}$ (possibly defined with parameters) and for all $\bar{z}_{1}\in Z_{1}$ and $\bar{z}_{2}\in Z_{2}$, we have $\phi(\bar{z}_{1})\leq \gamma$ and $\phi(\bar{z}_{2})\geq\delta$. Moreover $X\cap Z_{1}=X_{1}$ and $X\cap Z_{2}=X_{2}$.
        \item There exist $\epsilon>0$ and $\Tilde{\Theta}(\bar{x},\bar{y})\in\pi_{\mathrm{KP}}(\bar{x},\bar{y})$ such that $\Tilde{\Theta}(\bar{y}_{1},\bar{y}_{2})\geq \epsilon$ for all $\bar{y}_{1}\in Z_{1}$ and $\bar{y}_{2}\in Z_{2}$.
        \item (Given $\Tilde{\Theta}$ and $\epsilon$.) There exist a zeroset $Y_{\Tilde{\Theta},\epsilon}\supseteq X$ and $\Theta(\bar{x},\bar{y})\in\pi_{\mathrm{KP}}(\bar{x},\bar{y})$ that refines $\Tilde{\Theta}(\bar{x},\bar{y})$ on $Y_{\Tilde{\Theta},\epsilon}$ so that $\Theta(\bar{x},\bar{y})$ is $\epsilon$-thick on $Y_{\Tilde{\Theta},\epsilon}$.
        \item There exist $\epsilon>0$, a zeroset $Y\supseteq X$ and a formula $\Theta(\bar{x},\bar{y})\in\pi_{\mathrm{KP}}(\bar{x},\bar{y})$ such that the following holds:
        \begin{itemize}
            \item $Y=Y_{1}\dot{\cup}Y_{2}$ for zerosets $Y_{1}, Y_{2}$ with $X\cap Y_{1}=X_{1}$ and $X\cap Y_{2}=X_{2}$.
            \item $\Theta(\bar{x},\bar{y})$ is $\epsilon$-thick on $Y$.
            \item $\Theta(\bar{y}_{1},\bar{y}_{2})\geq \epsilon$ for all $\bar{y}_{1}\in Y_{1}$ and $\bar{y}_{2}\in Y_{2}$.
            \item $\Theta(\bar{x},\bar{x})\leq\epsilon$ on $Y$.
        \end{itemize}
    \end{enumerate}
\textit{Proof of} (A): From $X=X_{1}\dot{\cup}X_{2}$ we obtain $\phi(\bar{x})$ and $\tau>0$ such that for all $\bar{x}_{1}\in X_{1}$ and $\bar{x}_{2}\in X_{2}$ we have $\phi(\bar{x}_{1})=0$ and $\phi(\bar{x}_{2})\geq \tau$. Then, by compactness we find some zeroset $Z\supseteq X$ as well as $0\leq\gamma<\delta\leq\tau$ such that $Z$ is the disjoint union of the sets defined by $\phi(\bar{x})\leq\gamma$ and by $\phi(\bar{x})\geq\delta$.\\ 
\textit{Proof of} (B):
This follows directly from compactness. (Possibly after intersecting the $Z$ from (A) with some other zeroset $Y^{\prime}\supseteq X$.)\\
\textit{Proof of} (C): This is Lemma \ref{lemmaformulasinKPtypearethick}.\\
\textit{Proof of} (D): Take $Z, \epsilon>0$ and $\Tilde{\Theta}(\bar{x},\bar{y})$ as given by (B). Then, apply (C) to obtain $Y_{\Tilde{\Theta},\epsilon}$ and $\Theta(\bar{x},\bar{y})$ and set $Y$ to be given by $Y:=Z\cap Y_{\Tilde{\Theta},\epsilon}$. Similarly, we set $Y_{1}:=Z_{1}\cap Y_{\Tilde{\Theta},\epsilon}$ and $Y_{2}:=Z_{2}\cap Y_{\Tilde{\Theta},\epsilon}$. If necessary, we intersect with the zeroset defined by $\Theta(\bar{x},\bar{x})\leq\epsilon$.\\
We continue with the proof of the lemma with the data obtained from (D).
Let $F(\bar{x},\bar{y})$ be the transitive closure of the relation $\Theta(\bar{x},\bar{y})\leq \epsilon$ on $Y$, then as $\Theta(\bar{x},\bar{y})$ is defined over $\emptyset$, by Lemma \ref{lemmatransitiveclosureofthickdefinable} $F(\bar{x},\bar{y})$ is defined over $\emptyset$ as well. (Here we use that $Y$ is defined over $\emptyset$, whereas $Y_{1}, \,Y_{2}$ not necessarily.) Moreover, on the one hand, $F(\bar{x},\bar{y})$ has only finitely many classes on $Y$ since otherwise $\Theta(\bar{x},\bar{y})$ could not be $\epsilon$-thick. On the other hand, for all $\bar{a}\in Y_{1}$ and $\bar{b}\in Y$ we find that $\Theta(\bar{a},\bar{b})\leq\epsilon$ implies that $\bar{b}\in Y_{1}$ and similarly for $\bar{a}\in Y_{2}$. It follows that an element of $X_{1}$ cannot be in the same $F$-class as an element of $X_{2}$ and vice versa, which completes the proof.
\end{proof}

\begin{lemma}\label{lemmaequivalencetocompactopentopology}
   If $C, C^{\prime}\subseteq X_{\mathrm{KP}}$ are closed, then the subset $K_{C,C^{\prime}}\subseteq G_{\mathrm{KP}}^{X}$ given by $K_{C,C^{\prime}}:=\{g\in G_{\mathrm{KP}}^{X}\,|\,g(C)\cap C^{\prime}=\emptyset\}$ is open, and moreover the sets $K_{C,C^{\prime}}$ form a basis of the topology. In other words, $G_{\mathrm{KP}}^{X}$ carries the compact-open topology.
\end{lemma}
\begin{proof}
    Let $C, C^{\prime}\subseteq X_{\mathrm{KP}}$ be as in the statement. If $K_{C,C^{\prime}}=\emptyset$ there is nothing to show. Thus, we may suppose that there is $g\in K_{C,C^{\prime}}$. By compactness, we find zerosets $\mathcal{D}$ and $\mathcal{D}^{\prime}$ (defining subsets of $X$) such that their projections $D$ and $D^{\prime}$ to $X_{\mathrm{KP}}$ include $C,C^{\prime}$, respectively, and $g(D)\cap D^{\prime}=\emptyset$ holds. To prove the theorem, it will suffice to find an open neighbourhood of $g$ with the same property. Consider a lift $\Tilde{g}$ of $g$ to $X$. 
    For any formula $\Theta(\bar{x},\bar{y})$ and $\epsilon>0$ let $S_{\Theta,\epsilon}$ denote the zeroset (in $X$) of the formula $\phi(\bar{x})$ given by
    \[\phi_{\Theta,\epsilon}(\bar{x})=\inf_{\bar{y}}\inf_{\bar{z}}\left( (\Theta(\bar{x},\bar{y})\dotminus\epsilon)+(\Theta(\bar{y},\bar{z})\dotminus\epsilon)+\mathcal{D}^{\prime}(\bar{z}) \right). \]
    By compactness, we find some $\epsilon>0$ and $\Theta(\bar{x},\bar{y})\in\pi_{\mathrm{KP}}(\bar{x},\bar{y})$ such that $\Tilde{g}(\mathcal{D})$ is disjoint from $S_{\Theta,\epsilon}$ and by Lemma \ref{lemmaformulasinKPtypearethick} we can assume that $\Theta(\bar{x},\bar{y})$ is $\epsilon$-thick. Let $\{\Tilde{d}_{1},\dots,\Tilde{d}_{n}\}\subseteq\mathcal{D}$ be a maximal $\epsilon$-antichain of $|\bar{x}|$-tuples for $\Theta(\bar{x},\bar{y})$ in $\mathcal{D}$ and $d_{i}$ the image of $\Tilde{d}_{i}$ under the projection $X\rightarrow X_{\mathrm{KP}}$ for any $1\leq i\leq n$. Let $\Tilde{U}\subseteq \mathrm{Aut}(\mathcal{M})$ be the set given as follows. 
    \[\Tilde{U}:=\{\Tilde{f}\in \mathrm{Aut}(\mathcal{M})\;|\;\Tilde{f}(\Tilde{d}_{i})\not\in S_{\Theta,\epsilon}\text{ for all }1\leq i\leq n\}.\] Let $U$ be the image of $\Tilde{U}$ under the projection $\mathrm{Aut}(\mathcal{M})\rightarrow G_{\mathrm{KP}}^{X}$. Then $U$ is an open neighbourhood of $g$.
    Finishing the proof now reduces to showing that for any $f\in U$ we have $f(D)\cap D^{\prime}=\emptyset$. Let $f\in U$ be arbitrary and $\Tilde{f}\in\Tilde{U}$ a lift of $f$ to $X$. Let $\Tilde{d}$ be a tuple in $\mathcal{D}$, then $\Theta(\Tilde{d},\Tilde{d}_{i})\leq\epsilon$ for some $1\leq i\leq n$ and consequently $\Theta(\Tilde{f}(\Tilde{d}),\Tilde{f}(\Tilde{d}_{i}))\leq\epsilon$. In the case where $\Theta(\Tilde{f}(\Tilde{d}),\Tilde{e})\leq\epsilon$ for some tuple $\Tilde{e}\in\mathcal{D}^{\prime}$, it would follow that $\Tilde{f}(\Tilde{d}_{i})\in S_{\Theta,\epsilon}$, which yields a contradiction. Since $\Tilde{d}\in \mathcal{D}$ and $\Tilde{e}\in\mathcal{D}^{\prime}$ were arbitrary, it follows that $f(D)\cap D^{\prime}=\emptyset$ for any $f\in U$.
    The moreover part then follows using Theorem 2 in \cite{topologyforspacesoftransformations} (and the fact that the natural action $G_{\mathrm{KP}}^{X}\times X_{\mathrm{KP}}\rightarrow X_{\mathrm{KP}}$ is continuous as we established in Lemma \ref{lemmalocalKPgroupastopologicalgroup}).
\end{proof}    

We will use the above to characterise the Shelah Galois group with respect to the Kim-Pillay Galois group topologically. To do so, we will need the following basic lemma on topological groups.
\begin{lemma}\label{lemmaconnectedcoponenttoplogicalcharacterisation}
    Any $f\in G_{\mathrm{KP}}^{X}$ is in the connected component of the identity if and only if $f$ preserves all connected components of $X_{\mathrm{KP}}$.
\end{lemma}
\begin{proof}
    Assume that $f$ does not preserve the connected components, then we find some $a,b\in X_{\mathrm{KP}}$ and a clopen set $A\subseteq X_{\mathrm{KP}}$ with $a\in A$ and $b\in X_{\mathrm{KP}}\backslash A$ such that $f(a)=b$. Now the set $\{g\in G_{\mathrm{KP}}^{X}\,|\,g(a)\in X_{\mathrm{KP}}\backslash A\}$ is clopen by Lemma \ref{lemmaequivalencetocompactopentopology} and witnesses that $f$ is not in the connected component of the identity.
    For the other direction, note that a clopen set in $X_{\mathrm{KP}}$ projects to a clopen set in every product over a finite subset of the coordinates where we consider $G_{\mathrm{KP}}^{X}\subseteq X_{\mathrm{KP}}^{X_{\mathrm{KP}}}$. (The projections are open since $G_{\mathrm{KP}}^{X}$ carries the product topology and closed since $G_{\mathrm{KP}}^{X}$ is compact.) Any such clopen in a finite product is a finite union of clopen boxes. 
    Let $B\subseteq G_{\mathrm{KP}}^{X}$ be a clopen set that witnesses that $f$ is not in the connected component of the identity, that is, $B$ contains the identity and $f\notin B$. By compactness $B$ is a finite union of basic open sets (for the product topology). Using the above, we can then assume that $B$ is a positive Boolean combination of sets of the form $
    \{g\in G_{\mathrm{KP}}^{X}\,|\,g(a)\in A\}$ for some clopen $A\subseteq X_{\mathrm{KP}}$ and $a\in A$ (here we use that $B$ contains the identity). It follows that there is some clopen $A\subseteq X_{\mathrm{KP}}$ such that
    $f(a)\notin A$ for some $a\in A$ and thus $f$ does not preserve all connected components.
\end{proof}
\begin{lemma}\label{lemmarestrictiondefequivcomesfromequiv}
Let $X$ be a definable set in a sort $S$ with definable complement. Then $E_{\mathrm{Sh}}^{X}$ is given by $E_{\mathrm{Sh}}^{S}\restriction_{X}$. Moreover, if a sort $S^{\prime}$ projects onto $S$ via $\pi: S^{\prime}\rightarrow S$, then $E_{\mathrm{Sh}}^{S^{\prime}}$ refines (or equals) $\pi^{-1}(E_{\mathrm{Sh}}^{S})$.
\end{lemma}
\begin{proof}
The first part is immediate by extending every definable equivalence relation with finitely many classes by one class outside of $X$. For the moreover part, simply note that for $E^{S}$ definable with finitely many classes on $S$, $\pi^{-1}(E^{S})$ is definable and has finitely many classes as well.
\end{proof}

\begin{theorem}\label{theoremlocalconnectedcomponentshelahgroup}
    Let $X$ be type-definable. Consider the action of $G_{\mathrm{KP}}^{X}$ on $X_{\mathrm{Sh}}$ induced by the map $\zeta: X_{\mathrm{KP}}\rightarrow X_{\mathrm{Sh}}$ (together with the action of $G_{\mathrm{KP}}^{X}$ on $X_{\mathrm{KP}}$). Let $H_{X}$ be the point-wise stabiliser of $X_{\mathrm{Sh}}$ under this action. Then $G_{\mathrm{Sh}}^{X}\cong G_{\mathrm{KP}}^{X}/H_{X}$ and $H_{X}$ is the connected component of the identity in $G_{\mathrm{KP}}^{X}$ and consequently $G_{\mathrm{Sh}}^{X}$ the largest totally disconnected quotient of $G_{\mathrm{KP}}^{X}$.
\end{theorem}

\begin{proof}
By Lemma \ref{lemmaconnectedcomponentofKPspace} we have that the group $H_{X}$ is the largest subgroup of $G_{\mathrm{KP}}^{X}$ that stabilises the connected components of $X_{\mathrm{KP}}$, then by Lemma \ref{lemmaconnectedcoponenttoplogicalcharacterisation} it is the connected component of the identity in $G_{\mathrm{KP}}^{X}$. 
\end{proof}


\begin{remark}\label{remarkdiscussiononShelahspace}
Lemma \ref{lemmarestrictiondefequivcomesfromequiv} can already fail for zerosets $X$. We provide an example using the theory $\mathrm{PF}^{+}$ that will only be introduced in the next chapter. In a sufficiently saturated model $\mathcal{M}\models \mathrm{PF}^{+}$ consider for $r_{1}\neq r_{2}\in S^{1}$ the zeroset $X_{r_{1},r_{2}}:=\{x\in M\,|\,\Psi(x)\in\{r_{1},r_{2}\}\}$. Then the preimages of $r_{1}$ and of $r_{2}$ under $\Psi$ are both clopen sets in $X_{r_{1},r_{2}}$ for every pair $(r_{1},r_{2})$, i.e., separated by a definable equivalence relation with finitely many classes (relative to $X_{r_{1},r_{2}}$). However, this cannot be the case in $M$ since otherwise the character $\Psi$ would be definable in classical first-order logic which is not the case as explained at the end of Section 4 in \cite{Hrushovski2021AxsTW}.
\end{remark}

\begin{corollary}\label{corollaryconnectedcomponentasinverselimit}
We consider the directed systems of all finite products of the sorts $S$ together with all projections. Let $H_{\mathrm{KP}}(M)$ be the inverse limit of the $H_{S}$ (as defined in Theorem \ref{theoremlocalconnectedcomponentshelahgroup}). Then, $\mathrm{Gal}_{\mathrm{KP}}(M)/H_{\mathrm{KP}}(M)\cong \mathrm{Gal}_{\mathrm{Sh}}(M)$ and $H_{\mathrm{KP}}(M)$ is the connected component of the identity in $\mathrm{Gal}_{\mathrm{KP}}(M)$.\footnote{Similarly, we can express $\mathrm{Gal}_{\mathrm{Sh}}(M)$ and $H_{\mathrm{KP}}(M)$ as the inverse limit of the groups $G_{\mathrm{Sh}}^{X},\,H_{X}$ for $X$ definable with definable complement where the directed system comes from projections of sorts and (contravariantly) from inclusions of type-definable sets in the same sort.} 
\end{corollary}
\begin{proof}
   First, $\mathrm{Gal}_{\mathrm{KP}}$ is given as the inverse limit of the local Kim-Pillay groups $G_{\mathrm{KP}}^{S}$ (Lemma \ref{lemmakimpillaygroupasinverselimit}).
   By Lemma \ref{lemmarestrictiondefequivcomesfromequiv} this is also the case for $\mathrm{Gal}_{\mathrm{Sh}}(M)$ and $G_{\mathrm{Sh}}^{S}$. Then it also follows that $H_{\mathrm{KP}}(M)$ is the inverse limit of the $H_{S}$: Taking inverse limits is left-exact, hence it suffices to establish surjectivity of the induced map $\phi:(\varprojlim G_{\mathrm{KP}}^{S})\rightarrow(\varprojlim G_{\mathrm{Sh}}^{S})$. But by compactness $\phi$ coincides with the natural map $\mathrm{Gal}_{\mathrm{KP}}(M)\rightarrow \mathrm{Gal}_{\mathrm{Sh}}(M)$ after identifying $\mathrm{Gal}_{\mathrm{KP}}(M)\cong\varprojlim G_{\mathrm{KP}}^{S}$ and $\mathrm{Gal}_{\mathrm{Sh}}\cong\varprojlim G_{\mathrm{Sh}}^{S}$. By Theorem \ref{theoremlocalconnectedcomponentshelahgroup} we can then conclude that $H_{\mathrm{KP}}(M)$ is the connected component of the identity in $\mathrm{Gal}_{\mathrm{KP}}(M)$.
\end{proof}
Finally, we observe that both the Shelah and the Kim-Pillay group do not depend on the choice of a monster model.
\begin{remark}\label{remarkKPgroupdoesnotdependonmonster}
   If $\mathcal{M}\prec\mathcal{M}^{\prime}$ for $\mathcal{M}$ and $\mathcal{M}^{\prime}$ both sufficiently saturated monster models of $T$, then the restriction map $\mathrm{Aut}(M^{\prime})\rightarrow \mathrm{Aut}(M)$ induces homeomorphisms $\mathrm{Gal}_{\mathrm{KP}}(M^{\prime})\rightarrow \mathrm{Gal}_{\mathrm{KP}}(M)$ and $\mathrm{Gal}_{\mathrm{Sh}}(M^{\prime})\rightarrow \mathrm{Gal}_{\mathrm{Sh}}(M)$.
\end{remark}
\begin{proof}
    We can use the exact same proof as that given in Corollary 42 in \cite{CategoryofModelsLascar}.
\end{proof}
\begin{convention}
In light of Remark \ref{remarkKPgroupdoesnotdependonmonster} we will from now on drop the $M$ in the notation of $\mathrm{Gal}_{\mathrm{KP}}$, $\mathrm{Gal}_{\mathrm{Sh}}$ and $H_{\mathrm{KP}}$ as well as in the corresponding notation for automorphism groups. 
\end{convention}

\begin{corollary}\label{corollarydescriptionconnectedcompasAutomgroup}
    The connected component $H_{\mathrm{KP}}$ of the Kim-Pillay group is isomorphic to the group $\mathrm{Aut}(\mathrm{acl}_{\mathrm{CL}}^{\mathrm{eq}}(\emptyset)/\mathrm{acl}_{\mathrm{DL}}^{\mathrm{eq}}(\emptyset))$.
\end{corollary}
\begin{proof}
    By Fact \ref{factKPgroupasautomgroup} we have $\mathrm{Gal}_{\mathrm{KP}}\cong \mathrm{Aut}(\mathrm{acl}_{\mathrm{CL}}^{\mathrm{eq}}(\emptyset))$. By Fact \ref{factShelahgroupasautomgroup} $\mathrm{Gal}_{\mathrm{Sh}}\cong \mathrm{Aut}(\mathrm{acl}_{\mathrm{DL}}(\emptyset))$ holds.
    Thus, we obtain from Corollary \ref{corollaryconnectedcomponentasinverselimit} that $H_{\mathrm{KP}}\cong \mathrm{Aut}(\mathrm{acl}_{\mathrm{CL}}^{\mathrm{eq}}(\emptyset)/\mathrm{acl}_{\mathrm{DL}}^{\mathrm{eq}}(\emptyset))$.
\end{proof}

\subsection{The Kim-Pillay group in simple theories}
In Section \ref{sectionautomgrouptorsors} we are in the situation of being given a set $\mathcal{D}$ of type-definable equivalence relations and we want to know whether it already suffices to be equivalent for all relations in $\mathcal{D}$ to have the same Kim-Pillay-strong type. We will now see how we can use the Independence Theorem as a criterion thereof.
\begin{definition}\label{definitionDstrongtype}
    Let $\mathcal{D}$ be a set of bounded $A$-type-definable equivalence relations that contain all definable equivalence relations (with finitely many classes). We say for two tuples $\bar{a},\bar{b}$ of the same length that they have the same $\mathcal{D}$-strong type, written $\bar{a}\equiv^{\mathcal{D}}\bar{b}$ if $\tau(\bar{a})=(\bar{b})$ for some $\tau\in \mathrm{Aut}(M/A)$ that preserves all $E$-classes for all $E\in\mathcal{D}$.
\end{definition}

\begin{theorem}\label{theorewhendstrongtypeequivlascarstrongtype}
The following are equivalent in a simple theory $T$ over some set $A$:
\begin{enumerate}
	\item The Kim-Pillay-strong type is the same as the $\mathcal{D}$-strong type over $A$.
	\item The Independence Theorem holds for $\mathcal{D}$-strong types over $A$. That is, whenever we have $\bar{a}_{0}\ind_{A} \bar{a}_{1},\, \bar{b}_{0}\ind_{A}\bar{a}_{0},\,\bar{b}_{1}\ind_{A}\bar{a}_{1}$ and $\bar{b}_{0}\equiv_{A}^{\mathcal{D}}\bar{b}_{1}$, then we can find some $\bar{b}\ind_{A}\bar{a}_{0}\bar{a}_{1}$ such that $\bar{b}\equiv_{A\bar{a}_{0}}^{\mathcal{D}}\bar{b}_{0}$ and $\bar{b}\equiv_{A\bar{a}_{1}}^{\mathcal{D}}\bar{b}_{1}$. 
\end{enumerate}
\end{theorem}
\begin{proof}
$(1.)\implies (2.)$ This follows directly from the Independence Theorem for Kim-Pillay-strong types.\\
$(2.)\implies (1.)$ Assume that $\equiv^{\mathrm{KP}}_{A}$ is strictly stronger than $\equiv^{\mathcal{D}}_{A}$. Then, we find $\bar{b}_{0},\bar{b}_{1}$ with $\bar{b}_{0}\equiv^{\mathcal{D}}_{A}\bar{b}_{1}$ but $\bar{b}_{0}\not\equiv^{\mathrm{KP}}_{A}\bar{b}_{1}$. Let $\bar{a}_{0}$ and $\bar{a}_{1}$ be such that $\bar{a}_{0}\ind_{A} \bar{a}_{1},\, \bar{b}_{0}\ind_{A}\bar{a}_{0},\,\bar{b}_{1}\ind_{A}\bar{a}_{1}$ and, moreover, $\bar{b}_{0}\equiv_{A}^{\mathrm{KP}}\bar{a}_{0}$ and $\bar{b}_{1}\equiv_{A}^{\mathrm{KP}}\bar{a}_{1}$ hold. (This is possible in simple theories by Proposition 5.5 in \cite{Kim199simpletheories} whose proof works also in our context.)
Now, there cannot exist any $\bar{b}$ with $\bar{b}\equiv_{A\bar{a}_{0}}\bar{b}_{0}$ and $\bar{b}\equiv_{A\bar{a}_{1}}\bar{b}_{1}$ because $\bar{b}\equiv_{A\bar{a}_{0}}\bar{b}_{0}$ now already implies that $\bar{b}\equiv_{A}^{\mathrm{KP}}\bar{b}_{0}$ and similarly we obtain $\bar{b}\equiv_{A}^{\mathrm{KP}}\bar{b}_{1}$, a contradiction.
\end{proof}

\begin{definition}
    We define $I$ to be the directed system consisting of all pairs $(S,E)$ where $S$ is a sort and $E$ an $A$-type-definable equivalence relation on $S$ with a bounded number of classes. The partial order is defined as $(S,E)\geq (S^{\prime},E^{\prime})$ if and only if $S^{\prime}$ projects onto $S$, say via the map $\pi:S^{\prime}\rightarrow S$, and $E^{\prime}$ refines (or equals) $\pi^{-1}(E)$ on $S^{\prime}$.
    Let $\mathcal{D}$ be a set of bounded $A$-type-definable equivalence relations that is closed under finite products.\footnote{I.e., if $E_{1},E_{2}\in\mathcal{D}$ and $E_{1}$ is defined on $S_{1}$, $E_{2}$ is defined on $S_{2}$, then we assume that the equivalence relation $E_{1}\times E_{2}$ on the product sort $S_{1}\times S_{2}$ is contained in $\mathcal{D}$.} We denote by $I_{\mathcal{D}}$ the directed sub-system consisting of all $(S,E)$ for $E\in\mathcal{D}$.
\end{definition}

\begin{notation}
    Below $\varprojlim$ will be the inverse limit over a system indexed by $I$, whereas $\varprojlim_{\mathcal{D}}$ will indicate that the system is indexed by $I_{\mathcal{D}}$.
\end{notation}
\begin{lemma}\label{lemmaKPgroupasinverselimitoverallspaces}
    For any $(S,E)\in I$ let $G^{S}_{E}$ (as in Definition \ref{definitionlocalKPgroup}) be equipped with the compact-open topology with respect to its action on $S/E$ (which carries the logic topology). Then $I$ induces a directed system of topological groups, and moreover $\varprojlim G^{S}_{E}$ and $\mathrm{Gal}_{\mathrm{KP}}$ are isomorphic as topological groups.
\end{lemma}
\begin{proof}
First note that $G^{S}_{E}$ is a compact Hausdorff topological group (using the exact same construction as for the $G_{\mathrm{KP}}^{X}$ in Definition \ref{definitionlocalKPgroup} and Lemma \ref{lemmalocalKPgroupastopologicalgroup}) that can equivalently be equipped with the compact-open topology for its action on $S/E$ (with the exact same proof as in Lemma \ref{lemmaequivalencetocompactopentopology}). It is then immediate that the $G^{S}_{E}$ form a directed system of topological groups and as pure groups $\varprojlim G^{S}_{E}\cong \mathrm{Gal}_{\mathrm{KP}}$. Clearly, for any sort $S$ we have $(S,E_{\mathrm{KP}}^{S})\in I$. Moreover, the projection $G^{S}_{\mathrm{KP}}\rightarrow G^{S}_{E}$ is a continuous map (using $G^{S}_{\mathrm{KP}}$ and, moreover, $\varprojlim G^{S}_{E}$ projects continuously to any $G^{S}_{\mathrm{KP}}$. Thus, $\varprojlim G^{S}_{E}$ and $\mathrm{Gal}_{\mathrm{KP}}$ are isomorphic as topological groups. 
\end{proof}

\begin{lemma}\label{lemmaKPgruppeinverserlimesmitPKstarkemsystem}
   Let $\mathcal{D}$ be a set of bounded $A$-type-definable equivalence relations that is closed under finite products. If the Kim-Pillay type is the same as the $\mathcal{D}$-strong type over $A$, it follows that $\mathrm{Gal}_{\mathrm{KP}}^{A}\cong \varprojlim_{\mathcal{D}} G^{S}_{E}$ as topological groups.
\end{lemma}
\begin{proof}
    The isomorphism as pure groups is immediate. Also, the natural isomorphism $\varprojlim G^{S}_{E}\rightarrow\varprojlim_{\mathcal{D}} G^{S}_{E}$ is continuous and both groups are inverse limits of compact Hausdorff groups and thus compact Hausdorff themselves. Consequently, the statement follows with Lemma \ref{lemmaKPgroupasinverselimitoverallspaces}.
\end{proof}

\section{Appendix II: Structure expansions on stably embedded sets}\label{sectionstructureexpansionsonstablyembeded}
The purpose of this section is to provide a relative quantifier elimination statement for structure expansions on stably embedded sets (Corollary \ref{corollaryabstractQEconition}). We work in a classical logic theory but the expansion allows for continuous logic predicates. Again, the statements are certainly well known to experts but to our knowledge the results we require are in this form only written down in the purely classical context. We use \cite{enrichedpredicate} and \cite{extdefbsetsanddependendpairs2} as our main sources and adapt the results to the continuous context.\\

In the following definition, we adapt the notion of stable embeddedness as presented in the Appendix of \cite{acfa} to the context of continuous logic. This was already presented in Chapter A.2.6 of the preprint \cite{chevalier2022piecewise}. 

\begin{definition}
Given some complete theory $T$ in some (possibly continuous) language $\mathcal{L}$. Let $\mathcal{M}\models T$ 
and $D\subseteq M$ be an $\emptyset$-definable set. Then we say that $D$ is \textit{stably embedded} if for any $\emptyset$-definable predicate  $P(\bar{z},\bar{x})$ and $\epsilon>0$ there is a $\emptyset$-definable predicate $\Tilde{P_{\epsilon}}(\bar{y},\bar{x})$ such that for any $\bar{b}\in M^{n}$ there is some $\bar{a}\in D^{m}$ such that \[\sup_{\bar{x}\in D^{k}}|P(\bar{b},\bar{x})-\Tilde{P_{\epsilon}}(\bar{a},\bar{x})|<\epsilon.\]
\end{definition}

\begin{remark}
From the definition it follows that if $D$ is stably embedded, then the restriction of every definable predicate to $D^{k}$ can be uniformly approximated (on $D^{k}$) by $D$-definable predicates and thus the restriction of the predicate itself is $D$-definable. 
\end{remark}



From now on let $\mathcal{M}$ be a saturated and strongly homogeneous model of size $\kappa$ for some regular cardinal $\kappa>|T|$. See Example 3.2 in \cite{halevikaplansaturatedmodelsformodeltheorists} for an explanation of why this does not pose any set-theoretic problem.
The following can be shown with the same proof as in Lemma 1 in the appendix of \cite{acfa}. See Lemma A.14 in \cite{chevalier2022piecewise} for an account thereof.


\begin{lemma}\label{lemmastableembeddednessequivalence}
    Let $\mathcal{M}$ be a saturated and strongly homogeneous model of size $\kappa$ for some regular cardinal $\kappa>|T|$ and let $D\subseteq M$ be an $\emptyset$-definable subset. The following assertions are equivalent:
    \begin{enumerate}
        \item $D$ is stably embedded.
        \item Every automorphism of $D$ lifts to an automorphism of $M$.
        \item For every $a\in M$ there is a small $D_{0}\subseteq D$ such that $\mathrm{tp}(a/D_{0})\vdash tp (a/D)$.
    \end{enumerate}
\end{lemma}

Now we turn to expansions of the structure on a stably embedded set. Our setting will be restricted to the following. We will work in a classic logic structure with a stably embedded set whose enrichment is allowed to use continuous logic predicates.
We closely follow the corresponding results in classical logic with \cite{enrichedpredicate} and \cite{extdefbsetsanddependendpairs2} as our main sources.
Throughout this section, we fix an $\mathcal{L}$-theory $T$ and a monster model $\mathcal{M}\models T$ with stably embedded definable set $D\subseteq M$.

\begin{definition}\label{definitioninducedstructure}
    Assume that $\mathcal{U}$ is an $\mathcal{L}^{\prime}$-structure on the set $D$ for some language $\mathcal{L}^{\prime}$. We say that $\mathcal{U}$ carries or has the $A$-induced structure from $\mathcal{M}$, if the following holds: Every $P(\bar{x})$ is $\emptyset$-definable in $\mathcal{U}$ if and only if it is given by the restriction (to $D^{|\bar{x}|}$) of an $A$-definable predicate in $M$. If $A=\emptyset$, we will later just say that $\mathcal{U}$ has the induced structure from $\mathcal{M}$.
\end{definition}
    
    

\begin{notation}
    From now on we work in the following situation: $\mathcal{D}_{\mathrm{ind}}$ is an $\mathcal{L}_{D-\mathrm{ind}}$-structure on $D$ that carries the $\emptyset$-induced structure from $\mathcal{M}$. We assume that $\mathcal{L}_{D-\mathrm{ind}}\subseteq \mathcal{L}$. $\mathcal{D}$ is an expansion of $\mathcal{D}_{\mathrm{ind}}$ in the language $\mathcal{L}_{D-\mathrm{exp}}$ where we assume $\mathcal{L}_{D-\mathrm{exp}}\backslash \mathcal{L}_{D-\mathrm{ind}}$ to be disjoint from the language $\mathcal{L}$. 
\end{notation}

\begin{convention}
    For the rest of the section we assume that $\mathcal{L}$ is a classical logic language, i.e., that predicates (including the metric) take their values in $\{0,1\}$ and as before we fix an $\mathcal{L}$-theory $T$ and a model $\mathcal{M}\models T$ with stably embedded set $D\subseteq M$. Further, we assume $\mathcal{L}_{D-\mathrm{exp}}$ to be a relational expansion (possibly in continuous logic) of $\mathcal{L}_{D-{\mathrm{ind}}}$. 
\end{convention}

    \begin{definition}\label{definitionextendedwholetheorystableembedded}
         Let $\mathcal{L}_{\mathcal{D}}=\mathcal{L}\cup\mathcal{L}_{D-\mathrm{exp}}$. Let $\mathcal{M}_{\mathcal{D}}$ be the $\mathcal{L}_{\mathcal{D}}$-structure obtained by enriching $D\subseteq M$ by $\mathcal{D}$ and extending every predicate from $\mathcal{L}_{D-\mathrm{exp}}$ by $0$ outside of $D$.
        We define the $\mathcal{L}_{\mathcal{D}}$-theory $T_{\mathcal{D}}$ to be $\mathrm{Th}(\mathcal{M}_{\mathcal{D}})$.

    \end{definition}

The following statement is Proposition 2.7 in \cite{enrichedpredicate}. The proof given there works in our context as well.
\begin{lemma}\label{lemmaenrichedstableembeddedstructure}
Let $\mathcal{M},\mathcal{N}\models T_{\mathcal{D}}$ be saturated models of the same cardinality. Then any $\mathcal{L}_{D-\mathrm{exp}}$-isomorphism from $D(M)$ to $D(N)$ extends to an $\mathcal{L}_{\mathcal{D}}$-isomorphism from $\mathcal{M}$ to $\mathcal{N}$.   
\end{lemma}
\begin{proof}
Let $f$ be an $\mathcal{L}_{D-\mathrm{exp}}$-isomorphism from $D(M)$ to $D(N)$ and $g$ be an $\mathcal{L}$-isomorphism from $M$ to $N$ which exists by saturation. Now $(g\restriction_{D(M)})^{-1}f$ is an automorphism of the $\mathcal{L}$-induced structure on $D(M)$ and by stable embeddedness of $D$ (in $T$) and Lemma \ref{lemmastableembeddednessequivalence} we can extend it to an $\mathcal{L}$-automorphism $h$ of $M$. Then it results that the map $gh:M\rightarrow N$ is an $\mathcal{L}$-isomorphism that extends $f$ and consequently $gh$ is an $\mathcal{L}_{\mathcal{D}}$-isomorphism. 
\end{proof}

\begin{corollary}\label{corollarypropertiesofenrihcedstableembeddedstructures}
Let $\mathcal{M}\models T_{\mathcal{D}}$. The theory $T_{\mathcal{D}}$ is complete, the definable set $\mathcal{D}$ is stably embedded in $T_{\mathcal{D}}$ and the $\mathcal{L}_{D}$-induced structure on $D(M)$ is $\mathcal{D}$.
\end{corollary}
\begin{proof}
We can assume $\mathcal{M}$ to be saturated. It follows directly from Lemma \ref{lemmaenrichedstableembeddedstructure} that $T_{\mathcal{D}}$ is complete. Moreover, Lemma \ref{lemmaenrichedstableembeddedstructure} yields the second point of the equivalence in Lemma \ref{lemmastableembeddednessequivalence} and thus $\mathcal{D}$ is stably embedded (in $T_{\mathcal{D}}$). Finally, let $P$ be an $\mathcal{L}_{\mathcal{D}}$-definable (over $\emptyset$) predicate. The restriction $P\restriction_{D^{n}}$ has to be invariant under every $\mathcal{L}_{D-\mathrm{exp}}$-automorphism $\sigma$ of $D(M)$ since otherwise $\sigma$ would lift by Lemma \ref{lemmaenrichedstableembeddedstructure} to an $\mathcal{L}_{\mathcal{D}}$-automorphism of $\mathcal{M}$ that does not leave $P$ invariant and thus contradicts the assumption of definability (of $P$). Thus, by saturation, we can conclude that $P\restriction_{D^{n}}$ was already $\mathcal{L}_{D-\mathrm{exp}}$-definable over $\emptyset$.
\end{proof}

\begin{remark}
    The results that we will obtain in the rest of this chapter will also work correspondingly for $\mathcal{L}$ a continuous language and $D$ a stably embedded definable set whose complement is also definable. In this context, we then need to assume elimination of CL-imaginaries in Corollary \ref{corollaryabstractQEconition} (and in Fact \ref{factdefbfunctionofcanocparamter}).
\end{remark}

\begin{definition}
    We denote by $\Tilde{\mathcal{L}}_{D-\mathrm{exp}}$ the language consisting of $\mathcal{L}_{D-\mathrm{exp}}$ together with the function symbols from $\mathcal{L}$. (That is, if $\mathcal{L}$ is a relational expansion of $\mathcal{L}_{\mathrm{ind}}$, then $\Tilde{\mathcal{L}}_{D-\mathrm{exp}}=\mathcal{L}_{D-\mathrm{exp}}$.) 
\end{definition}

The following Definition \ref{definitionDboundedformula} as well as Lemma \ref{lemmareductionqfformulastoDbdd} and Lemma \ref{lemmareductiontoDboundedformulas} will be an adaptation of Lemma 46 in \cite{extdefbsetsanddependendpairs2} to our context. Note that we do not only allow for a continuous logic expansion, but we will also slightly weaken one of the assumptions in \cite{extdefbsetsanddependendpairs2}. Whereas in \cite{extdefbsetsanddependendpairs2} the language $\mathcal{L}_{\mathcal{D}}$ is assumed to be relational, we will only require $\mathcal{L}_{D-\mathrm{exp}}$ to be a relational expansion of $\mathcal{L}_{D_{\mathrm{ind}}}$. Our notion of $D$-boundedness will a priori be slightly more restrictive (only allowing atomic relational $\mathcal{L}_{D-\mathrm{exp}}$-formulas) which will be of use later on. However, this is the reason why we will need Lemma \ref{lemmareductionqfformulastoDbdd}.

\begin{definition}\label{definitionDboundedformula}
    We say that an $\mathcal{L}_{\mathcal{D}}$-formula $\phi$ is $D$-bounded if $\phi$ is equivalent to a formula of the form
    \[\psi(\bar{x})=Q_{1}z_{1}\in D\dots Q_{n}z_{n}\in D\;h(\xi_{1}(\bar{x},\bar{z}),\dots,\xi_{m}(\bar{x},\bar{z}),\gamma_{1}(\bar{x},\bar{z}),\dots,\gamma_{l}(\bar{x},\bar{z}))\]
    where $h$ is some connective, the $Q_{i}$ are quantifiers, the $\xi_{i}(\bar{x},\bar{z})$ are atomic $\mathcal{L}$-formulas and the $\gamma_{i}(\bar{x},\bar{z})$ are atomic relational (i.e., no function symbol occurs) $\mathcal{L}_{D-\mathrm{exp}}$-formulas.
\end{definition}

\begin{lemma}\label{lemmareductionqfformulastoDbdd}
     We assume that $\mathcal{L}_{D-\mathrm{exp}}$ is a relational expansion of $\mathcal{L}_{D_{\mathrm{ind}}}$ and that $\mathcal{L}$ has quantifier elimination. Then, if two tuples $\bar{a}$ and $\bar{a}^{\prime}$ have the same $D$-bounded type, it follows that they have the same quantifier-free type in $\mathcal{L}_{\mathcal{D}}$, i.e.,  $\bar{a}\equiv^{\mathrm{qf}}_{\mathcal{L}_{\mathcal{D}}}\bar{a}^{\prime}$.
\end{lemma}
\begin{proof}
    Let $\phi(\bar{x})=h(\zeta_{1}(\bar{x}),\dots,\zeta_{k}(\bar{x}))$ be a quantifier-free $\mathcal{L}_{\mathcal{D}}$-formula with the $\zeta_{i}(\bar{x})$ being atomic. Clearly, it suffices to show that the $D$-bounded type of $\bar{a}$ determines the value of $\zeta_{i}(\bar{a})$ for any $1\leq i\leq k$. The $\zeta_{i}(\bar{x})$ are by definition either atomic $\mathcal{L}$-formulas or atomic $\Tilde{\mathcal{L}}_{D-\mathrm{exp}}$-formulas. We only have to take care of the latter case. Since $\mathcal{L}_{D-\mathrm{exp}}$ is a relational expansion of $\mathcal{L}_{D-\mathrm{ind}}$, we can thus assume that $\zeta_{i}$ is of the form $P(t_{i_{1}}(\bar{x}),\dots,t_{i_{l}}(\bar{x}))$ where $P$ is an $l$-ary predicate symbol from $\mathcal{L}_{D-\mathrm{exp}}\backslash\mathcal{L}_{D-\mathrm{ind}}$ and the $t_{i_{1}}(\bar{x}),\dots,t_{i_{l}}(\bar{x})$ are $\mathcal{L}$-terms.
    Let $\rho(\bar{x},\bar{z})$ be an $\mathcal{L}$-formula that takes value $0$, if $(t_{i_{1}}(\bar{x}),\dots,t_{i_{l}}(\bar{x}))=\bar{z}$ and $\bar{z}\in D^{n}$ hold and otherwise takes value $1$. The formulas $\phi_{1}(\bar{x})=\inf_{\bar{z}}\rho(\bar{x},\bar{z})$ and $\phi_{2}(\bar{x})=\inf_{\bar{z}}\min\{1,(\rho(\bar{x},\bar{z})+P(\bar{z}))\}$ are clearly $D$-bounded. Now, if $\phi_{1}(\bar{a})=0$, we have $\zeta_{i}(\bar{a})=\phi_{2}(\bar{a})$ and otherwise $\zeta_{i}(\bar{a})=0$.    
\end{proof}

\begin{lemma}\label{lemmareductiontoDboundedformulas}
    We assume that the $\mathcal{L}$-theory $T$ has quantifier elimination and that $\mathcal{L}_{D-\mathrm{exp}}$ is a relational expansion of $\mathcal{L}_{D-{\mathrm{ind}}}$. Then, if two tuples $\bar{a}$ and $\bar{a}^{\prime}$ have the same $D$-bounded type, it follows that $\bar{a}\equiv_{\mathcal{L}_{\mathcal{D}}}\bar{a}^{\prime}$.
\end{lemma}
\begin{proof}
    We assume $\mathcal{M},\mathcal{N}\models T_{\mathcal{D}}$ to be saturated models and show the statement by a back-and-forth argument. Let $\bar{a}\in M^{|\bar{a}|}$ and $\bar{a}^{\prime}\in N^{|\bar{a}|}$ have the same $D$-bounded type, and let $b\in M$ be arbitrary.
    First, consider the case $b\in D$. Let $p(x,\bar{a})$ denote the $D$-bounded type of $b$ over $\bar{a}$. Since for any formula $\phi(x,\bar{a})=0$ from $p(x,\bar{a})$ we have $\inf_{x\in D}\phi(x,\bar{a})=0$, which is a $D$-bounded formula, it follows from saturation that there is $b^{\prime}\in D(N)$ such that $\bar{a}b$ has the same $D$-bounded type as $\bar{a}^{\prime}b^{\prime}$.
    Now assume $b\notin D$. We can apply the above case to enlarge the tuple $\bar{a}$ (and similarly $\bar{a}^{\prime}$) with elements from $D$ such that $\mathrm{tp}_{\mathcal{L}}(\bar{a}b/D)$ is definable over $\bar{c}=\bar{a}\cap D$.\footnote{See Definition 8.1.4 in \cite{Tent_Ziegler_2012} for the definition of a definable type. By the third point of Lemma \ref{lemmastableembeddednessequivalence} $\mathrm{tp}_{\mathcal{L}}(\bar{a}b/D)$ is definable over a small set in $D$.}
    Take $b^{\prime}\in N$ such that $\bar{a}b\equiv_{\mathcal{L}}\bar{a}^{\prime}b^{\prime}$. Then $\mathrm{tp}_{\mathcal{L}}(\bar{a}^{\prime}b^{\prime}/D)$ is definable over $\bar{c}^{\prime}=\bar{a}^{\prime}\cap D$ as well. Let $\psi(\bar{x})$ be a $D$-bounded formula as in Definition \ref{definitionDboundedformula}, that is,
 \[\psi(\bar{x})=Q_{1}z_{1}\in D\dots Q_{n}z_{n}\in D\;h(\xi_{1}(\bar{x},\bar{z}),\dots,\xi_{m}(\bar{x},\bar{z}),\gamma_{1}(\bar{x},\bar{z}),\dots,\gamma_{l}(\bar{x},\bar{z}))\]
    where $h$ is some connective, the $\xi_{i}(\bar{x},\bar{z})$ are atomic $\mathcal{L}$-formulas and the $\gamma_{i}(\bar{x},\bar{z})$ are atomic relational $\mathcal{L}_{D-\mathrm{exp}}$-formulas. Again using stable embeddedness, we find $\mathcal{L}$-formulas  $\iota_{1}(\bar{y},\bar{z}),\dots,\iota_{m}(\bar{y},\bar{z})$ such that for every $1\leq i\leq m$ and every $\bar{z}\in D^{n}$ we have that $\iota_{i}(\bar{c},\bar{z})$ holds if and only if $\xi_{i}(\bar{a}b,\bar{z})$ holds. Here, $\iota_{i}$ is allowed to be constant (in $\bar{z}$). (If $\xi_{i}(\bar{x},\bar{z})$ does not depend on $\bar{z}$, then $\xi_{i}(\bar{a}b)=\xi_{i}(\bar{a}^{\prime}b^{\prime})$ already holds as $\bar{a}b\equiv_{\mathcal{L}}\bar{a}^{\prime}b^{\prime}$.) 
    Also, we have for all $1\leq j\leq l$ that either $\gamma_{j}(\bar{a}b,\bar{z})=\Tilde{\gamma}_{j}(\bar{a},\bar{z})$ for any $\bar{z}\in D^{n}$ and some atomic relational $\mathcal{L}_{D-\mathrm{exp}}$-formula $\Tilde{\gamma}_{i}$ (simply obtained by leaving out dummy-variables), or $\gamma_{j}(\bar{a}b,\bar{z})=0$ for any $\bar{z}\in D^{n}$ as $b\notin D$. We can thus simply replace the $\xi_{i}$ and $\gamma_{j}$ according to the above and we obtain some $D$-bounded formula $\phi(\bar{y})$ such that
    \[\psi(\bar{a}b)=\phi(\bar{a})=\phi(\bar{a}^{\prime})=\psi(\bar{a}^{\prime}b^{\prime}).\]
    This completes the proof.
\end{proof}
Recall that $\mathcal{L}$ is a classical logic language. In particular, elimination of imaginaries for the $\mathcal{L}$-theory $T$ will be meant in the sense of classical logic.

\begin{fact}\label{factdefbfunctionofcanocparamter}(See Lemma 2.4 (b) in \cite{enrichedpredicate}.)
    Assume that the $\mathcal{L}$-theory $T$ has uniform elimination of imaginaries. 
    For any $\mathcal{L}$-formula $\psi(\bar{x},\bar{y})$ there exist some $\mathcal{L}$-formula $\phi_{\psi}(\bar{x},\bar{z})$ and an $\mathcal{L}$-definable function $f_{\psi,\phi}(\bar{y})$ into $D(M)^{|\bar{z}|}$ such that for any $\bar{a}\in M^{|\bar{y}|}$ and tuple $\bar{x}$ from $D(M)$ we have $\phi_{\psi}(\bar{x},f_{\psi,\phi}(\bar{a}))$ holds if and only if $\psi(\bar{x},\bar{a})$ holds.
\end{fact}

\begin{definition}
    Assume that $T$ has uniform elimination of imaginaries. We denote by $\mathcal{L}_{\mathcal{D}}^\mathrm{qe}$ the language consisting of $\mathcal{L}_{\mathcal{D}}$ together with the functions $f_{\psi,\phi}$ of Fact \ref{factdefbfunctionofcanocparamter} for every $\mathcal{L}$-formula $\psi$. We can enrich any  $\mathcal{M}\models T_{\mathcal{D}}$ to an $\mathcal{L}_{\mathcal{D}}^\mathrm{qe}$-structure where we interpret the $f_{\psi,\phi}$ as in Fact \ref{factdefbfunctionofcanocparamter}. We will write $T_{D}$ again for the corresponding $\mathcal{L}_{\mathcal{D}}^\mathrm{qe}$-theory.
\end{definition}

\begin{corollary}\label{corollaryabstractQEconition}
    If the $\mathcal{L}$-theory $T$ has quantifier elimination and uniform elimination of imaginaries and the $\mathcal{L}_{D-\mathrm{exp}}
    $-theory $\mathrm{Th}(\mathcal{D})$ has quantifier elimination, then the theory $T_{\mathcal{D}}$ has quantifier elimination in $\mathcal{L}_{\mathcal{D}}^\mathrm{qe}$.
\end{corollary}
\begin{proof}
    First we can apply Lemma \ref{lemmareductiontoDboundedformulas} to obtain that every $\mathcal{L}_{\mathcal{D}}$-formula in $T_{\mathcal{D}}$ is approximable by continuous combinations of $D$-bounded ones. We consider a $D$-bounded formula $\psi(\bar{x})$ as in Definition \ref{definitionDboundedformula}, i.e.,\[\psi(\bar{x})=Q_{1}z_{1}\in D\cdots Q_{n}z_{n}\in D\;h(\xi_{1}(\bar{x},\bar{z}),\dots,\xi_{m}(\bar{x},\bar{z}),\gamma_{1}(\bar{x},\bar{z}),\dots,\gamma_{l}(\bar{x},\bar{z}))\]
    where $h$ is some connective, the $\xi_{i}(\bar{x},\bar{z})$ are atomic $\mathcal{L}$-formulas and the $\gamma_{j}(\bar{x},\bar{z})$ are atomic relational $\mathcal{L}_{D-\mathrm{exp}}$-formulas.
    We work in the notation of the proof of Lemma \ref{lemmareductiontoDboundedformulas}. For every $1\leq i\leq m$ let $f_{i}:=f_{\xi_{i},\iota_{i}}$ be the function corresponding to $\xi_{i}$ as given by Fact \ref{factdefbfunctionofcanocparamter}. We can then rewrite $\psi(\bar{x})$ as
\[\psi(\bar{x})=\bar{Q}\bar{z}\in D^{n}\;h(\iota_{1}(f_{1}(\bar{x}),\bar{z}),\dots,\iota_{m}(f_{m}(\bar{x}),\bar{z}),\gamma_{1}(\bar{x},\bar{z}),\dots,\gamma_{l}(\bar{x},\bar{z})).\]Next, for every $J\subseteq\{1,\dots,k\}$ define the $\{0,1\}$-valued quantifier-free $\mathcal{L}$-formula $\beta_{J}(\bar{x})$ via $\beta_{J}(\bar{x})=1$, if and only if $\{i\;|\;x_{i}\in D\}=J$. We obtain
\[\psi(\bar{x})=\sum_{J\subseteq\{1,\dots,k\}}\beta_{J}(\bar{x})\cdot\psi(\bar{x}).\]We define the formula $\phi(\bar{x},\bar{y})$ to be obtained from $\psi(\bar{x})$ by replacing the $f_{i}(\bar{x})$ by new variables, i.e.,
\[\phi(\bar{x},\bar{y})=\bar{Q}\bar{z}\in D^{n}\;h(\iota_{1}(\bar{y}_{1},\bar{z}),\dots,\iota_{m}(\bar{y}_{m},\bar{z}),\gamma_{1}(\bar{x},\bar{z}),\dots,\gamma_{l}(\bar{x},\bar{z})).\]
For every $J\subseteq\{1,\dots,k\}$ consider the formula $\phi_{J}(\bar{x},\bar{y})$ obtained from $\phi(\bar{x},\bar{y})$ by replacing $\gamma_{u}(\bar{x},\bar{z})$ by $0$ for every $1\leq u\leq l$ such that there is some $j\in\{1,\dots,k\}\backslash J$ such that $x_{j}$ occurs syntactically in $\gamma_{u}(\bar{x},\bar{z})$. It follows that
\[\phi(\bar{x},\bar{y})=\sum_{J\subseteq\{1,\dots,k\}}\beta_{J}(\bar{x})\cdot\phi(\bar{x},\bar{y})=\sum_{J\subseteq\{1,\dots,k\}}\beta_{J}(\bar{x})\cdot\phi_{J}(\bar{x},\bar{y}).\]
By quantifier elimination for $\mathrm{Th}(\mathcal{D})$ it follows that $\phi_{J}(\bar{x},\bar{y})$ is approximable by quantifier-free formulas $\theta_{J}(\bar{x},\bar{y})$ on $D^{|\bar{x}|+|\bar{y}|}$. (Here we included dummy-variables; it would suffice to work with $D^{|J|+|\bar{y}|}$.)
Let $\epsilon>0$ and choose for every $J\subseteq\{1,\dots,k\}$ some quantifier-free formula $\theta_{J,\epsilon}(\bar{x},\bar{y})$ such that $\sup_{(\bar{x},\bar{y})\in D^{|\bar{x}|+|\bar{y}|}}|\phi_{J}(\bar{x},\bar{y})-\theta_{J,\epsilon}(\bar{x},\bar{y})|\leq\epsilon$. 
Then, the formula 
\[\Tilde{\theta}_{\epsilon}(\bar{x}):=\sum_{J\subseteq\{1,\dots,k\}}\beta_{J}(\bar{x})\cdot\theta_{J,\epsilon}(\bar{x},f_{1}(\bar{x}),\dots,f_{m}(\bar{x}))\]is a quantifier-free $\mathcal{L}_{\mathcal{D}}^\mathrm{qe}$-formula and $\sup_{\bar{x}}|\psi(\bar{x})-\Tilde{\theta}_{\epsilon}(\bar{x})|\leq\epsilon$. Since $\epsilon>0$ was arbitrarily chosen, quantifier elimination for $T_{\mathcal{D}}$ in $\mathcal{L}_{\mathcal{D}}^\mathrm{qe}$ follows.
\end{proof}

\begin{lemma}\label{lemmaconservativeextonstablembset}
    We work with the assumptions of Corollary \ref{corollaryabstractQEconition}.
    Assume that $\mathcal{D}$ is a conservative expansion for definable sets, i.e., every set definable (with parameters) in  $\mathcal{D}$ was already definable (with parameters) in $D(M)$ as an $\mathcal{L}_{D-\mathrm{ind}}$-structure. It follows that any model $\mathcal{M}_{\mathcal{D}}\models T_{\mathcal{D}}$ is a conservative expansion for definable sets of the corresponding model $\mathcal{M}\models T$.
\end{lemma}
\begin{proof}
As the metric is discrete and $T_{D}$ has quantifier elimination in $\mathcal{L}_{\mathcal{D}}^\mathrm{qe}$ the definable sets in $\mathcal{M}_{\mathcal{D}}^{n}$ form a boolean algebra generated by sets definable by atomic $\mathcal{L}_{\mathcal{D}}^\mathrm{qe}$-formulas. Assuming that the statement of the lemma does not hold, we thhen find some set $X$ definable in $\mathcal{M}_{\mathcal{D}}$ by an atomic formula that was not definable in $\mathcal{M}$. Since $X$ was not definable in $\mathcal{M}$ and $\mathcal{D}$ was a relational expansion of $D(M)$, we can assume that the formula defining $X$ is given by an atomic formula of the form $P(t_{1}(\bar{x}),\dots, t_{n}(\bar{x}))$ for $P$ a predicate symbol from $\mathcal{L}_{D-\mathrm{exp}}\backslash\mathcal{L}_{D-\mathrm{ind}}$ and $t_{1}(\bar{x}),\dots, t_{n}(\bar{x})$ $\mathcal{L}$-terms.
The set $B:=\{(t_{1}(\bar{a}),\dots,t_{n}(\bar{a}))\;|\;\bar{a}\in \mathcal{M}_{\mathcal{D}}^{|\bar{x}|}\}\subseteq D^{n}$ is $\mathcal{L}$-definable with parameters in $D$ by stable embeddedness. But then there is $\epsilon>0$ such that $P(\bar{b})=0$, whenever $\bar{b}\in B_{X}:=\{(t_{1}(\bar{a}),\dots,t_{n}(\bar{a}))\;|\;\bar{a}\in X\}$ and $P(\bar{b})\geq\epsilon$ for $\bar{b}\in B\backslash B_{X}$. As $\mathcal{D}$ was a conservative expansion for definable sets, it follows that $B_{X}$ is already $\mathcal{L}$ definable and consequently $X$ is $\mathcal{L}$-definable, contradicting the assumption.
\end{proof}

\bibliography{bibpfplustimes}
\bibliographystyle{abbrv}
\end{document}